\newcommand\BB{{\mathbb B}}
\newcommand\CC{{\mathbb C}}
\newcommand\cA{{\cal A}}

\newcommand\cC{{\cal C}}
\newcommand\cD{{\cal D}}
\newcommand\cE{{\cal E}}
\newcommand\cF{{\cal F}}

\newcommand\cI{{\cal I}}

\newcommand\cK{{\cal K}}
\newcommand\cL{{\cal L}}

\newcommand\cO{{\cal O}}
  
\newcommand\cQ{{\cal Q}}

\newcommand\cS{{\cal S}}
\newcommand\cT{{\cal T}}
\newcommand\cU{{\cal U}}
\newcommand\cV{{\cal V}} 
\newcommand\cW{{\cal W}}
\newcommand\cX{{\cal X}} 
\newcommand\cY{{\cal Y}}
\newcommand\cZ{{\cal Z}}

\newcommand\DD{{\mathbb D}}

\newcommand\es{\emptyset}

\newcommand\GG{{\mathbb G}}

\newcommand\hra{\hookrightarrow}
\newcommand\la{\langle}
\newcommand\lag{\mathbb{LG}}

\newcommand\lagr{\mathbb{LG}(\bigwedge^3 V)}
\newcommand\lagrdual{\mathbb{LG}(\bigwedge^3 V^{\vee})}
\newcommand\lagre{\mathbb{LG}(\cE_W)}

\newcommand\LL{{\mathbb L}}
\newcommand\lra{\longrightarrow}
\newcommand\n{\noindent}
\newcommand\NN{{\mathbb N}}
\newcommand\ov{\overline}
\newcommand\PP{{\mathbb P}}

\newcommand\ra{\rangle}

\newcommand\sF{{\mathsf F}}

\newcommand\XX{{\mathbb X}}

\newcommand\wt{\widetilde}

\newcommand{\Gr}{\mathrm{Gr}}
\documentclass[10pt,a4paper]{article}
\usepackage{amsthm}
\usepackage{amsmath}
\usepackage{amssymb}
\usepackage{makeidx}         % permette di generare l'indice
\usepackage{graphicx}        % pacchetto grafico standard LaTeX 
\usepackage{multicol}        % usato per l'indice a due colonne
\usepackage[all,ps]{xy}
\usepackage{layout}
\usepackage{booktabs}

\theoremstyle{plain}

\newtheorem{thm}{Theorem}[section]
\newtheorem{clm}[thm]{Claim}

\newtheorem{crl}[thm]{Corollary}

\newtheorem{lmm}[thm]{Lemma}
\newtheorem{prp}[thm]{Proposition}
\newtheorem{prp-dfn}[thm]{Proposition-Definition}

\theoremstyle{definition}

\newtheorem{dfn}[thm]{Definition}

\theoremstyle{remark}

\newtheorem{rmk}[thm]{Remark}

\DeclareMathOperator{\Ann}{Ann}

\DeclareMathOperator{\chord}{chord}
\DeclareMathOperator{\cod}{cod}

\DeclareMathOperator{\coker}{coker}
\DeclareMathOperator{\cork}{cork}

\DeclareMathOperator{\Hom}{Hom}
\DeclareMathOperator{\Id}{Id}
\DeclareMathOperator{\im}{im}

\DeclareMathOperator{\mult}{mult}

\DeclareMathOperator{\Pic}{Pic}

\DeclareMathOperator{\rk}{rk}

\DeclareMathOperator{\sing}{sing}

\DeclareMathOperator{\supp}{supp}
\DeclareMathOperator{\Sym}{Sym}
\DeclareMathOperator{\vol}{vol}

\usepackage{ifthen}
\newcommand{\cit}[1]{{\rm \textbf{#1}}}
\newcommand{\Ref}[2]{\cit{%
\ifthenelse{\equal{#1}{thm}}{Theorem}{}%
\ifthenelse{\equal{#1}{ass}}{Assumption}{}%
\ifthenelse{\equal{#1}{chp}}{Chapter}{}%
\ifthenelse{\equal{#1}{prp}}{Proposition}{}%
\ifthenelse{\equal{#1}{lmm}}{Lemma}{}%
\ifthenelse{\equal{#1}{crl}}{Corollary}{}%
\ifthenelse{\equal{#1}{dfn}}{Definition}{}%
\ifthenelse{\equal{#1}{expl}}{Example}{}%
\ifthenelse{\equal{#1}{hyp}}{Hypothesis}{}%
\ifthenelse{\equal{#1}{rmk}}{Remark}{}%
\ifthenelse{\equal{#1}{clm}}{Claim}{}%
\ifthenelse{\equal{#1}{exe}}{Exercise}{}%
\ifthenelse{\equal{#1}{sec}}{Section}{}%
\ifthenelse{\equal{#1}{subsec}}{Subsection}{}%
\ifthenelse{\equal{#1}{univ}}{Universal Property}{}%
\ifthenelse{\equal{#1}{trm}}{Terminology}{}%
\ifthenelse{\equal{#1}{tbl}}{Table}{}%
\  \ref{#1:#2}%
}}

\usepackage{multirow}

 \makeindex
 
\begin{document}
 \title{EPW-sextics: taxonomy}
 \author{Kieran G. O'Grady\thanks{Supported by
 PRIN 2007}\\\\
\lq\lq Sapienza\rq\rq Universit\`a di Roma}
\date{November 24  2010}
 \maketitle
 \tableofcontents
 \section{Introduction}\label{sec:prologo}
 \setcounter{equation}{0}
EPW-sextics are special sextic hypersurfaces in $\PP^5$ which come equipped with a double cover ramified over their singular locus (generically a smooth surface). They were introduced by Eisenbud, Popescu and Walter~\cite{epw} in order to give examples of a \lq\lq quadratic sheaf\rq\rq (on a hypersurface)  which does not admit a symmetric resolution. We proved~\cite{og2}  that if the EPW-sextic is generic then the double cover is a hyperk\"ahler (HK) $4$-fold deformation deformation equivalent to the Hilbert square of a $K3$, moreover  the family of (smooth) double EPW-sextics is a locally complete family of projective HK's. We recall that three other locally complete families of projective HK's of dimension greater than $2$ are known, those introduced by Beauville and Donagi~\cite{beaudon},  Debarre and Voisin~\cite{debvoi}, Iliev and Ranestad~\cite{iliran1,iliran2}; in all of the above examples the HK manifolds are deformations of  the Hilbert square of a $K3$ and they are distinguished by the value of the Beauville-Bogomolov form on the polarization class (it equals $2$ in the case of double EPW-sextics and $6$, $22$ and $38$ in the other cases). EPW-sextics are defined as follows. Let $V$ be a $6$-dimensional complex vector space - this notation will be in force 
throughout the paper. We choose a volume-form on $V$ 
\begin{equation}\label{volumone}
\vol\colon\bigwedge^6 V\overset{\sim}{\lra}\CC
\end{equation}
 and we equip $\bigwedge^3 V$ with the symplectic form
\begin{equation}
  (\alpha,\beta)_V:=\vol(\alpha\wedge\beta).
\end{equation}
Let $\lagr$ be the symplectic Grassmannian parametrizing Lagrangian subspaces of $\bigwedge^3 V$ - of course $\lagr$ does not depend on the choice of volume-form. 
Given a non-zero $v\in V$ we let 
\begin{equation}
F_v:=\{\alpha\in\bigwedge^3 V\mid v\wedge\alpha=0\}\index{$F_v$}
\end{equation}
be the sub-space of $\bigwedge^3 V$ consisting of multiples of $v$. Notice that $(,)_V$ is zero on $F_v$ and   $\dim(F_v)=10$; thus $F_v\in \lagr$. Let 
\begin{equation}\label{eccoeffe}
F\subset\bigwedge^3 V\otimes\cO_{\PP(V)}
\end{equation}
be  the sub-vector-bundle with fiber $F_v$ over $[v]\in\PP(V)$. A straightforward computation gives that
\begin{equation}\label{ciunoeffe}
\det F\cong\cO_{\PP(V)}(-6).
\end{equation}
Given $A\in\lagr$ we let
\begin{equation}
Y_A=\{[v]\in\PP(V)\mid F_v\cap A\not=\{0\}\}.
\end{equation}
Thus $Y_A$ is the degeneracy locus of the map
\begin{equation}\label{diecidieci}
F\overset{\lambda_A}{\lra}(\bigwedge^3 V/A)\otimes\cO_{\PP(V)}
\end{equation}
where $\lambda_A$ is given by Inclusion~\eqref{eccoeffe} 
followed by the quotient map 
\begin{equation*}
\bigwedge^3 V\otimes\cO_{\PP(V)}\to 
(\bigwedge^3 V/A)\otimes\cO_{\PP(V)}.
\end{equation*}
 Since the vector-bundles appearing in~\eqref{diecidieci} have equal rank  the determinat of $\lambda_A$ makes sense and of course $Y_A=V(\det\lambda_A)$; this formula shows that $Y_A$ has a natural structure of closed subscheme of $\PP(V)$.  By~\eqref{ciunoeffe} we have  $\det\lambda_A\in H^0(\cO_{\PP(V)}(6))$ and hence  $Y_A$ is either a sextic hypersurface or $\PP(V)$. An {\it EPW-sextic}\index{EPW-sextic} is a sextic hypersurface in $\PP^5$ which is projectively equivalent to $Y_A$ for some $A\in\lagr$.  One verifies readily that EPW-sextics exist; in fact 
given $[v]\in\PP(V)$ there exists $A\in\lagr$ such that $A\cap F_v=\{0\}$ and hence $[v]\notin Y_A$. (On the other hand there do exist $A\in\lagr$ such that $Y_A=\PP(V)$ e.g.~$A=F_w$ for $[w]\in\PP(V)$.)
 Let
\begin{eqnarray}
\Sigma:= & \{A\in\lagr\mid \text{$\exists W\in{\mathbb G}r(3,V)$ s.~t.~$\bigwedge^3 W\subset A$}\},\\
\Delta:=  &  \{A\in\lagr\mid \text{$\exists[v]\in\PP(V)$ s.~t.~$\dim(A\cap F_v)\ge 3$}\}\,.
\end{eqnarray}
(We will denote $\Sigma$ by $\Sigma(V)$  whenever we will need to keep track of $V$, and similarly for $\Delta$).
Then $\Sigma$ and $\Delta$ are closed subsets of $\lagr$; a straightforward computation  shows that $\Sigma$ and $\Delta$ are irreducible of codimension $1$ - see~\Ref{sec}{singplan} for the case of $\Sigma$.
Let
\begin{equation}\label{eccozero}
\lagr^0:=  \lagr\setminus\Sigma\setminus\Delta\,.
\end{equation}
Thus $\lagr^0$ is open dense in $\lagr$. 
In~\cite{og2} we proved  the following results. If $A\in\lagr^0$ then $Y_A\not=\PP(V)$ and there exists a finite degree-$2$ map $f_A
\colon X_A\to Y_A$   unramified over the smooth locus of $Y_A$ with $X_A$  a HK $4$-fold deformation equivalent to $(K3)^{[2]}$.
For $A\in\lagr^0$ let $h_A:=c_1(f_A^{*}\cO_{Y_A}(1))$. We proved that the family of polarized $4$-folds 
 \begin{equation*}
 \{ (X_A,h_A)\}_{A\in \lagr^0}
\end{equation*}
 is  locally complete.  Let us compare the family of double EPW-sextics and the family of HK $4$-folds introduced by Beauville and Donagi~\cite{beaudon}. Donagi and Beauville consider a cubic $4$-fold $Z\subset\PP^5$ and the family $F(Z)$ parametrizing lines in $Z$; they proved that if $Z$ is smooth then $F(Z)$ is a HK $4$-fold deformation equivalent to the Hilbert square of a $K3$. Moreover they showed that the primitive weight-$2$ integral Hodge structure of $F(Z)$ is isomorphic to the integral primitive weight-$4$ Hodge structure of $Z$ (after a Tate twist) and that the isomorphism takes the Beauville-Bogomolov quadratic form on $H^2(F(Z))_{pr}$ to the opposite of the intersection from on $H^4(Z)_{pr}$. Thus the period map for the family $\{F(Z)\}$ may be studied via the period map for cubic $4$-folds. Periods of cubic $4$-folds were first studied by Voisin~\cite{claire} who proved the Global Torelli Theorem. More recently  Laza~\cite{laza1,laza2} and Looijenga~\cite{eddyloo} proved various results, in particular they gave a complete description of the periods of  smooth cubics. 
 
 This  is the first in a series of papers on moduli and periods of double EPW-sextics.  In order to present the  results of the present paper we introduce the following notation: given $A\in\lagr$ 
 we let
\begin{equation}\label{eccoteta}
\Theta_A:=\{W\in{\mathbb Gr}(3,V)\mid \bigwedge^3 W\subset A\}.
\end{equation}
Our main result  is a classification of those $A$ such that $\Theta_A$ has strictly positive dimension (in particular $A\in\Sigma$). Why are we concerned with such $A$ ? The period map $\lagr^0\to\DD$ extends to a rational map  $\lagr\dashrightarrow\DD^{BB}$ where $\DD^{BB}$ is the  Baily-Borel~\cite{bb} compactification of $\DD$: if $\dim\Theta_A>0$ then either the period map is not regular at $A$ or it goes to the boundary of $\DD^{BB}$. Moreover many of the non-stable (in the sense of GIT)   $A\in\lagr$ have positive-dimensional $\Theta_A$ -  the results of the present work will shed light on the description of the GIT-stable  points in $\lagr$ that will be appear in a forthcoming paper. 
Let us look at the analogous case of cubic $4$-folds. We claim that the prime divisor $D\subset |\cO_{\PP^5}(3)|$ parametrizing singular cubics is analogous to $\Sigma$. As is well-known $F(Z)$ is smooth if and only if $Z\in(|\cO_{\PP^5}(3)|\setminus D)$ and   if $Z$ is a singular cubic $4$-fold then $\sing F(Z)$ has dimension at least $2$ (generically it is a $K3$ of degree $6$). Moreover  the period map extends across the generic $Z\in D$ but it does not lift to the relevant classifying space: in order to lift it 
one needs  first to take a (local) double  cover ramified over $D$.  In the case of interest to us similar  results hold. Let $A\in(\lagr\setminus\Sigma)$; then $X_A$ is either smooth (if $A\in\lagr^0$) or the contraction of a finite union of (disjoint) copies  of $\PP^2$ in a $4$-fold $X^{\epsilon}_A$ with a holomorphic symplectic form\footnote{If $A$ is generic in $(\Delta\setminus\Sigma)$ then  $X^{\epsilon}_A$ is projective but it might not be K\"ahler for a particular $A$.}. On the other hand if $A\in\Sigma$ (and $Y_A\not=\PP(V)$) then $X_A$ has singular locus of dimension at least $2$ (generically a $K3$ of degree $2$).  What about periods ? The period map extends regularly on $(\Delta\setminus\Sigma)$ and it lifts to the classifying space. On the other hand let  $A\in\Sigma$ be  generic: the period map  extends across $A$ but 
  in order to lift it  to the relevant classifying space 
one needs  first to take a (local) double cover ramified over $\Sigma$. Thus one might view the $A$ such that $\dim\Theta_A>0$
 as analogues of cubic $4$-folds whose singular locus is of strictly positive dimension - we notice that such cubics play a prominent r\^ole in Laza's papers~\cite{laza1,laza2}. The  following simple remark is very useful when analyzing cubics with positive dimensional singular locus: if $Z\subset\PP^5$ is a cubic $4$-fold and $p,q\in Z$ are distinct points then the line joining $p$ and $q$ is contained in $Z$. The elementary remark below might be considered as an analogue in our context.
\begin{rmk}\label{rmk:adueadueinc}
Let $\Theta\subset{\mathbb Gr}(3,V)$. The following statements are equivalent:
\begin{itemize}
\item[(1)]
 $\dim(W_1\cap W_2)>0$ for any $W_1,W_2\in\Theta$.
\item[(2)]
The symplectic form $(,)_V$ vanishes on the subspace $\la\la\Theta\ra\ra\subset \bigwedge^3 V$ spanned by $\bigwedge^3 W$ for $W\in\Theta$.
\end{itemize}
In particular if $A\in\lagr$ then  $\PP(W_1)\cap \PP(W_2)\not=\es$ for any $W_1,W_2\in\Theta_A$.
\end{rmk}
Morin~\cite{morin} classified maximal families parametrizing  pairwise incident planes in $\PP^5$. Modulo  projectivities there are $6$ such families: $3$  elementary (or Schubert) families, namely planes containing a fixed point, planes contained in a hyperplane and planes whose intersection with a fixed plane has dimension at least $1$, and $3$ non-elementary families, namely planes contained in a smooth quadric hypersurface, planes tangent to a Veronese surface   and planes intersecting a Verones surface in a conic,
   see~\Ref{thm}{teomorin}. The non-elementary families give rise to EPW-sextics which are a triple quadric (the first case) and a double discriminant cubic (the second and third case); they are in the indeterminacy locus of the period map and they correspond to double EPW-sextics approaching HK $4$-folds with a  (pseudo)polarization defining a map which is no longer $2$-to-$1$ onto its image - see~\cite{ferretti} for a discussion of the first case. Building on Morin's theorem we will classify the possible positive-dimensional  irreducible components of $\Theta_A$. 

The paper is organized as follows. In the first section we will prove some basic results on EPW-sextics. In particular we will show that $\Theta_A$ determines how pathological   $Y_A$ might be  - for example $Y_A=\PP(V)$ if and only if the planes in $\Theta_A$ sweep out all of $\PP(V)$. We will also show how to produce a triple smooth quadric, a double discriminant cubic and the union of $6$ independent hyperplanes as EPW-sectics. In the last subsection we will show that EPW-sextics have a \lq\lq classical\rq\rq description as discriminant loci of certain linear systems of quadrics in $\PP^9$ (see~\cite{ilimani} for related results). The second section begins with some dimension counts for natural subsets of $\Sigma$ and standard infinitesimal computations. The main body  of that section is devoted to a classification of the elements of
\begin{equation}\label{sigdiesis}
\Sigma_{\infty}:=  \{A\in\lagr\mid \dim\Theta_A>0\}.
\end{equation}
In particular we will describe the irreducible components of  $\Sigma_{\infty}$ and we will compute their dimension.  
Going back to the analogy with the family of cubic $4$-folds: the family of double EPW-sextics has a more elaborate geometry, in fact there are $12$  irreducible components of $\Sigma_{\infty}$ while the set of cubic $4$-folds with positive dimensional singular locus has $5$ irreducible components - see Theorem~6.1 of~\cite{laza1}. 
\vskip 3mm
\n
{\bf Notation and conventions:} 
Let $W$ be a finite-dimensional complex vector-space. The span of a subset $S\subset W$ is denoted by $\la S\ra$. Let $S\subset\bigwedge^q W$.  The {\it support of $S$} is the smallest subspace $U\subset W$ such that $S\subset\im(\bigwedge^{q}U\lra \bigwedge^q W)$: we denote it by $\supp(S)$, if $S=\{\alpha\}$ is a singleton we let $\supp(\alpha)=\supp(\{\alpha\})$ \index{$\supp(S)$, $\supp\alpha$} (thus if $q=1$ we have $\supp(\alpha)=\la\alpha\ra$).  We define the support of a set of symmetric tensors analogously. If  $\alpha\in\bigwedge^q W$ or $\alpha\in \Sym^d W$ the {\it rank of $\alpha$} is the  dimension of $\supp(\alpha)$. 
An element of $\Sym^2 W^{\vee}$ may be viewed either  
as a symmetric map or as a quadratic form: we will denote the former by $\wt{q},\wt{r},\ldots$ and the latter by $q,r,\ldots$ respectively.
\vskip 2mm
\noindent
Let $U$ be a vector space. The  
{\it wedge subspace} of $\bigwedge^d U$ associated to a collection of subspaces $U_1,\ldots,U_{\ell}\subset U$ and a partition $i_1+\cdots+ i_{\ell}=d$ is defined as the span
\begin{equation}
(\bigwedge^{i_1} U_1)\wedge\cdots \wedge(\bigwedge^{i_{\ell}} U_{\ell})
:=\la\alpha_1\wedge\cdots\wedge\alpha_{\ell}\mid \alpha_s\in\bigwedge^{i_s}U_s\ra
\index{$(\bigwedge^{i_1} U_1)\wedge\cdots \wedge(\bigwedge^{i_{\ell}} U_{\ell})$}
\index{wedge-subspace}
\end{equation}
\vskip 2mm
\noindent
Let $W$ be a finite-dimensional complex vector-space. We will adhere to pre-Grothendieck conventions: $\PP(W)$ is the set of $1$-dimensional vector subspaces of $W$. Given  a non-zero $w\in W$ 
we will denote the span of $w$ by $[w]$ rather than $\la w\ra$; this agrees with  standard notation. 
Suppose that $T\subset\PP(W)$. Then $\la T\ra\subset\PP(W)$\index{$\la T\ra$} is the {\it projective span of  $T$}\index{projective span} i.e.~the intersection of all linear subspaces of $\PP(W)$ containing $T$ while
 $\la\la T\ra\ra\subset W$\index{$\la\la \cdot\ra\ra$} is the {\it vector-space span of  $T$}\index{vector-space span}   i.e.~the span of all $w\in(W\setminus\{0\})$ such that $[w]\in T$. 
\vskip 2mm
\noindent
Schemes  are defined over $\CC$, the topology is the Zariski topology unless we state the contrary. Let $W$ be finite-dimensional complex vector-space: $\cO_{\PP(W)}(1)$ is the line-bundle on $\PP(W)$ with fiber $L^{\vee}$ on the point $L\in\PP(W)$. Let $F\in \Sym^d W^{\vee}$: we let  $V(F)\subset\PP(W)$ be the subscheme defined by vanishing of $F$. If $E\to X$ is a vector-bundle we denote by $\PP(E)$ the projective fiber-bundle with fiber $\PP(E(x))$ over $x$ and we define $\cO_{\PP(W)}(1)$ accordingly.
If $Y$ is a subscheme of $X$ we let $Bl_Y X\lra X$ be the blow-up of $Y$.
\vskip 3mm
\n
{\bf Acknowledgments:} It is a pleasure to thank 
Claudio Procesi for conversations relating to  this work  and Ciro Ciliberto for pointing out  Morin's Theorem~\cite{morin} and for help with the proof of~\Ref{prp}{supgrass}. 
 \section{EPW-sextics}\label{sec:sistole}
 \setcounter{equation}{0}
\subsection{Symplectic Grassmannians}
\setcounter{equation}{0}
Let ${\bf H}$  be a complex vector-space of dimension $2n$ equipped with a symplectic form $(,)_{\bf H}$. We let 
$\mathbb{LG}({\bf H})\subset {\mathbb G}r(n,{\bf H})$\index{$\mathbb{LG}({\bf H})$} be the symplectic
Grassmannian parametrizing Lagrangian
subspaces of ${\bf H}$. We will go through some well-known results regarding $\mathbb{LG}({\bf H})$. 
 Let $A\in\LL\GG({\bf H})$: the symplectic form gives an isomorphism  
\begin{equation}\label{abdual}
\begin{matrix}
{\bf H}/A & \overset{\sim}{\lra} & A^{\vee}\\
\ov{x} & \mapsto & (a\mapsto (a,x)_{\bf H})
\end{matrix}
\end{equation}
and hence
we have a  canonical inclusion
\begin{equation}\label{dentrom}
T_{A}\LL\GG({\bf H})\subset T_A\Gr(n,{\bf H})=\Hom(A,{\bf H}/A)=A^{\vee}\otimes A^{\vee}
\end{equation}
 Let $\mathbb{LG}({\bf H})\hra \PP(\bigwedge^{n}{\bf H})$ be the Pl\"ucker embedding: the pull-back of the ample generator of $\Pic(\PP(\bigwedge^ {n}{\bf H}))$ is 
the {\it Pl\"ucker line-bundle} on $\LL\GG({\bf H})$. The following result is well-known; one reason for providing a proof is to introduce notation that will be used throughout the paper. 
\begin{prp}\label{prp:globcalcul}
Keep notation and hypotheses as above. 
\begin{itemize}
\item[(1)]
$\mathbb{LG}({\bf H})$ is  smooth, irreducible and Inclusion~\eqref{dentrom} identifies $T_{A}\LL\GG({\bf H})$ with $\Sym^2 A^{\vee}$.  
\item[(2)]
The Picard group of $\LL\GG({\bf H})$ is generated by the class of the Pl\"ucker line-bundle. 
\end{itemize}
\end{prp}
\begin{proof} 
The symplectic group $Sp({\bf H})$ acts transitively on $\mathbb{LG}({\bf H})$ and hence $\mathbb{LG}({\bf H})$ is smooth. Given $B\in \mathbb{LG}({\bf H})$ we let
 \begin{equation}\label{eccoua}
U_B:=\{C\in \mathbb{LG}({\bf H}) \mid B\cap C=\{0\}\}\,.
\end{equation}
Clearly $U_B$ is  open in $\mathbb{LG}({\bf H})$.
One defines  a (non canonical) isomorphism of varieties
\begin{equation}\label{identifico}
\Sym^2 B\lra U_B
\end{equation}
as follows.
Choose $C\in U_B$. The direct-sum decomposition 
${\bf H}=C\oplus B$ defines an isomorphism $C\overset{\sim}{\lra}{\bf H}/B$; composing with the isomorphism ${\bf H}/B \overset{\sim}{\lra} B^{\vee}$ (see~\eqref{abdual}) we get an isomorphism
$\iota\colon C  \overset{\sim}{\lra}  B^{\vee}$. Let 
$\wt{q}\in \Sym^2 B$ and view $\wt{q}$ as a symmetric map $B^{\vee}\to B$; the graph $\Gamma_{\wt{q}}$ of $\wt{q}$ lies in 
$B^{\vee}\oplus B$ and hence
\begin{equation}
(\iota,\Id_B)^{-1}\Gamma_{\wt{q}}\subset C\oplus B={\bf H}\,.
\end{equation}
Moreover $(\iota,\Id_B)^{-1}\Gamma_{\wt{q}}$ is Lagrangian because $\wt{q}$ is symmetric and 
 it belongs to $U_B$ because $\Gamma_{\wt{q}}$ is a graph. We define~\eqref{identifico} by sending $\wt{q}$ to $(\iota,\Id_B)^{-1}\Gamma_{\wt{q}}$.   Now choose $B$ transversal to $A$. Then $A\in U_B$ and hence we may choose $C=A$. We have defined an isomorphism $\iota\colon A  \overset{\sim}{\lra}  B^{\vee}$  and hence~\eqref{identifico} gives an isomorphism $\Sym^2 A^{\vee}\lra U_B$: the differential at $0$ 
 equals~\eqref{dentrom} and this proves that~\eqref{dentrom} identifies $T_{A}\LL\GG({\bf H})$ with
 $\Sym^2 A^{\vee}$. Irreducibility of  $\mathbb{LG}({\bf H})$ follows from the following two facts:
first the open sets $U_B$ for $B$ varying in $\mathbb{LG}({\bf H})$ form a covering of $\mathbb{LG}({\bf H})$
 and secondly $U_{B}\cap U_{B'}$ is non-empty for arbitrary $B,B'\in\mathbb{LG}({\bf H})$.  
 Let's prove Item~(2). Given $A\in \mathbb{LG}({\bf H})$ we let
 \begin{equation}\label{divpluck}
D_A:=  \{B\in \mathbb{LG}({\bf H}) \mid A\cap B\not=\{0\}\}=
 (\mathbb{LG}({\bf H})\setminus U_A)\,.
\end{equation}
One checks easily that $D_A$ is of pure codimension $1$ in $\mathbb{LG}({\bf H})$ and hence it may be viewed as an effective divisor: in fact it belongs to the Pl\"ucker linear system.
We have an exact sequence of Chow groups (see Proposition~(1.8) of~\cite{fulton})
\begin{equation}
CH^0(D_A)\lra CH^1(\mathbb{LG}({\bf H}))\lra CH^1(U_A)\lra 0\,.
\end{equation}
Isomorphism~\eqref{identifico} gives that $CH^1(U_A)=0$ and hence $CH^1(\mathbb{LG}({\bf H}))$ is generated by the classes of irreducible components of $D_A$. Since $D_A$ is  irreducible and it it belongs to the Pl\"ucker linear system  we get Item~(2).
\end{proof}
\begin{crl}\label{crl:picuno}
Let $V$ be a $6$-dimensional complex vector-space. Then  $\lagr$ is  irreducible, smooth of dimension $55$ and  its Picard group is generated by the class of the Pl\"ucker line-bundle. 
\end{crl}
\subsection{Degeneracy loci attached to $A\in\lagr$}\label{subsec:luogodeg}
\index{$Y_A[k]$, $Y_A(k)$}
\setcounter{equation}{0}
Let $A\in\lagr$. We let
\begin{equation}
Y_A[k]=\{[v]\in\PP(V)\mid \dim(A\cap F_v)\ge k\}\,.\index{$Y_A[k]$}
\end{equation}
Thus $Y_A[0]=\PP(V)$ and $Y_A[1]=Y_A$. We will show that  $Y_A[k]$ has a natural  structure of closed sub-scheme of $\PP(V)$. First we associate to $B\in\lagr$  the open subset
 $\cU_B\subset\PP(V)$  defined by
\begin{equation}\label{calub}
\cU_B:=\{[v]\in\PP(V)\mid  F_v\cap B=\{0\}\}\,.
\end{equation}
(In other words $\cU_B$ is the intersection of $U_B$ and $\PP(V)$ embedded in $\lagr$ by the map  $[v]\mapsto F_v$.)  Choose $B$ transversal to $A$: we will write $Y_A[k]\cap\cU_B$ as the $k$-th degeneracy locus of a symmetric map of vector-bundles. 
We have a direct-sum decomposition $\bigwedge^3 V=A\oplus B$ and for  $[v]\in\cU_B$ the Lagrangian subspace    $F_{v}$ is transversal to $B$; thus $F_v$ is the graph of a symmetric map 
\begin{equation}\label{mappagrafo}
\tau_A^B([v])\colon A\to B\cong A^{\vee}\,.\index{$\wt{q}(v)$}
\end{equation}
(The symplectic form $(,)_V$ together with the decomposition $\bigwedge^3 V=A\oplus B$ provides us with  an isomorphism $B\cong A^{\vee}$ - see the proof of~\Ref{prp}{globcalcul}.)
Since $\tau_A^B\colon\cU_B\to \Sym^2 A^{\vee}$ is a regular map we may  define a closed subscheme  $Y_A^B[k]\subset\cU_B$  by setting
\begin{equation}
Y_A^B[k]:=V(\bigwedge^{(11-k)}\tau_A^B)\,.
\end{equation}
The support of $Y_A^B[k]$ is equal to $Y_A[k]\cap\cU_B$. If $B'\subset\bigwedge^3 V$ is another Lagrangian  subspace transversal to $A$ then the restrictions of $Y_A^B[k]$ and $Y_A^{B'}[k]$ to $\cU_B\cap\cU_{B'}$ are equal. The open sets $\cU_B$ with $B$ transversal to $A$ form a covering of $\PP(V)$. Thus the collection of $Y_A^B[k]$'s  glue together to give a closed subscheme of $\PP(V)$ whose 
 support  is equal to $Y_A[k]$.  It follows immediately   from  the definitions that the scheme $Y_A[1]$ is equal to the scheme $Y_A$  defined in~\Ref{sec}{prologo}. 
 By~\Ref{prp}{tanquad} we have
 \begin{equation}\label{codykappa}
\cod(Y_A[k],\PP(V))\le \frac{k(k+1)}{2}\quad \text{if}
\quad Y_A[k]\not=\es\,.
\end{equation}
We set
 \begin{equation}\label{yerre}
Y_A(k):=Y_A[k]\setminus Y_A[k+1]\,.
\end{equation}
\subsection{Local equation of $Y_A$}\label{subsec:espsarto}
\setcounter{equation}{0}
Let $A\in\lagr$ and $[v_0]\in \PP(V)$:  we will analyze  $Y_A$ in a neighborhood of $[v_0]$. Let $V_0\subset V$  be a  subspace complementary to $[v_0]$. 
 We identify $V_0$ with the open affine  $(\PP(V)\setminus\PP(V_0))$ via the isomorphism
\begin{equation}\label{apertoaffine}
\begin{matrix}
 V_0 & \overset{\sim}{\lra} & \PP(V)\setminus\PP(V_0) \\
v & \mapsto &  [v_0+v].
\end{matrix}
\end{equation}
(Thus $0\in V_0$ corresponds to $[v_0]$.) Since $Y_A$ is a sextic hypersurface we have
\begin{equation}\label{taylor}
Y_A\cap V_0=V(f_0+f_1+\cdots+f_6),\qquad f_i\in \Sym^i V_0^{\vee}\,,
\end{equation}
where the $f_i$'s are determined up to a common multiplicative non-zero constant. 
We will describe explicitly
the polynomials $f_i$ of~\eqref{taylor} for $i\le \dim(A\cap F_{v_0})$.
First some preliminaries.  Given $v\in V$ we define 
a  quadratic form $\phi^{v_0}_v$ on $F_{v_0}$ as follows. Let $\alpha\in F_{v_0}$; then $\alpha=v_0\wedge\beta$ for some $\beta\in\bigwedge^2 V$. We set
\begin{equation}\label{quadricapluck}
\phi^{v_0}_v(\alpha):=\vol(v_0\wedge v\wedge\beta\wedge\beta).
\end{equation}
The above equation gives a well-defined quadratic form on $F_{v_0}$ because  $\beta$ is determined up to addition by an element of $F_{v_0}$.
\begin{prp}\label{prp:effekappa}
Let $A\in\lagr$. Let $[v_0]\in \PP(V)$ and $V_0\subset V$ be a  subspace complementary to $[v_0]$. Let  $f_i\in \Sym^i V_0^{\vee}$   for $i=0,\ldots,6$ be the  polynomials appearing in~\eqref{taylor}. 
Let  $K:=A \cap F_{v_0}$ and  $k:=\dim K$ . Then
\begin{itemize}
\item[(1)]
$f_i=0$ for $i<k$, and
\item[(2)]
 there exists $\mu\in\CC^{*}$ such that
\begin{equation}\label{sistquad}
f_k(v)=\mu\det(\phi^{v_0}_v|_{K}),\quad v\in V_0,
\end{equation}
where $\phi^{v_0}_v$ is the quadratic form defined by~\eqref{quadricapluck}.
\end{itemize}
\end{prp} 
\begin{proof}
Let $B\in\lagr$ be transversal both to $A$ and $F_{v_0}$. Let $\cV\subset V_0$ be the open subset of $v$ such that $[v_0+v]\in\cU_B$ where $\cU_B$ is given by~\eqref{calub}. Notice that $0\in\cV$. 
For $v\in\cV$ we let $\wt{q}(v):=\tau_A^B([v_0+v])$ where $\tau_A^B$ is given by~\eqref{mappagrafo}.
 Let $q(v)\colon A\to\CC$ be the quadratic form associated to $\wt{q}(v)$. By definition of $Y_A$ we have
\begin{equation}\label{piadet}
Y_A\cap\cV=V(\det q)\,.
\end{equation}
We have $\ker q(0)=A\cap F_{v_0}=K$; by~\Ref{prp}{conodegenere} it follows that $\det q\in\mathfrak{m}^k_0$ where $\mathfrak{m}_0\subset\cO_{\cV,0}$ is the maximal ideal.  This proves Item~(1). 
Let's prove Item~(2). Let  $(\det q)_k\in \left(\mathfrak{m}_0^k/\mathfrak{m}_0^{k+1}\right)\cong \Sym^{k}V_0^{\vee}$ be the \lq\lq initial\rq\rq\, term of $\det q$; by~\eqref{piadet} we have
\begin{equation}\label{giacinto}
f_k=c (\det q)_k,\qquad c\in\CC^{*}\,.
\end{equation}
By~\Ref{prp}{conodegenere}  there exists $\theta\in\CC^{*}$ such that
\begin{equation}\label{carini}
(\det q)_k(v)=\theta\det\left(\left.\frac{d\left(q(vt)|_{K}\right)}{dt}\right|_{t=0}\right).
\end{equation}
A straighforward computation (see Equation~(2.26) of~\cite{og2}) gives  that
\begin{equation}\label{derivata}
\left.\frac{d\left(q(vt)|_{K}\right)}{dt}\right|_{t=0}=\phi^{v_0}_v|_{K}.
\end{equation}
  Item~(2) follows from~\eqref{giacinto}, \eqref{carini} and~\ref{derivata}. 
\end{proof}
In order to apply the above proposition we will need a geometric description of the right-hand side of~\eqref{sistquad}.
Let $[v_0]\in\PP(V)$ and $V_0\subset V$ be complementary to $[v_0]$; we let
\begin{equation}\label{trezeguet}
\begin{matrix}
\lambda^{v_0}_{V_0}\colon\bigwedge^2 V_0 & 
\overset{\sim}{\lra} & F_{v_0}\\
\beta & \mapsto & v_0\wedge\beta
\end{matrix}
\end{equation}
Without choossing a complementary subspace we get an isomorphism
\begin{equation}\label{iaquinta}\index{$\lambda^{v_0}$}
\begin{matrix}
\lambda^{v_0}\colon\bigwedge^2 (V/[v_0]) & 
\overset{\sim}{\lra} & F_{v_0}\\
\ov{\beta} & \mapsto & v_0\wedge\beta
\end{matrix}
\end{equation}
(Here $\ov{\beta}$ is the class represented by $\beta$; the point being that $v_0\wedge\beta$ is indeopendent of the representative.)
Taking inverses we get isomorphisms
\begin{equation}\label{ilpacciani}
 F_{v_0}\overset{\rho^{v_0}_{V_0}}{\lra}\bigwedge^2 V_0,\qquad
 F_{v_0}\overset{\rho^{v_0}}{\lra}\bigwedge^2 (V/[v_0])\,.
\end{equation}
Via $\rho^{v_0}_{V_0}$ we may view $\phi^{v_0}_v$ as a Pl\"ucker  quadratic form on $\bigwedge^2 V_0$. More precisely: given  $v\in V_0$ let $q_v$ be the quadratic form on $\bigwedge^2 V_0$ defined by
\begin{equation}\label{quadpluck}
\begin{matrix}
\bigwedge^2 V_0 & \overset{q_v}{\lra} & \CC \\
\alpha & \mapsto & \vol(v_0\wedge v\wedge \alpha\wedge\alpha)
\end{matrix}
\end{equation}
Then $q_v$ is  a Pl\"ucker quadratic form and we have an isomorphism
\begin{equation}\label{evidente}
\begin{matrix}
V_0 & \overset{\sim}{\lra} & H^0(\cI_{\Gr(2,V_0)}(2)) \\
v & \mapsto & q_v
\end{matrix}
\end{equation}
\begin{rmk}\label{rmk:pluckintr}
We may view $q_v$ as a (Pl\"ucker) quadratic form on $V/[v_0]$ because given $\ov{\alpha}\in\bigwedge^2(V/[v_0])$ the value $\vol(v_0\wedge v\wedge \alpha\wedge\alpha)$ is independent of the representative $\alpha\in\bigwedge^2 V$ of $\ov{\alpha}$.
\end{rmk}
Clearly we have the following relation between  $\phi^{v_0}_v$ and $q_v$:
\begin{equation}\label{eccopluck}
\Sym^2(\rho^{v_0}_{V_0})(\phi^{v_0}_v)=q_v,\qquad v\in V_0 \,.
\end{equation}
Since $\Gr(2,V_0)$ is cut out by quadrics we get that
\begin{equation}\label{intekappa}
\PP(\rho^{v_0}_{V_0}(\bigcap_{v\in V_0}V(\phi^{v_0}_v|_{K})  ))=
{\mathbb Gr}(2,V_0)\cap\PP(\rho^{v_0}_{V_0}(K)).
\end{equation}
\begin{crl}\label{crl:nonsvanisce}
Keep hypotheses  as in~\Ref{prp}{effekappa}.   Then the following hold:
\begin{itemize}
\item[(1)]
Suppose that $A\cap F_{v_0}$ does not contain a non-zero decomposable element of $\bigwedge^3 V$. Then
$Y_A\not=\PP(V)$ and $\mult _{[v_0]}Y_A=\dim(A\cap F_{v_0})$. If moreover $A\cap F_{v_0}$ is one-dimensional, say $A\cap F_{v_0}=\la v_0\wedge\beta\ra$,  then the projective tangent space of $Y_A$ at $[v_0]$ is 
\begin{equation}\label{tangepw}
T_{[v_0]}Y_A=\PP(\supp(v_0\wedge\beta))\,.
\end{equation}
\item[(2)]
If  $A\cap F_{v_0}$  contains a non-zero decomposable element of $\bigwedge^3 V$ then $Y_A$ is singular at $[v_0]$ unless $Y_A=\PP(V)$. 
\end{itemize}
\end{crl}
\begin{proof}
Let's prove Item~(1). Let $K:=A\cap F_{v_0}$ and $k:=\dim K$; we let $f_k$ be the degree-$k$ polynomial appearing in~\eqref{taylor}. By  hypothesis $\Gr(2,V_0)\cap\PP(\rho^{v_0}_{V_0}(K))=\es$. By~\eqref{intekappa}, Bertini's Theorem and~\eqref{evidente}  we get that if $v\in V_0$ is generic the quadratic form $\phi^{v_0}_v |_{K}$ is non-degenerate. Thus
  $f_k\not=0$ by~\eqref{sistquad}. If $\dim(A\cap F_{v_0})=1$ Formula~\eqref{sistquad} gives~\eqref{tangepw}.
  Let's prove Item~(2). Suppose that $Y_A\not=\PP(V)$. Then one at least of the polynomials $f_i$ appearing in~\eqref{taylor} is non-zero; thus it suffices to prove that $f_1=0$. If $k\ge 2$ then we are done  by Item~(1) of~\Ref{prp}{effekappa}. Now assume that $k=1$. Then $\PP(\rho^{v_0}_{V_0}(K))$ is a point contained in $\Gr(2,V_0)$ and hence
  $f_1=0$ by~\eqref{intekappa} and~\eqref{evidente}. 
\end{proof}
 In order to  prove sharper results we will analyze the rational map
\begin{equation}\label{mapachi}
\PP(\bigwedge^2 V_0)\overset{\Phi}{\dashrightarrow}
|\cI_{\Gr(2,V_0)}(2)|^{\vee}\cong\PP(V_0^{\vee})\,.
\end{equation}
Let $Z\subset \PP(\bigwedge^2 V_0)\times\PP(V_0^{\vee})$ be the incidence subvariety defined by
\begin{equation}
Z:=\{([\alpha],[\phi])\mid \phi(\supp\alpha)=0\}\,.
\end{equation}
We have a triangle
\begin{equation}\label{risolvo}
\xymatrix{  
& Z \ar^{\wt{\Phi}}[dr] \ar_{\Psi}[dl]  & \\
 \PP(\bigwedge^2 V_0) &\overset{\Phi}{\dashrightarrow} & \PP( V_0^{\vee})}
\end{equation}
where $\Psi$ and $\wt{\Phi}$ are the restrictions to $Z$ of the two projections of $\PP(\bigwedge^2 V_0)\times\PP(V_0^{\vee})$. 
We will be using the following result; the easy proof is left to the reader.
\begin{lmm}\label{lmm:geomappa}
Keep notation as above - in particular $\dim V_0=5$. Then the following hold:
\begin{itemize}
\item[(1)]
The map $\Psi$ appearing in~\eqref{risolvo} is the blow-up of $\Gr(2,V_0)$.
\item[(2)]
\eqref{risolvo} is a commutative diagram; in fact $Z$ is the graph of the rational map $\Phi$.
\end{itemize}
\end{lmm}
In particular the lemma above states  the following: if $\alpha\in\bigwedge^2 V_0$ is not decomposable then  $\Phi([\alpha])=\supp\alpha$. 
\begin{lmm}\label{lmm:inietto}
Keep notation as above - in particular $\dim V_0=5$. Let ${\bf K}\subset\PP(\bigwedge^2 V_0)$ be a projective subspace which does not intersect $\Gr(2,V_0)$. The map
\begin{equation}\label{mapparis}
\begin{matrix}
{\bf K} & \lra & \Phi({\bf K}) \\
[\alpha] & \mapsto & \Phi([\alpha])
\end{matrix}
\end{equation}
is an isomorphism.
\end{lmm}
\begin{proof}
Let $U\in\Phi({\bf K})$ where $U\subset V_0$ is a subspace of codimension $1$. Then $\PP(\bigwedge^2 U)\subset\PP(\bigwedge^2 V_0)$. We claim that
\begin{equation}\label{solosolo}
\text{$\PP(\bigwedge^2 U)\cap{\bf K}$ is a point.}
\end{equation}
In fact suppose the contrary. Then there exists a line $L\subset(\PP(\bigwedge^2 U)\cap{\bf K})$. Since $\Gr(2,U)\subset\PP(\bigwedge^2 U)$ is a (quadric) hypersurface we get that $L\cap \Gr(2,U)\not=\es$ and this contradicts the hypothesis that ${\bf K}\cap\Gr(2,V_0)=\es$. This proves~\eqref{solosolo}). By commutativity of~\eqref{risolvo} we get that Map~\eqref{mapparis} is bijective with injective differential; it follows that Map~\eqref{mapparis} is an isomorphism.
\end{proof}
\begin{lmm}\label{lmm:rigatacubica}
Keep notation as above - in particular $\dim V_0=5$. Let ${\bf K}\subset\PP(\bigwedge^2 V_0)$ be a projective subspace intersecting $\Gr(2,V_0)$ in a unique point $p_0$ and such that
\begin{equation}
{\bf K}\cap T_{p_0}\Gr(2,V_0)=\{p_0\}\,.
\end{equation}
The restriction of $\Phi$ to ${\bf K}$ is identified with the natural map 
\begin{equation}
{\bf K}  \dashrightarrow |\cI_{p_0}(2)|^{\vee}\,.
\end{equation}
\end{lmm}
\begin{proof}
Suppose the contrary. Then there exists a proper projective subspace ${\bf P}\subset  |\cI_{p_0}(2)|$ such that  the
 restriction of $\Phi$ to ${\bf K}$  is identified with the natural map
 ${\bf K}\dashrightarrow {\bf P}^{\vee}$. It follows that there exists a subscheme $\cZ\subset({\bf K}\setminus\{p_0\})$ of length $2$ over which $\Phi$ is constant. Let $L\subset{\bf K}$ be the line containing $\cZ$. Arguing as in the proof of~\Ref{lmm}{inietto} we get that $L$ intersects $\Gr(2,V_0)$ in two points or is tangent to $\Gr(2,V_0)$; that contradicts our hypothesis.
\end{proof}
\begin{prp}\label{prp:cicciobello}
Let $A\in\lagr$.  Suppose that $[v_0]\in Y_A(2)$ and   that $A\cap F_{v_0}$ does not contain a non-zero decomposable element of $\bigwedge^3 V$.  Then $Y_A\not=\PP(V)$ and the following hold:
\begin{itemize}
\item[1]
$\mult _{[v_0]} Y_A=2$ and the quadric cone $C_{[v_0]} Y_A$  has rank $3$.
\item[2]
 $Y_A[2]$ is  smooth two-dimensional at $[v_0]$.
\end{itemize}
\end{prp}
\begin{proof}
By~\Ref{crl}{nonsvanisce} we know that $Y_A\not=\PP(V)$.
Let $K:= A\cap F_{v_0}$. We claim that the map
\begin{equation}\label{stafilococco}
\begin{matrix}
V_0 & \lra & \Sym^2 K^{\vee} \\
v & \mapsto & \phi^{v_0}_v|_{K}
\end{matrix}
\end{equation}
is surjective. In fact let ${\bf K}:=\PP(\rho^{v_0}_{V_0}(K))$. By hypothesis ${\bf K}$ does not intersect the indeterminacy locus of 
 the map $\Phi$ given by~\eqref{mapachi}. Thus pull-back by $\Phi$ defines a map
 \begin{equation}\label{emmatore}
H^0(\cO_{\PP(V_0^{\vee})}(1))\overset{\Phi^{*}}{\lra} H^0(\cO_{\bf K}(2))\,.
\end{equation}
 By~\Ref{lmm}{inietto} the restriction of $\Phi$ to ${\bf K}$ is injective; since ${\bf K}\cong\PP^1$ it follows that~\eqref{emmatore} is surjective.  By Isomorphism~\eqref{evidente} we get that~\eqref{stafilococco} is surjective. Items~(1), (2) follow from  surjectivity of~\eqref{stafilococco} together with~\eqref{sistquad}, and \eqref{derivata}, \ref{norquad} respectively.
\end{proof}
The following result was proved in~\cite{og2}.
\begin{crl}[{\rm  Proposition~(2.8) of~\cite{og2}}]\label{crl:ybello}
If $A\in\lagr^0$ then $Y_A\not=\PP(V)$ and the following hold:
\begin{itemize}
\item[(1)]
$\sing Y_A$ is a smooth  pure-dimensional surface of degree $40$.
\item[(2)]
$\mult _{[v]}Y_A=2$ for every $[v]\in \sing Y_A$ and the quadric cone $C_{[v]} Y_A$  has rank $3$.
\end{itemize}
\end{crl}
\begin{proof}
\Ref{crl}{nonsvanisce}  gives that $Y_A\not=\PP(V)$ and that $\sing Y_A=Y_A[2]$.
Since $Y_A$ is a global Lagrangian degeneracy locus  there is a Porteous formula that gives the expected cohomology class of $Y_A[k]$ for every integer $k$ - see~\cite{fulpra}. In fact
 Formula~(6.7)
of~\cite{fulpra} and the equation
\begin{equation}\label{cherndieffe}
  c(F)=1-6 c_1(\cO_{\PP(V)}(1))+18 c_1(\cO_{\PP(V)}(1))^2
  -34 c_1(\cO_{\PP(V)}(1))^3+\ldots
\end{equation}
give that the expected cohomology class of $Y_A[2]$ is 
\begin{equation}\label{classew}
\text{exp.class of $Y_A[2]$}=
  2c_3(F)-c_1(F)c_2(F)=40 c_1(\cO_{\PP(V)}(1))^3
\end{equation}
Since the above class is non-zero it follows that $Y_A[2]\not=\es$. By~\Ref{prp}{cicciobello} we get that $Y_A[2]$ is a smooth pure-dimensional surface. Surjectivity of Map~\eqref{stafilococco} gives that the expected class of $Y_A[2]$ is the cohomology class of the smooth surface $Y_A[2]$; thus $\deg Y_A[2]=40$. This proves Item~(1). Item~(2) follows from~\Ref{prp}{cicciobello}.
\end{proof}
We notice that the converse  of~\Ref{crl}{ybello} holds but we will not prove it here.

We close the present subsection by showing how to detect the most pathological  $A\in\lagr$.  Let
\begin{equation}
\NN(V):= \{A\in\lagr \mid Y_A=\PP(V)\}\,.\index{$\NN(V)$}
\end{equation}
We say that a closed  $Z\subset{\mathbb Gr}(3,V)$ is {\it invasive}\index{invasive} if
\begin{equation}
\bigcup\limits_{W\in Z}W=V.
\end{equation}
\begin{clm}\label{clm:seinvadi}
Let $A\in\lagr$. Then $A\in\NN(V)$ if and only if $\Theta_A$ is invasive. In particular if $\dim\Theta_A\le 2$ then $Y_A\notin\PP(V)$. 
\end{clm}
\begin{proof}
Suppose that $A\in\NN(V)$; then $\Theta_A$ is invasive by Item~(1) of~\Ref{crl}{nonsvanisce}. 
Now suppose that $\Theta_A$ is invasive. 
Let $[v]\in\PP(V)$. Since $\Theta_A$ is invasive there exists $W\in\Theta_A$ containing $v$. Then $\bigwedge^3W\subset(A\cap F_v)$ and hence $[v]\in Y_A$. Since $[v]$ was arbitrary we get that $Y_A=\PP(V)$.  This proves that $A\in\NN(V)$ if and only if $\Theta_A$ is invasive. Now suppose that $\dim\Theta_A\le 2$. Then $\dim(\cup_{W\in \Theta_A}W)\le 5$ and hence $\Theta_A$ is not invasive; thus $Y_A\notin\PP(V)$ by the first part of the claim. 
\end{proof}  
\subsection{Morin's Theorem}\label{subsec:vivamorin}
\setcounter{equation}{0}
 Suppose that $\Theta\subset\Theta_A$ for some $A\in\lagr$; by~\Ref{rmk}{adueadueinc} any two $W_1,W_2\in\Theta$ have non-trivial intersection or equivalently $\PP(W_1)$ and $\PP(W_2)$ are incident.
Ugo Morin~\cite{morin} classified maximal irreducible families of pairwise incident  planes in a projective space. In order to recall Morin's result we introduce certain subsets of 
$\Gr(3,V)$.
The first three  are Schubert cycles, namely 
\begin{eqnarray}
J_{v_0}:=& \{W\in\Gr(3,V)\mid v_0\in W\},\quad [v_0]\in\PP(V)
\index{$J_{v_0}$}\label{geivu}\\
I_U:= & \{W\in\Gr(3,V)\mid \dim(W\cap U)\ge 2\},
\quad U\in\Gr(3,V),
\index{$I_U$}\label{iconu}
\end{eqnarray}
and $\Gr(3,E)$ where $E\in\Gr(5,V)$. Next let $\cQ,\cV\subset\PP(V)$ be a smooth quadric hypersurface and a Veronese surface respectively; set
\begin{eqnarray}
F_{+}(\cQ),F_{-}(\cQ):=&\text{irred.~compt.'s of 
$\{W\in\Gr(3,V)\mid \PP(W)\subset\cQ\}$,}
\index{$F(\cQ)$}\label{fanodiqu}\\
C(\cV):= & \{W\in\Gr(3,V)\mid \text{$\PP(W)\cap\cV$ is a conic}\},
\index{$C(\cV)$}\label{cidivu}\\
T(\cV):= & \{W\in\Gr(3,V)\mid \text{$\PP(W)=T_p\cV$ 
for some $p\in\cV$}\}.\index{$T(\cV)$}\label{tidivu}
\end{eqnarray}
(Here $T_p\cV$ is the projective tangent plane to $\cV$ at $p$.) 
As is easily checked Item~(1) of~\Ref{rmk}{adueadueinc} holds for $\Theta$ equal to  one of the six subsets listed above (of course there is no intrinsic difference between $F_{+}(\cQ)$ and $F_{-}(\cQ)$)  - the first three are the {\it elementary complete systems} in Morin's terminology.  
\begin{thm}\label{thm:teomorin}{\rm[U.~Morin~\cite{morin}]}
Let $\Theta\subset\Gr(3,V)$ be a closed irreducible subset such that $\dim (W_1\cap W_2)>0$ for every $W_1,W_2\in\Theta$.  Then 
one of the following holds:
\begin{itemize}
\item[(a)]
There exists a smooth quadric $\cQ\subset\PP(V)$ such that $\Theta\subset F_{\pm}(\cQ)$.
\item[(b)]
There exists a Veronese surface $\cV\subset\PP(V)$ such that
$\Theta$ is contained in $C(\cV)$ or in $T(\cV)$.
\item[(c)]
$\Theta\subset J_{v_0}$ for a certain $[v_0]\in\PP(V)$.
\item[(d)]
$\Theta\subset \Gr(3,E)$ for a certain $E\in\Gr(5,V)$.
\item[(e)]
$\Theta\subset I_U$ for a certain $U\in\Gr(3,V)$.
\end{itemize}
\end{thm}
Let us examine $J_{v_0}$ and $I_U$ in greater detail. 
Let
\begin{equation}\label{robarvu}
\ov{\rho}^{v_0}_{V_0}\colon\PP(F_{v_0})\overset{\sim}{\lra}\PP(\bigwedge^2 V_0)
\end{equation}
be the isomorphism induced by~\eqref{trezeguet}. Restricting $\ov{\rho}_{v_0}$ to $J_{v_0}\subset\PP(F_{v_0})$ we get an isomorphism
\begin{equation}
\begin{matrix}
J_{v_0} & \overset{\sim}{\lra} & \Gr(2,V_0) \\
W & \mapsto & W\cap V_0
\end{matrix}
\end{equation}
Next we examine $I_U$. Given a subspace $U\subset V$ of arbitrary dimension we let
\begin{equation}\label{esseu}
S_U:= (\bigwedge^2 U)\wedge V.\index{$S_U$}
\end{equation}
Of course $S_U=0$ if $\dim U\le 1$ and thus $S_U$ is of interest only for $\dim U\ge 2$. Notice that $S_U\supset \bigwedge^3 U$. We let
\begin{equation}\label{tiutiu}
T_U:= S_U/ \bigwedge^3 U\cong \bigwedge^2 U\otimes (V/U),\index{$T_U$}
\end{equation}
and
\begin{equation}\label{rouma}
\rho_U\colon S_U \lra T_U
\end{equation}
be the quotient map. Now assume that $\dim U=3$ and hence $I_U$ is defined and $I_U\subset\PP(S_U)$. Let
\begin{equation}\label{robaruma}
\ov{\rho}_U\colon \PP(S_U) \dashrightarrow P(T_U)
\end{equation}
be the rational map corresponding to~\eqref{rouma} i.e.~projection from $\PP(\bigwedge^3 U)$.  The following claim gives an explicit description of $I_U$; the easy proof is left to the reader.
\begin{clm}\label{clm:conosegre}
Let $U\in\Gr(3,V)$. Then $I_U$ is
 the cone with vertex $\bigwedge^3 U$ over the Segre variety $\PP(\bigwedge^2 U)\times\PP(W)\subset\PP(T_U)$.
\end{clm}
\subsection{Menagerie}\label{subsec:sespec}
\setcounter{equation}{0}
We show that the following are EPW-sextics: 
$3\cQ$ where $\cQ$ is a smooth quadric, $2\chord(\cV)$ where $\chord(\cV)$ is the chordal variety of a Veronese surface (i.e.~the discriminant cubic parametrizing degenerate plane conics) and the union of six independent hyperplanes. These special EPW-sextics are analogues of certain  cubic $4$-folds and  plane sextic curves which have a special r\^ole in the works of Laza~\cite{laza1,laza2} and  Shah~\cite{shah} - Table~\eqref{analoghi} gives the dictionary between the three cases. 
\begin{table}[tbp]
\caption{Analogous special EPW-sextics, cubic $4$-folds and plane sextics}\label{analoghi}
\vskip 1mm
\centering
\renewcommand{\arraystretch}{1.60}
\begin{tabular}{l l l}
\toprule
  EPW-sextic & cubic $4$-fold & plane sextic \\
\midrule
 triple quadric, double & \multirow{2}{*}{discriminant cubic} & \multirow{2}{*}{triple conic} \\
 discriminant cubic &  & \\
\midrule
union of $6$ independent & \multirow{2}{*}{$V(x_0 x_1 x_2+x_3 x_4 x_5)$} & 
\multirow{2}{*}{$V(x_0^2 x_1^2 x_2^2)$} \\
hyperplanes & & \\
\bottomrule 
\end{tabular}
\end{table} 
 \vskip 2mm
 \n
 $\boxed{\text{\it Triple smooth quadric}}$
 Write $V=\bigwedge^2 U$ where $U$ is a complex vector-space of dimension $4$. Thus $\cQ:=\Gr(2,U)\subset\PP(\bigwedge^2 U)$ is a smooth quadric hypersurface. We have embeddings
\begin{equation}\label{piumenomap}
\qquad\qquad
\begin{matrix}
\PP(U) & \overset{i_{+}}{\hra} & \Gr(3,\bigwedge^2 U)\\
[u_0] & \mapsto & \{u_0\wedge u\mid u\in U\}
\end{matrix}
\qquad
\begin{matrix}
\PP(U^{\vee}) & \overset{i_{-}}{\hra} & \Gr(3,\bigwedge^2 U)\\
[f_0] & \mapsto & \bigwedge^ 2 ker(f_0).
\end{matrix}
\end{equation}
Thus  referring to~\eqref{fanodiqu} we may set 
\begin{equation}\label{piumenofano}
F_{\pm}(\cQ)=\im(i_{\pm})\,. 
\end{equation}
Let $A_{\pm}(U)\subset\bigwedge^3(\bigwedge^2 U)$ be the 
subspaces defined by 
\begin{equation}\label{aconiconj}
A_{\pm}(U)=\la\la \im(i_{\pm})\ra\ra
\end{equation}
\begin{clm}\label{clm:tripquad}
Keep notation as above. Then $A_{\pm}(U)\in \mathbb{LG}(\bigwedge^3(\bigwedge^2 U))$ and furthermore
\begin{equation}\label{trequad}
Y_{A_{\pm}(U)}=3\cQ.
\end{equation}
\end{clm}
\begin{proof}
Since any two planes in $\im(i_{\pm})$ are incident we get that $A_{\pm}(U)$  is isotropic for $(,)_V$ - see~\Ref{rmk}{adueadueinc}. 
Let $\cL$ be the Pl\"ucker(ample) line-bundle on  $\Gr(3,\bigwedge^2 U)$; one checks easily that
\begin{equation}\label{paccheri}
i_{+}^{*}\cL\cong\cO_{\PP(U)}(2),\qquad
i_{-}^{*}\cL\cong\cO_{\PP(U^{\vee})}(2). 
\end{equation}
Since $i_{\pm}$ is $SL(U)$-equivariant and $\Sym^2 U,\Sym^2 U^{\vee}$ are irreducible  $SL(U)$-modules we get surjections
\begin{equation}\label{granrett}
\begin{array}{ccccc}
H^0(\cO_{\PP(U)}(2)) & \overset{i_{+}^{*}}{\longleftarrow} &
H^0(\Gr(3,\bigwedge^2 U)) &
\overset{i_{-}^{*}}{\lra} &  H^0(\cO_{\PP(U^{\vee})}(2)) \\
 \vert\vert & &  \vert\vert & &  \vert\vert \\
\Sym^2 U^{\vee} & \longleftarrow &
\bigwedge^3(\bigwedge^2 U^{\vee}) &
\lra & \Sym^2 U 
\end{array}
\end{equation}
It follows that we have  isomorphisms of $SL(U)$-modules
\begin{equation}\label{babucce}
 A_{+}(U)\cong \Sym^2 U,\qquad A_{-}(U)\cong \Sym^2 U^{\vee}.
\end{equation}
In particular we get that $2\dim A_{\pm}(U)=\dim\bigwedge^3(\bigwedge^2 U)$ and hence  $A_{\pm}(U)$  is Lagrangian.
Lastly we will prove~\eqref{trequad}.
First we claim that 
\begin{equation}\label{tetapiu}
\Theta_{A_{\pm}(U)}=\im(i_{\pm}).
\end{equation}
 Suppose that $W\in \Theta_{A_{+}(U)}$; let's prove that $\PP(W)\subset\cQ$. If $\PP(W)\not\subset\cQ$ then 
 \begin{equation}\label{solocurva}
 \dim(\PP(W)\cap\cQ)=1. 
\end{equation}
On the other hand $\PP(W)$ is incident to every plane parametrized by a point of $\im(i_{+})$ (see~\Ref{rmk}{adueadueinc}); by~\eqref{solocurva} we get that if $p\in\PP(W)\cap\cQ$ then 
\begin{equation}
\dim\{\Lambda\in \im(i_{+})\mid p\in\Lambda\}\ge 2.
\end{equation}
That is absurd: if $p\in\cQ$ the set of planes $\Lambda\subset\cQ$ containing $p$ is a line. This proves  that $\PP(W)\subset\cQ$. 
Thus $W\in (\im(i_{+})\cup \im(i_{-}))$. Since $\PP(W)$ is incident to every plane parametrized by $\im(i_{+})$ we get that $W\in \im(i_{+})$. This proves~\eqref{tetapiu} for $A_{+}(U)$ - the proof for $A_{-}(U)$ is the same mutatis mutandis. By~\eqref{tetapiu} we get that
 $\Theta_{A_{\pm}(U)}$  is not invasive and hence $Y_{A_{\pm}(U)}\not=\PP(\bigwedge^2 U)$ by~\Ref{clm}{seinvadi}.
 On the other hand $Y_{A_{\pm}(U)}$ is $SL(U)$-invariant i.e.~invariant for the action of $SO(\bigwedge^2 U,q)$ where  the zero-locus of the quadratic form $q$ is  $\cQ$. It follows that $Y_{A_{\pm}(U)}$ is a multiple of $\cQ$; since $\deg Y_{A_{\pm}(U)}=6$ we  get~\eqref{trequad}.  
\end{proof}
We notice that $A_{+}(U)\not=A_{-}(U)$ because if $(\Lambda_{+},\Lambda_{-})\in \im(i_{+})\times \im(i_{-})$ is generic then $\Lambda_{+}\cap \Lambda_{-}=\es$; it follows that  the  irreducible decomposition of the $SL(U)$-module 
$\bigwedge^3(\bigwedge^2 U)$ is given by
\begin{equation}\label{decirr}
\bigwedge^3(\bigwedge^2 U)=A_{+}(U)\oplus A_{-}(U)\cong
\Sym^2 U\oplus \Sym^2 U^{\vee}.
\end{equation}
\begin{rmk}\label{rmk:contromano}
Referring to~\eqref{decirr}: the $SL(U)$-action on $\Sym^2 U^{\vee}$ is the inverse of the contragradient action.
\end{rmk}
\vskip 2mm
\n
 $\boxed{\text{\it Double discriminant cubic}}$
Write $V=\Sym^2 L$ where $L$ is a complex vector-space of dimension $3$. Let $(\Sym^2 L)_i\subset \Sym^2 L$ be the subset of symmetric tensors of rank at most $i$. Then $\cV:=\PP((\Sym^2 L)_1)$ is a Veronese surface. The   chordal variety of $\cV$ is the  discriminant cubic
\begin{equation}
\PP((\Sym^2 L)_2)=\chord(\cV)=\bigcup_{\text{$\ell$ chord of $\cV$}} \ell.
\end{equation}
(Tangents to $\cV$ are included  among chords of $\cV$.)
We have embeddings
\begin{equation}\label{kappacca}
\qquad
\begin{matrix}
\PP(L) & \overset{k}{\hra} & \Gr(3,\Sym^2 L)\\
[l_0] & \mapsto & \{l_0\cdot l\mid l\in L\}
\end{matrix}
\qquad
\begin{matrix}
\PP(L^{\vee}) & \overset{h}{\hra} & \Gr(3,\Sym^2 L)\\
[f_0] & \mapsto & \{q\mid  f_0\in ker q\}.
\end{matrix}
\end{equation}
Let $\cL$ be the Pl\"ucker(ample) line-bundle on  $\Gr(3,\bigwedge^2 U)$; one checks easily that
\begin{equation}\label{mentuccia}
k^{*}\cL\cong\cO_{\PP(L)}(3),\qquad
h^{*}\cL\cong\cO_{\PP(L^{\vee})}(3). 
\end{equation}
Let $A_k(L),A_h(L)\subset\bigwedge^3(\Sym^2 L)$ be 
the subspaces defined by 
\begin{equation}\label{aconkconh}
A_k(L)=\la\la \im(k)\ra\ra,\qquad A_h(L)=\la\la \im(h)\ra\ra.
\end{equation}
Arguing as in the proof of~\Ref{clm}{tripquad} one proves the following result.
\begin{clm}\label{clm:discub}
Keep notation as above. Then $A_k(L),A_h(L)\in \mathbb{LG}(\bigwedge^3(\Sym^2 L))$ and 
\begin{equation}\label{duecub}
Y_{A_k(L)}=Y_{A_h(L)}=2\chord(\cV)\,.
\end{equation}
Moreover
\begin{equation}\label{accappa}
\Theta_{A_{k}(L)}=\im(k),\qquad \Theta_{A_{h}(L)}=\im(h).
\end{equation}
\end{clm}
A remark: Equation~\eqref{mentuccia} gives 
the irreducible decomposition of the $SL(L)$-module 
$\bigwedge^3(\Sym^2 L)$
\begin{equation}
\bigwedge^3(\Sym^2 L)\cong
\Sym^3 L\oplus \Sym^3 L^{\vee}.
\end{equation}
\vskip 2mm
\n
 $\boxed{\text{\it Union of six independent hyperplanes}}$
 The following example was worked out together with C.~Procesi.
 Let $\{v_0,\ldots,v_5\}$ be a basis of $V$ and $\{X_0,\ldots,X_5\}$ be the dual basis of $V^{\vee}$. Our special $A_{III}$ ($III$ refers to Type $III$ degeneration) is spanned by decomposable vectors  $v_{i_1}\wedge v_{i_2}\wedge v_{i_3}$ for a suitable collection of $\{i_1,i_2,i_3\}$'s. To simplify notation we denote 
 $v_{i_1}\wedge v_{i_2}\wedge v_{i_3}$ by the characteristic function of $\{i_1,i_2,i_3\}$, i.e.~the string composed of three $0$'s and three $1$'s which has $1$ at place $j$ (we start counting from  $0$) if and only if $j\in\{i_1,i_2,i_3\}$.  With the above notation $A_{III}$ is given by
\begin{equation}\label{combinatoria}
  \begin{matrix}
1  & 1 & 1 & 0 & 0& 0 \\
 1  & 1 & 0 & 1 & 0& 0 \\
1  & 0 & 1 & 0 & 1 & 0 \\
0  & 1 & 1 & 0 & 0& 1 \\
0  & 0 & 1 & 1 & 0& 1 \\
0  & 0 & 1 & 1 & 1& 0 \\
0  & 1 & 0 & 0 & 1 & 1 \\
0  & 1 & 0 & 1 & 1 & 0 \\
1  & 0 & 0 & 1 & 0 & 1 \\
1  & 0 & 0 & 0 & 1 & 1 \\
  \end{matrix}
\end{equation} 
Notice that $A_{III}$ is fixed by the maximal torus $T$ of $SL(V)$ diagonalized by $\{v_0,\ldots,v_5\}$ and that $T$ acts trivially on $\bigwedge^{10}A$. In particular $Y_{A_{III}}=V(cX_0\cdot X_1\cdots X_5)$ for a constant $c$. An explicit computation shows that $[v_0+v_1+\cdots+v_5]\not\in Y_A$ and hence $Y_{A_{III}}=V(X_0\cdot X_1\cdots X_5)$. The reader can check that the following holds.
\begin{clm}
Let $T$ be a maximal torus  of $SL(V)$. Suppose that $A
\in\lagr$  is fixed by $T$ and that $T$ acts trivially on $\bigwedge^{10}A$.
Then $A$ is $SL(V)$-equivalent to $A_{III}$.
\end{clm} 
\subsection{Dual of an EPW-sextic}
\setcounter{equation}{0}
The volume-form~\eqref{volumone} defines a volume-form $\vol^{\vee}\colon\bigwedge^6 V^{\vee}\overset{\sim}{\lra}\CC$. Let  $(,)_{V^{\vee}}$ be the symplectic form 
on $\bigwedge^3 V^{\vee}$ defined by $(\alpha,\beta)_{V^{\vee}}:=\vol^{\vee}(\alpha\wedge\beta)$. We let $\lagrdual$ be the symplectic Grassmannian relative to $(,)_{V^{\vee}}$.  
 Let
\begin{equation}\label{alida}
\begin{matrix}
 \bigwedge^3 V & \overset{\delta_V}{\overset{\sim}{\lra}} & \bigwedge^3 V^{\vee} \index{$\delta_V$}\\
 \alpha & \mapsto & \left(\beta \mapsto (\alpha,\beta)_V\right)
\end{matrix}
\end{equation}
be the isomorphism induced by  $(,)_V$. As is easily checked  
\begin{equation}\label{formesimp}
(\alpha,\beta)_V=(\delta_V(\alpha),\delta_V(\beta))_{V^{\vee}},
\quad \alpha,\beta\in \bigwedge^3 V.
\end{equation}
Thus we have an isomorphism
\begin{equation}\label{specchio}
\begin{matrix}
 \lagr & \overset{\delta_V}{\overset{\sim}{\lra}} & \lagrdual \\
 A & \mapsto & \delta_V(A)=\Ann   A.
\end{matrix}
\end{equation}
The geometric meaning of $Y_{\delta_V(A)}$ is the following~\cite{og2,og4}: if $A$ is generic in $\lagr$   then 
\begin{equation}\label{dualey}
Y_{\delta_V(A)}=Y^{\vee}_A\,,
\end{equation}
where $Y^{\vee}_A$ is the dual of $Y_A$. We list below the images by $\delta_V$ of certain  subspaces of $\bigwedge^3 V$.
\begin{eqnarray}
\delta_V(F_{v_0}) & = & \bigwedge^3 \Ann  (v_0),\qquad\qquad
\qquad\qquad [v_0]\in\PP(V) \label{dualeffe}\\
 \delta_V(\bigwedge^3 W)&  = & \bigwedge^3 \Ann  (W),
\qquad\qquad\quad\qquad W\in\Gr(3,V),\label{dualedecom}\\
\delta_V(A_{+}(U)&  = & A_{-}(U^{\vee}),\qquad\qquad\qquad\qquad
\qquad\dim U=4\label{dualeapiu}\\
\delta_V(A_{k}(U) & = & A_{h}(U^{\vee}).\qquad\qquad\qquad
\qquad\qquad\dim U=3\label{dualeakappa}
\end{eqnarray}
(See~\Ref{subsec}{sespec} for the notation in the last two lines.) 
Let $A\in\lagr$. We notice that~\eqref{dualeffe} gives the following description of  $Y_{\delta_V(A)}$: given $E\in\Gr(5,V)$  then
 \begin{equation}\label{uaidelta}
\text{$E\in Y_{\delta_V(A)}$ if and only if 
$(\bigwedge^3 E)\cap A\not=\{0\}$.}
\end{equation}
Let us examine the action of  $\delta_V$ on $\Sigma(V)$. We have a  canonical identification $\Gr(3,V)=\Gr(3,V^{\vee})$ and~\eqref{dualedecom} gives that 
\begin{equation}\label{tetaduale}
\Theta_A=\Theta_{\delta_V(A)}.
\end{equation}
In particular 
\begin{equation}\label{dualesigma}
\delta_V(\Sigma(V))=\Sigma(V^{\vee}),\quad \delta_V(\Sigma_{\infty}(V))=\Sigma_{\infty}(V^{\vee}).
\end{equation}
\subsection{EPW-sextics as  discriminant loci}\label{subsec:ydiscrimina}
\setcounter{equation}{0}
Let $A\in\lagr$. In~\Ref{subsec}{luogodeg} we described $Y_A$ locally around $[v_0]\in\PP(V)$  as the discriminant locus of a symmetric map of vector-bundles. Recall that in order to do so we need to choose $B\in\lagr$ transversal to $A$ and to $F_{v_0}$. In this subsection we will write out explicitly the equation of $Y_A$ by choosing  $B$ as follows. Let $V_0\subset V$ be a codimension-$1$ subspace transversal to $[v_0]$ and $\cD$ be the direct-sum decomposition
\begin{equation}\label{trasuno}
 V=[v_0]\oplus V_0. 
\end{equation}
Assume moreover that
\begin{equation}\label{trasdue}
(\bigwedge^3 V_0)\cap A=\{0\} 
\end{equation}
(notice that the above condition is equivalent to $\delta_V(A)\notin\NN(V^{\vee})$); then $B=\bigwedge^3 V_0$ is indeed a Lagrangian subspace transversal to $A$ and to $F_{v_0}$.
 The open subset $\cU_{\bigwedge^3 V_0}\subset\PP(V)$ is readily seen to be equal to $(\PP(V)\setminus
\PP(V_0))$; we identify it with $V_0$ via Map~\eqref{apertoaffine}. 
Given $v\in V_0$ we have the map 
\begin{equation}\label{seattle}
\tau_A^{\bigwedge^3 V_0}([v_0+v])\colon A\to \bigwedge^3 V_0
\end{equation}
given by~\eqref{mappagrafo}.
We have the isomorphism 
\begin{equation}\label{volumino}
\begin{matrix}
\bigwedge^2 V_0 & \overset{\sim}{\lra} & \bigwedge^3 V_0^{\vee} \\
\alpha & \mapsto & (\xi\mapsto \vol(v_0\wedge\alpha\wedge\xi))\,.
\end{matrix}
\end{equation}
On the other hand~\eqref{abdual}  gives an isomorphism $\bigwedge^3 V^{\vee}_0\cong A$; composing with~\eqref{volumino} we 
 get an isomorphism 
\begin{equation}
\nu^{\cD}_A\colon\bigwedge^2 V_0 \overset{\sim}{\lra} A\,.
\end{equation}
(The superscript $\cD$ refers to Decomposition~\eqref{trasuno}.) Let
\begin{equation}\label{zinaida}
\wt{q}_A^{\cD}(v):=
\tau_A^{\bigwedge^3 V_0}([v_0+v])\circ\nu^{\cD}_A\colon\bigwedge^2 V_0\to\bigwedge^3 V_0.
\end{equation}
and $q_A^{\cD}(v)$ be the associated quadratic form; thus
\begin{equation}\label{eccoqul}
q_A^{\cD}(v)\in \Sym^2 \bigwedge^2 V_0^{\vee}\,.
\end{equation}
Identify $(\PP(V)\setminus
\PP(V_0))$ with $V_0$ via Map~\eqref{apertoaffine}; by definition we have
\begin{equation}\label{lidia}
Y_A\cap V_0=V(\det(q_A^{\cD}(v)))
\end{equation}
where $v\in V_0$.
We will write write out explicitly the maps introduced above. Given $
\alpha\in\bigwedge^2 V_0$ we have $v_0\wedge\alpha\in F_{v_0}$ and  there exists a unique decomposition
\begin{equation}\label{ondalunga}
v_0\wedge\alpha=\beta+\gamma,\qquad \beta\in A,
\quad \gamma\in\bigwedge^3 V_0\,.
\end{equation}
Wedging both sides of the above equation with elements of $\bigwedge^3 V_0$ we get that 
\begin{equation}\label{montalbano}
\nu_A^{\cD}(\alpha)=-\beta\,.
\end{equation}
Moreover we get that
\begin{equation}\label{montelusa}
\wt{q}_A^{\cD}(v)(\alpha)=-\gamma-v\wedge\alpha
=\wt{q}_A^{\cD}(0)(\alpha)+\wt{q}_v(\alpha)
\end{equation}
where $\wt{q}_v$ is  the symmetric map associated to the Pl\"ucker quadratic form~\eqref{quadpluck}. For future reference we record the following description of $\wt{q}_A^{\cD}(0)$:
\begin{equation}\label{giannibrera}
\wt{q}_A^{\cD}(0)(\alpha)=\gamma\iff
(v_0\wedge\alpha+\gamma)\in A. 
\end{equation}
By~\eqref{lidia} we have the following local description of $Y_A$. 
\begin{prp}\label{prp:ipsdisc}
Let $A\in\lagr$ and  $[v_0]\in\PP(V)$.
Suppose that there exists  a codimension-$1$ subspace $V_0\subset V$ such that~\eqref{trasuno}-\eqref{trasdue} hold. Identify $(\PP(V)\setminus
\PP(V_0))$ with $V_0$ via Map~\eqref{apertoaffine}. 
 Then
\begin{equation}
Y_A\cap V_0=V(\det(q_A^{\cD}(0)+q_v))\,.
\end{equation}
\end{prp}
 \section{EPW-sextics in $\Sigma$}\label{sec:singplan}
 \setcounter{equation}{0}
\subsection{Dimension computations for $\Sigma$}
\setcounter{equation}{0}
Let $\wt{\Sigma}\subset\Gr(3,V)\times\lagr$  be  defined by
\begin{equation}
\wt{\Sigma}:=  \{(W,A)\mid \bigwedge^3 W\subset A\}.
\end{equation}
Given $d\ge 0$ we let $\wt{\Sigma}[d]\subset\Gr(3,V)\times\lagr$  be defined by
\begin{equation}
\wt{\Sigma}[d]:= \{(W,A)\in\wt{\Sigma}\mid \dim(A\cap S_W)\ge (d+1)\}.
\end{equation}
Thus $\wt{\Sigma}:=\wt{\Sigma}[0]$. 
Let $\pi\colon\Gr(3,V)\times\lagr\to \lagr$ be projection; we let
\begin{equation}\label{imagosigma}
\Sigma[d]:=\pi(\wt{\Sigma}[d]).
\end{equation}
Thus $\Sigma:=\Sigma[0]$. 
For a geometric interpretation of $\Sigma[1]$ see~\Ref{prp}{singsig}. Let
\begin{equation}
\Sigma_{+}:=  \{A\in\Sigma\mid |\Theta_A|>1\}.
\end{equation}
\begin{prp}\label{prp:codsigma}
Keep notation as above. 
\begin{itemize}
\item[(1)]
Let $0\le d\le 9$. Then $\Sigma[d]$ is closed irreducible and 
\begin{equation}\label{romalazio}
\cod(\Sigma[d],\lagr)=(d^2+d+2)/2.
\end{equation}
In particular $\Sigma$  is closed irreducible 
of codimension $1$.
\item[(2)]
$\Sigma_{+}$ is an irreducible constructible subset of $\lagr$ of  codimension $2$,  moreover  if $A\in \Sigma_{+}$ is generic then $|\Theta_A|=2$ and $A\notin\Sigma[1]$. 
\end{itemize}
\end{prp}
\begin{proof}
 Since $\wt{\Sigma}[d]$ is closed and $\pi$ is projective we get  that  $\Sigma[d]$ is closed by~\eqref{imagosigma}.  Let $\rho_d\colon\wt{\Sigma}[d]\to \Gr(3,V)$  be (the restriction of) projection. 
 Let $W_0\in\Gr(3,V)$; one describes  $ \rho_d^{-1}(W_0)$
as follows.  Given  $W\in\Gr(3,V)$ let
\begin{equation}\label{econwu}
\cE_W:=(\bigwedge^3 W)^{\bot}/\bigwedge^3 W
\end{equation}
where orthogonality is with respect to  $(,)_V$. The symplectic form $(,)_V$ induces a symplectic form on $\cE_W$ and hence we have an associated symplectic grassmannain $\lagre$; notice that $T_W\in\lagre$ where $T_W$ is defined by~\eqref{tiutiu}. Let $T_W[d]:= \{B\in\lagre\mid \dim(B\cap T_W)\ge d\}$.
We have an isomorphism
\begin{equation}\label{sigmawu}
\begin{matrix}
\rho_d^{-1}(W_0) & \overset{\sim}{\lra} & T_{W_0}[d] \\
(W_0,A) & \mapsto & A/\bigwedge^3 W_0
\end{matrix}
\end{equation}
We claim that  
\begin{equation}\label{cotiwu}
\cod(T_W[d],\lagre)=d(d+1)/2.
\end{equation}
In fact let $B_0\in T_W[d]$. Let $C\in\lagre$ be transversal both to $B_0$ and $T_W$. Let $U_C\subset\lagre$ be the open set given by~\eqref{eccoua} (beware: the r\^oles of $B$ and $C$ have been exchanged). Then $B_0\in U_{C}$ and we have an isomorphism $U_{C}\cong \Sym^2 C$ given by~\eqref{identifico}. Via this isomorphism $T_W[d]$ is identified with the subset $(\Sym^2 C)_d  \subset \Sym^2 C$  of symmetric tensors of corank at least $d$; by~\Ref{prp}{tanquad} we have $\cod((\Sym^2 C)_d,  \Sym^2 C)=d(d+1)/2$ and hence we get that~\eqref{cotiwu} holds. By~\Ref{prp}{globcalcul} and~\eqref{cotiwu} we get that 
\begin{multline}\label{dimsigmat}
\dim\wt{\Sigma}[d]= \dim T_W[d]+\dim\Gr(3,V)=
\dim\lagre-\frac{d^2+d}{2}+\dim\Gr(3,V)= \\
=\dim\lagr-10-\frac{d^2+d}{2}+9=\dim\lagr-\frac{d^2+d+2}{2}.
\end{multline}
One checks also  that $T_W[d]$ is irreducible (recall that $0\le d\le 9$, in particular $T_W[d]$ is not empty); it follows that $\wt{\Sigma}[d]$ is irreducible and hence  $\Sigma[d]$ is irreducible as well. Summing up: we have proved that $\Sigma[d]$ is irreducible and that its codimension is at least the right-hand side of~\eqref{romalazio}. Moreover in order to finish the proof of Item~(1) it suffices to show that the restriction of projection $\wt{\Sigma}[d]\to\Sigma[d]$ is birational.
 Let $U\subset\Gr(3,V)\times \Gr(3,V)$ be 
\begin{equation}
U:=\{(W_1,W_2)\mid 0<\dim(W_1\cap W_2)<3\}
\end{equation}
and $\wt{\Sigma}_{+}[d]\subset U\times\lagr$ be
 \begin{equation}
\wt{\Sigma}_{+}[d]:=
\{(W_1,W_2,A)\mid (W_i,A)\in \wt{\Sigma}[d]\quad i=1,2\}.
\end{equation}
We laim that in order to prove that $\wt{\Sigma}[d]\to\Sigma[d]$ is birational it suffices to show that
 \begin{equation}\label{rubymacri}
\dim\wt{\Sigma}_{+}[d]<\dim\wt{\Sigma}[d].
\end{equation}
In fact let's grant~\eqref{rubymacri} and let's suppose that $\wt{\Sigma}[d]\to\Sigma[d]$ is not birational. By~\Ref{rmk}{adueadueinc} we get that there is an open dense subset of $\wt{\Sigma}[d]$ which is in the image of the forgetful map
\begin{equation}
\begin{matrix}
\wt{\Sigma}_{+}[d] & \overset {f_d}{\lra} & \wt{\Sigma}[d] \\
(W_1,W_2,A) & \mapsto & (W_1,A) 
\end{matrix}
\end{equation}
and that contradicts~\eqref{rubymacri}. Let's proceed to prove~\eqref{rubymacri}. 
Let  $\eta_d\colon\wt{\Sigma}_{+}[d]\to U$ be the (restriction of) projection. Let $(W_1,W_2)\in U$; the fiber $\eta_d^{-1}(W_1,W_2)$ is described as follows. Let
\begin{equation}
\cE_{W_1,W_2}:= (\bigwedge^3 W_1\oplus\bigwedge^3 W_2)^{\bot}/(\bigwedge^3 W_1\oplus\bigwedge^3 W_2).
\end{equation}
We have an inclusion
\begin{equation}
\begin{matrix}
\eta_d^{-1}(W_1,W_2) & \overset{\theta_d}{\hra} & \LL\GG(\cE_{W_1,W_2}) \\
(W_1,W_2,A) & \mapsto & A/(\bigwedge^3 W_1\oplus \bigwedge^3 W_2)
\end{matrix}
\end{equation}
The above map is bijective if and only if $d=0$. In order to describe the image for $d>0$ we let $\cT^{W_i}_{W_1,W_2}\subset \cE_{W_1,W_2}$ be defined as
\begin{equation*}
\cT^{W_i}_{W_1,W_2}:=\im((S_{W_i}\cap (\bigwedge^3 W_1\oplus\bigwedge^3 W_2)^{\bot})
 \lra \cE_{W_1,W_2}).
\end{equation*}
Let $\cT_{W_1,W_2}[d]:= \{B\in \LL\GG(\cE_{W_1,W_2}) \mid 
\dim(B\cap \cT^{W_i}_{W_1,W_2})\ge d\quad i=1,2\}$. Clearly $\im\theta_d=\cT_{W_1,W_2}[d]$. As is easily checked $\cT^{W_i}_{W_1,W_2}\in \LL\GG(\cE_{W_1,W_2})$; thus arguing as in the proof of~\eqref{cotiwu} we get that
\begin{equation*}
\cod(\im\theta_d,\LL\GG(\cE_{W_1,W_2}))=
\cod(\cT_{W_1,W_2}[d],\LL\GG(\cE_{W_1,W_2}))\ge d(d+1)/2.
\end{equation*}
We have $\dim U=17$; by~\Ref{prp}{globcalcul}  we get that 
\begin{multline}\label{sigmatildestar}
\dim\wt{\Sigma}_{+}[d]\le
 \dim \LL\GG(\cE_{W_1,W_2}))-\frac{d(d+1)}{2}+\dim U= \\
=\dim\lagr-19-\frac{d(d+1)}{2}+17= 
\dim\lagr-\frac{d^2+d+2}{2}-1.
\end{multline}
Thus~\eqref{rubymacri} follows from the above inequality and~\eqref{dimsigmat}. This finishes the proof of Item~(1).
 Let's prove Item~(2).
We have
$\Sigma_{+}=\pi\circ f_0(\wt{\Sigma}_{+}[0])$.
Since $\wt{\Sigma}_{+}[0]$ is constructible we get that $\Sigma_{+}$ is constructible and by~\eqref{sigmatildestar} we get that $\cod(\Sigma_{+},\lagr)\ge 2$. A dimension count similar to those performed above gives that $\cod(\Sigma_{+},\lagr)= 2$ and that the generic $A\in\Sigma_{+}$ has the properties stated in Item~(2). 
\end{proof}
\subsection{First order computations}
\setcounter{equation}{0}
In  proving~\Ref{prp}{codsigma} we have shown that $\wt{\Sigma}$ is a locally trivial fibration over $\Gr(3,V)$ with fiber $\lagre$ over $W$; thus $\wt{\Sigma}$ is smooth. Let $\rho:=\pi|_{\wt{\Sigma}}\colon \wt{\Sigma}\to\Sigma$. The differential of $\rho$ at $(W,A)\in\wt{\Sigma}$ is expressed as follows. 
Of course $T_{(W,A)}\wt{\Sigma}\subset T_W\Gr (3,V)\oplus T_A\lagr$.
Choosing a volume form on $W$ i.e.~a generator $\alpha$ of $\bigwedge^3 W$ we have isomorphisms
\begin{equation*}
T_W\Gr(3,V)=\Hom (W,V/W)=\bigwedge^2 W\otimes(V/W)=S_W/\bigwedge^3 W=T_W.
\end{equation*}
Let $\varphi\colon\bigwedge^3 V\to \bigwedge^3 V/A$ be the quotient map. Given $\tau\in S_W/\bigwedge^3 W$ we let $\wt{\varphi}(\tau):=\varphi(\wt{\tau})$ where $\wt{\tau}\in S_W$ is an element representing the equivalence class $\tau$;   this makes sense because $\bigwedge^3 W\subset A$. 
On the other hand the tangent space $T_A\lag$ is given by~\Ref{prp}{globcalcul}: we have a canonical identification
\begin{equation*}
T_A\lagr\cong\{\theta\colon A\to A^{\vee}\mid \theta^t=\theta\}.
\end{equation*}
Given the above identifications one has
\begin{equation}\label{diffro}
T_{(W,A)}\wt{\Sigma}=\{(\tau,\theta)\in T_W\times \Sym^2 A^{\vee} 
\mid \theta(\alpha)=\wt{\varphi}(\tau)\}.
\end{equation}
(The proof consists of a straightforward computation.)  In particular we get that
\begin{equation}\label{nucleoro}
\ker d\rho(W,A)=A\cap T_W.
\end{equation}
Thus we have the following interpretation of $\wt{\Sigma}[1]$:
\begin{equation}\label{critici}
\wt{\Sigma}[1]=
\{(W,A)\in\wt{\Sigma}\mid \text{$d\pi(W,A)$ is not injective}\}
\end{equation}
\begin{prp}\label{prp:singsig}
The irreducible decomposition of $\sing\Sigma$ is equal to
\begin{equation}
\sing\Sigma=\ov{\Sigma}_{+}\cup\Sigma[1]\,.
\end{equation}
Both irreducible components are of codimension $1$ in $\Sigma$.
\end{prp}
\begin{proof}
By~\Ref{prp}{codsigma} we know that   $\Sigma[1]$ is irreducible of dimension $53$ and that $\ov{\Sigma}_{+}\not=\Sigma[1]$. By Item~(2) of~\Ref{prp}{codsigma} the map
$\rho \colon\wt{\Sigma}\lra\Sigma$ is birational. Thus~\ref{critici} gives that $\Sigma[1]\subset \sing\Sigma$. Since $\Sigma_{+}\not\subset\Sigma_{\infty}$ we also have $\ov{\Sigma}_{+}\subset \sing\Sigma$. Lastly $\Sigma$ is smooth away from $\ov{\Sigma}_{+}\cup\Sigma[1]$ by~\eqref{nucleoro}.
\end{proof}
\begin{prp}\label{prp:trasversi}
Let $A\in\lagr$ and $W_1,W_2\in\Gr (3,V)$. Suppose that for $i=1,2$
\begin{equation}\label{interesse}
A\cap S_{W_i}=\bigwedge^3 W_i.
\end{equation}
Then $\im(d\rho(W_1,A)$ and $\im(d\rho(W_2,A)$ are codimension-$1$ transverse subspaces of $T_A\lagr$.
\end{prp}
\begin{proof}
Both $\im(d\rho(W_1,A)$ and $\im(d\rho(W_2,A)$ are codimension-$1$ subspaces of $T_A\lagr$ by~\eqref{nucleoro}; it remains to prove that they are distinct. Let $\alpha_i$ be a generator of $\bigwedge^3 W_i$; Formula~\eqref{diffro} gives that
\begin{equation*}
\im d\rho(W_i,A)=\{\theta\in  \Sym^2 A^{\vee} 
\mid \theta(\alpha_i)\in\varphi(S_{W_i})\}.
\end{equation*}
It follows that $\im(d\rho(W_1,A)=\im(d\rho(W_2,A)$ if and only if $\varphi(S_{W_i})=\ker(\alpha_{3-i})$. Since $S_{W_i}$ is lagrangian that is possible only if $\alpha_{3-i}\in S_{W_i}$; that is absurd by~\eqref{interesse}.
\end{proof}
The following result is a straightforward consequence of~\Ref{prp}{codsigma} and~\Ref{prp}{trasversi}.
\begin{crl}\label{crl:sigmapiu}
If $A\in\Sigma_{+}$ is generic then $\Sigma$ has normal crossings at $A$ with exactly two sheets.
\end{crl}
\subsection{One-dimensional components of $\Theta_A$}
\label{compuno}
\setcounter{equation}{0}
In this subsection we will classifiy couples $(A,\Theta)$ where $A\in\Sigma$ and $\Theta$ is a $1$-dimensional irreducible component of $\Theta_A$ - of course our point of departure is Morin's~\Ref{thm}{teomorin}.  
\begin{dfn}\label{dfn:isoiso}
Let $\Theta\subset\Gr(3,V)$ be closed: it is {\it isotropic} if $\la\la\Theta\ra\ra\subset\bigwedge^3 V$ is an isotropic subspace or equivalently $W_1\cap W_2\not=\{0\}$ for all $W_1,W_2\in\Theta$, it is {\it isolated isotropic} if in addition it is a union of irreducible components of $\la\Theta\ra\cap\Gr(3,V)$. 
\end{dfn}
The following is an immediate consequence of~\Ref{rmk}{adueadueinc}.
\begin{rmk}\label{rmk:compisola}
Suppose that $A\in\lagr$ and that $\Theta$ is an irreducible component of $\Theta_A$. Then $\Theta$ is isolated isotropic.
\end{rmk}
Before stating our main result on  isolated isotropic curves in $\Gr(3,V)$ we go through some elementary remarks on projective families of  planes in $\PP(V)$.
Let $\Theta\subset \Gr(3,V)$ be an irreducible closed subset. We let $\cE_{\Theta}\to\Theta$\index{$\cE_{\Theta}$} be the restriction of the tautological rank-$3$ vector-bundle on $\Gr(3,V)$ - thus the dual $\cE^{\vee}_{\Theta}$ is globally generated. We let $R_{\Theta}\subset\PP(V)$ be the variety swept out by the $2$-dimensional projective spaces parametrized by $\Theta$, i.e.
\begin{equation}
R_{\Theta}:=\bigcup_{W\in \Theta}\PP(W).
\end{equation}
Let
\begin{equation}\label{mappeffe}
 f_{\Theta}\colon\PP(\cE_{\Theta})\to R_{\Theta} 
\end{equation}
be the tautological surjective map. We may factor $f_{\Theta}$ as follows. The surjective evaluation map $H^0(\cE_{\Theta}^{\vee})\otimes\cO_{\Theta}\to \cE_{\Theta}^{\vee}$ defines a map 
\begin{equation}
h_{\Theta}\colon \PP(\cE_{\Theta})\to \PP(H^0(\cE_{\Theta}^{\vee})^{\vee}).
\end{equation}
 Let $T_{\Theta}:=\im(h_{\Theta})$;  there is a natural map
\begin{equation}
g_{\Theta}\colon T_{\Theta}\lra R_{\Theta},\qquad
f_{\Theta}=g_{\Theta}\circ h_{\Theta}.
\end{equation}
The pull-back by $g_{\Theta}$ of the hyperplane line-bundle on $\PP(V)$ is isomorphic to the hyperplane line-bundle on $T_{\Theta}$; thus $g_{\Theta}$ is either an isomorphism or it may be identified with a projection of $T_{\Theta}$. 
Now assume that $\dim\Theta=1$. Then $\dim R_{\Theta}=3$ and hence $f_{\Theta}$ is of finite degree; one has
\begin{equation}\label{gradorigata}
\deg\Theta=c_1(\cE^{\vee}_{\Theta})=\deg f_{\Theta}\cdot\deg R_{\Theta}.
\end{equation}
\begin{prp}\label{prp:gradoeffe}
Suppose that $\Theta\subset\Gr(3,V)$ is an isolated isotropic   irreducible curve. Then 
$\deg f_{\Theta}=\deg g_{\Theta}=1$ and $\deg\Theta=\deg R_{\Theta}=\deg T_{\Theta}$.
\end{prp}
\begin{proof}
Suppose that $\deg f_{\Theta}>1$; we will reach a contradiction. Since $\deg f_{\Theta}>1$ the generic point of $R_{\Theta}$ (i.e.~the generic point on the generic plane parametrized by $\Theta$) is contained in two distinct planes parametrized by $\Theta$. Since $\dim\Theta=1$ it follows that
two distinct planes parametrized by $\Theta$ meet in a line. Hence  either all planes in $\Theta$ contain a fixed line or else they are all contained in a fixed $3$-dimensional projective space $\PP(U)$. If the former holds then $\deg f_{\Theta}=1$, that is a contradiction. Thus we may assume that $\Theta\subset\Gr(3,U)$ where $U\subset V$ is of dimension $3$ and that $\Theta$ is not  a line, in particular $2\le\dim\la\Theta\ra$. On the other hand 
\begin{equation}
\la\Theta\ra\subset \PP(\bigwedge^3 U)=\Gr(3,U)\subset\Gr(3,V).
\end{equation}
Hence the linear space  $\la\Theta\ra$ is contained in $\Gr(3,V)$; since it has dimension at least $2$ we get that $\Theta$ is not an irreducible component of $\la\Theta\ra\cap \Gr(3,V)$, that
 contradicts~\Ref{dfn}{isoiso}. This proves that $\deg f_{\Theta}=1$; since $f_{\Theta}=g_{\Theta}\circ h_{\Theta}$ it follows that $\deg g_{\Theta}=1$ as well. The equality $\deg\Theta=\deg R_{\Theta}=\deg T_{\Theta}$ follows from $\deg f_{\Theta}=\deg g_{\Theta}=1$ together with  Equation~\eqref{gradorigata}.  
\end{proof}
Table~\eqref{catundim} lists families of curves in $\Gr(3,V)$  and assigns a Type to each family - notice that there are   calligraphic and boldface Types, see~\Ref{rmk}{eccoperche} for an explanation of the difference. 
A few comments on Table~\eqref{catundim}. 
In the last four rows of Table~\eqref{catundim} we refer to~\eqref{piumenomap} and~\eqref{kappacca}.
We notice that Table~\eqref{catundim} is preserved by duality.  
More precisely let 
\begin{equation}
\begin{matrix}
\Gr(3,V) & \overset{\wt{\delta}_V}{\lra} & \Gr(3,V^{\vee}) \\
W & \mapsto & \Ann  (W).
\end{matrix}
\end{equation}
If $\Theta$ belongs to one of the familes in Table~\eqref{catundim} then 
$\wt{\delta}_V(\Theta)$  belongs to one of the familes as well - for this to make sense 
we choose an isomorphism $\PP(V)\cong\PP(V^{\vee})$.
Notation in Table~\eqref{catundim} makes it clear what is the Type of $\wt{\delta}_V(\Theta)$ given the Type of $\Theta$. All the asserted dualities are clear except possibly  for the Types $\cE_2,\cE_2^{\vee}$. Let $\Theta_1$ be of Type $\cE_2$ and $\Theta_2:=\wt{\delta}_V(\Theta_1)$. We must check that $\cE_{\Theta_2}\cong\cO_{\PP^1}(-1)^3$. We have a natural exact sequence
\begin{equation*}
0\to \cE_{\Theta_1}\overset{\alpha}{\lra} \cO_{\PP^1}^6\lra \cE^{\vee}_{\Theta_2}\to 0.
\end{equation*}
(The map $\wt{\delta}_V$ identifies  $\Theta_1$  and $\Theta_2$, moreover they are both isomorphic to $\PP^1$.)
On the other hand $\coker{\alpha}\cong \cO_{\PP^1}^3(1)$ and the result follows.
\begin{clm}\label{clm:ducalberti}
 Let  $X$ be one of the  Types appearing in Table~\eqref{catundim}. Let $\Theta$ be of Type $X$: then $\Theta$ is isolated isotropic. Suppose in addition that $\Theta$ is generic of Type $X$
 (this makes sense because the relevant parameter spaces are irreducible): then  the scheme-theoretic intersection $\la\Theta\ra\cap\Gr(3,V)$ is a smooth irreducible curve, set-theoretically equal to $\Theta$.
\end{clm}
\begin{proof}
 Let $\Theta$ be of Type $X$. Then $\Theta$  is contained in one of the maximal irreducible families of pairwise incident planes in $\PP(V)$ listed in~\Ref{subsec}{vivamorin} and hence it is isotropic. This is 
trivially verified except possibly for $\Theta$ of Type ${\bf Q}$: in that case 
 notice that  the projection from $p$ of $\la h_{\Theta}(\PP(\cO_{\PP^1}(-1)^2))\ra$ (a plane, call it $\PP(U)$) intersects the projection of an arbitrary plane in $T_{\Theta}$ along a line and hence $\Theta\subset I_U$.  
 If $\Theta$ is a line (Type $\cF_1$) then the remaining statements of the claim are trivially true.  From now on we may assume that $X$ is one of the remaining Types. Let $\Theta$ be generic of Type $X$: 
 we must prove that we have equality of sets
\begin{equation}\label{soloteta}
\la\Theta\ra\cap \Gr(3,V)=\Theta
\end{equation}
and that the scheme-theoretic intersection on the left is reduced (it is clear that $\Theta$ is smooth). 
Suppose first that $\Theta$ is generic of Type $\cA$, $\cA^{\vee}$ or $\cC_2$. Then~\eqref{soloteta} holds tautologically. Moreover let $[v_0]\in\PP(V)$, $E\in\Gr(5,V)$ and $U\in\Gr(3,V)$; a straighforward computation with tangent spaces shows that the scheme-theoretic intersections  $\PP(F_{v_0})\cap \Gr(3,V)$, $\PP(\bigwedge^3 E)\cap \Gr(3,V)$ and $\PP(S_U)\cap \Gr(3,V)$ are smooth at every point except for the last case and the point $U$ itself. From this we get that  the scheme-theoretic intersection in~\eqref{soloteta} is reduced. Now suppose that $\Theta$ belongs to one of the remaining Types; then  it belongs  to one of $\im(i_{+})$, $\im(k)$, $\im(h)$. More precisely one of the following holds:
\begin{enumerate}
\item[(1)]
There exist an isomorphism $V\cong\bigwedge^2 U$ where $\dim U=4$ and a curve $C\subset\PP(U)$ such that $\Theta=i_{+}(C)$. Moreover $C$ is cut out scheme-theoretically by quadrics.
\item[(2)]
There exist an isomorphism $V\cong\Sym^2 L$ where $\dim L=3$ and a curve $C\subset\PP(L)$ such that $\Theta=k(C)$. Moreover $C$ is cut out scheme-theoretically by cubics.
\item[(3)]
There exist an isomorphism $V\cong\Sym^2 L^{\vee}$ where $\dim L=3$ and a curve $C\subset\PP(L^{\vee})$ such that $\Theta=h(C)$. Moreover $C$ is cut out scheme-theoretically by cubics.
\end{enumerate}
In fact one of the items above holds by definition if $\Theta$ is of Type ${\bf R}$, ${\bf S}$, ${\bf T}$ or ${\bf T}^{\vee}$. If $\Theta$ is a conic then Item~(1) holds with $C$ a line. If $\Theta$ is of Type $\cE_2$  then Item~(2) holds with $C$ a line, if $\Theta$ is of Type $\cE^{\vee}_2$  then Item~(3) holds with $C$ a line. Lastly if $\Theta$ is of Type ${\bf Q}$ 
then Item~(1) holds with $C$ a conic. 
Now suppose that Item~(1) holds: since $i_{+}$ is defined by the complete linear system of quadrics it follows that~\eqref{soloteta} holds. Similarly if Item~(2) or~(3) holds then we get~\eqref{soloteta} because $k$ and $h$ are defined by the complete linear system of cubics. It remains to show that the scheme-theoretic intersection~\eqref{soloteta} is reduced. Refering to Items~(1), (2) and~(3) above the reduced curve $C$ is the scheme-theoretic intersection of quadrics if Item~(1) holds and the scheme intersection of cubics if Item~(2) or (3) holds: it follows that it suffices to show that the intersections $\PP(A_{+}(U))\cap\Gr(3,V)$, $\PP(A_{k}(L))\cap\Gr(3,V)$ and $\PP(A_{h}(L^{\vee}))\cap\Gr(3,V)$ are reduced. Consider  the first intersection.
Let $W=i_{+}([u_0])\in \PP(A_{+}(U))\cap\Gr(3,V)$ and suppose that the intersection is not reduced at $W$. Acting with the stabilizer of $[u_0]$ in $PGL(U)$ we get that the tangent space at $W$ of the scheme theoretic intersection $\PP(A_{+}(U))\cap\Gr(3,V)$ is all of the tangent space of $\PP(A_{+}(U))$ at $W$. Since $\PP(S_W)$ is the projective tangent space to $\Gr(3,V)$ (embedded in $\PP(\bigwedge^3 V)$) we get that
$\PP(A_{+}(U))\subset\PP(S_W)$ and hence they are equal because they have the same dimension.  This holds for each $W\in\im(i_{+})$: that is absurd because if $W_1\not= W_2$ then $S_{W_1}\not= S_{W_2}$. A similar argument shows that the scheme-theoretic intersections $\PP(A_{k}(L))\cap\Gr(3,V)$ and $\PP(A_{h}(L^{\vee}))\cap\Gr(3,V)$ are reduced. 
\end{proof}
Below is one the main results of the present subsection. 
\begin{thm}\label{thm:tetadimuno}
An isolated isotropic  irreducible curve $\Theta\subset\Gr(3,V)$  belongs to one of the  Types  of Table~\eqref{catundim}. 
\end{thm}
\begin{table}[tbp]\scriptsize
\caption{Types of one-dimensional components of $\Theta_A$}\label{catundim}
\vskip 1mm
\centering
\renewcommand{\arraystretch}{1.60}
\begin{tabular}{llllll}
\toprule
$\Theta$ &   $\deg\Theta$ & $\cE_{\Theta}$ & $R_{\Theta}$ isomorphic to & $\dim\la\Theta\ra$ & 
Type    \\
\midrule
 line & $1$ & $\cO_{\PP^1}^2\oplus\cO_{\PP^1}(-1)$
  & $T_{\Theta}$ & $1$ &  $\cF_1$ \\
\midrule
conic  & $2$ & $\cO_{\PP^1}\oplus\cO_{\PP^1}(-1)^2$  & 
$T_{\Theta}$ & $2$ &  $\cD$ \\
\midrule
rat'l normal cubic  & $3$  & $\cO_{\PP^1}\oplus\cO_{\PP^1}(-1)\oplus\cO_{\PP^1}(-2)$  &  $T_{\Theta}$ & $3$ & $\cE_2$ \\
\midrule
\multirow{2}{*}{rat'l normal cubic} & \multirow{2}{*}{$3$} & \multirow{2}{*}{$\cO_{\PP^1}(-1)^3$}  &  
proj.~of $T_{\Theta}$ & \multirow{2}{*}{$3$} & \multirow{2}{*}{$\cE^{\vee}_2$} \\
   & & & from a point  & \\
\midrule
rational normal  & \multirow{2}{*}{$4$} &  \multirow{2}{*}{$\cO_{\PP^1}(-1)^2\oplus\cO_{\PP^1}(-2)$}  &  
proj.~of $T_{\Theta}$ from  & \multirow{2}{*}{$4$} & 
\multirow{2}{*}{$\bf Q$} \\
quartic  & & &  $p\in\la h_{\Theta}(\PP(\cO_{\PP^1}(-1)^2))\ra$  & \\
\midrule
 in $J_{v_0}$, $[v_0]\in \PP(V)$  & $5$ &   &   & $4$ & $\cA$  \\
\midrule
 in $\Gr(3,E)$,  & \multirow{2}{*}{$5$} &   &   & \multirow{2}{*}{$4$} & \multirow{2}{*}{$\cA^{\vee}$}  \\
$E\in\Gr(5,V)$  & &   & & \\
\midrule
 in $(I_U\setminus\{U\})$,  & \multirow{2}{*}{$6$} &   &   & \multirow{2}{*}{$5$} & \multirow{2}{*}{$\cC_2$}  \\
$U\in \Gr(3,V)$  & &   & & \\
\midrule
$i_{+}(C)$, $C\subset\PP(U)$ a  & \multirow{2}{*}{$6$} &   &   & \multirow{2}{*}{$6$} & \multirow{2}{*}{$\bf R$} \\
rat'l normal cubic  & &   & & \\
\midrule
$i_{+}(C)$, $C\subset\PP(U)$ a  & \multirow{2}{*}{$8$} &  &   & \multirow{2}{*}{$7$} & \multirow{2}{*}{$\bf S$} \\
c.i.~of $2$ quadrics  & &   & & \\
\midrule
$k(C)$, $C\subset\PP(L)$ a & \multirow{2}{*}{$9$} &  &   & \multirow{2}{*}{$8$} & \multirow{2}{*}{$\bf T$} \\
cubic  & &   & & \\
\midrule
$h(C)$, $C\subset\PP(L^{\vee})$ a & \multirow{2}{*}{$9$} &  &   & \multirow{2}{*}{$8$} & 
\multirow{2}{*}{$\bf T^{\vee}$} \\
cubic  & &   & & \\
\bottomrule 
\end{tabular}
\end{table} 
Before proving~\Ref{thm}{tetadimuno} we will  give a series of preliminary results. 
\begin{lmm}\label{lmm:zerozerono}
Let  $\Theta\subset\Gr(3,V)$ be  isolated isotropic.
If we have an inclusion of vector-bundles $\cO_{\Theta}^2\subset\cE_{\Theta}$ 
then $\Theta$ is a linear space.
\end{lmm}
\begin{proof}
By hypothesis there exists  $U\in\Gr(2,V)$ such that 
\begin{equation}
\Theta\subset\{W\in\Gr(3,V)\mid U\subset W\}\cong\PP(V/U).
\end{equation}
It follows that $\la\Theta\ra\subset\Gr(3,V)$. By~\Ref{dfn}{isoiso} we get that  $\Theta=\la\Theta\ra$.
\end{proof}
\begin{prp}\label{prp:quandoconica}
Let  $\Theta\subset\Gr(3,V)$ be  isolated isotropic and suppose that it is a conic;  then  it is of Type $\cD$. 
\end{prp}
\begin{proof}
By hypothesis $\Theta\cong\PP^1$ and hence $\cE_{\Theta}\cong\cO_{\PP^1}(-a_1)\oplus\cO_{\PP^1}(-a_2)\oplus\cO_{\PP^1}(-a_3)$ where $0\le a_i$ and $\sum_i a_i=2$.
By~\Ref{lmm}{zerozerono}  we get that
\begin{equation}\label{pinnapapa}
\cE_{\Theta}\cong\cO_{\PP^1}(-1)^2\oplus\cO_{\PP^1}.
\end{equation}
 It follows that $R_{\Theta}$ is isomorphic either  to $T_{\Theta}$ (a $3$-dimensional quadric of rank $4$)  or to a projection of such a quadric. If the latter holds then $\deg f_{\Theta}=2$ contradicting~\Ref{prp}{gradoeffe}. 
\end{proof}
\begin{prp}\label{prp:cubgobba}
Let  $\Theta\subset\Gr(3,V)$ be  isolated isotropic and suppose that it is a  cubic rational normal curve; then it is either of Type $\cE_2$ or of Type $\cE^{\vee}_2$. 
\end{prp}
\begin{proof}
By hypothesis $\Theta\cong\PP^1$.  Arguing as in the proof of~\Ref{prp}{quandoconica} and invoking~\Ref{lmm}{zerozerono}  we get that
\begin{equation}\label{giucas}
\cE_{\Theta}\cong
\begin{cases}
\text{$\cO_{\PP^1}\oplus\cO_{\PP^1}(-1)\oplus\cO_{\PP^1}(-2)$, or} \\
\cO_{\PP^1}(-1)^3.
\end{cases}
\end{equation}
Suppose that  the first isomorphism holds. Then $R_{\Theta}$ is isomorphic  either to $T_{\Theta}$  or to  a projection of 
$T_{\Theta}$. If the former holds then $\Theta$ is of Type $\cE_2$. 
Suppose that the latter holds; we will reach a contradiction. In fact 
the trivial addend in~\eqref{giucas} gives that $\Theta\subset J_{v_0}$ for some $[v_0]\in\PP(V)$. Thus we have an embedding
\begin{equation}
\iota\colon\Theta\hra \Gr(2,V/[v_0]).
\end{equation}
We have $\dim\la T_{\Theta}\ra=5$ and by assumption $R_{\Theta}$ is isomorphic to a projection of $T_{\Theta}$; thus  $\dim\la R_{\Theta}\ra=4$ i.e.~there exists $U\subset(V/[v_0])$ with $\dim U=4$ such that
\begin{equation}
 \iota(\Theta)\subset \Gr(2,U)\subset\PP(\bigwedge^2 U).
\end{equation}
By hypothesis $\Theta\cong\iota(\Theta)$ is a cubic rational normal curve and hence $\dim\la\iota(\Theta)\ra=3$;
since $\Gr(2,U)$ is a quadric hypersurface in $\PP(\bigwedge^2 U)$ we get that $\la\iota(\Theta)\ra\cap \Gr(2,U)$ has pure dimension $2$. It follows that  $\la\Theta\ra\cap \Gr(3,V)$ has pure dimension $2$ as well. Thus $\Theta$ is not a component of $\la\Theta\ra\cap \Gr(3,V)$  contradicting~\Ref{dfn}{isoiso}. This proves that if  the first isomorphism of~(\ref{giucas}) holds then $\Theta$ is of Type $\cE_2$.
Now suppose that the second isomorphism of~(\ref{giucas}) holds. 
Then $R_{\Theta}$ is not isomorphic to $T_{\Theta}$ (which is $\PP^1\times\PP^2$ embedded by the Segre map)  because any two distinct planes in $T_{\Theta}$ are disjoint. Hence $R_{\Theta}$ is isomorphic to a projection of $T_{\Theta}$. Since $\dim \la T_{\Theta}\ra=5$  the center of projection is either a point or a line. If the latter holds then $\deg g_{\Theta}=3$, that contradicts~\Ref{prp}{gradoeffe}.  
Thus $R_{\Theta}$ is isomorphic to a projection of $T_{\Theta}$ with center of projection a point i.e.~$\Theta$ is of type $\cE^{\vee}_2$. 
\end{proof}
\begin{prp}\label{prp:quartnorm}
Let  $\Theta\subset\Gr(3,V)$ be  isolated isotropic and suppose that it is a    quartic rational normal curve. 
Suppose in addition that $R_{\Theta}$ is not a cone and is non-degenerate (i.e.~$\dim\la R_{\Theta} \ra=5$).
Then   $\Theta$ is  of Type $\bf Q$. 
\end{prp}
\begin{proof}
By hypothesis $\Theta\cong\PP^1$ and hence $\cE_{\Theta}\cong\cO_{\PP^1}(-a_1)\oplus\cO_{\PP^1}(-a_2)\oplus\cO_{\PP^1}(-a_3)$ where $0\le a_i$ and $\sum_i a_i=4$. Since $R_{\Theta}$ is not a cone $a_i>0$ for $i=1,2,3$. Thus
\begin{equation}
\cE_{\Theta}\cong\cO_{\PP^1}(-1)^2\oplus\cO_{\PP^1}(-2).
\end{equation}
Since $R_{\Theta}$ is non-degenerate it is isomorphic to the projection of $T_{\Theta}$ from a point $p$. In order to prove the proposition it remains to show that  
\begin{equation}\label{sotteso}
p\in\la h_{\Theta}(\PP(\cO_{\PP^1}(-1)^2))\ra.
\end{equation}
One verifies easily that if $p\notin\la h_{\Theta}(\PP(\cO_{\PP^1}(-1)^2))\ra$ then the projections of two generic planes in $T_{\Theta}$ are disjoint - thus~(\ref{sotteso}) holds. 
\end{proof}
\begin{rmk}
Suppose that $\Theta$ is 
of Type $\bf Q$. The proof of~\Ref{prp}{quartnorm} provides the following description of $R_{\Theta}$. There exists disjoint planes $P_1,P_2\subset\PP(V)$ and embeddings $\iota_1\colon \Theta\hra P^{\vee}_1$, $\iota_2\colon \Theta\hra P_2$ with $\im(\iota_1)$ and $\im(\iota_2)$ a conic such that
\begin{equation}
R_{\Theta}=\bigcup_{x\in\Theta}\la \iota_1(x),\iota_2(x)\ra.
\end{equation}
(Of course $\la \iota_1(x),\iota_2(x)\ra$ is the plane corresponding to $x$.) 
\end{rmk}
The following result shows that there are (at least) three interesting constructions of a $\Theta$ of Type $\bf R$.  
\begin{clm}\label{clm:bellacoinc}
Let $U$ be a $4$-dimensional vector space and $V:=\bigwedge^2 U$. Let $i_{+}$ be as in~(\ref{piumenomap}). 
Let   $\Theta\subset\Gr(3,V)$ be given by $\Theta=i_{+}(C)$ where $C\subset\PP(U)$ is a rational normal cubic.  There exist an isomorphism $V=\Sym^2 L$ where $\dim L=3$ and  conics $C\subset\PP(L)$,  $C'\subset\PP(L^{\vee})$ such that $\Theta=k(C)=h(C')$. (Here $k,h$ are given by~(\ref{kappacca}).)
\end{clm}
\begin{proof}
We have a map
\begin{equation}
\begin{matrix}
C^{(2)} & \overset{f}{\lra} & \Gr(2,U)\subset\PP(V) \\
P+Q & \mapsto & \la P+Q\ra
\end{matrix}
\end{equation}
 Of course $C^{(2)}\cong \PP^2$ and one may identify  $f$ (up to projectivities) with the natural map  $\PP^2\to |\cO_{\PP^2}(2)|^{\vee}$; thus $\im(f)$ is the Veronese surface $\cV$. Let $P\in C$; the plane $i_{+}(p)$ intersects $\cV$ in a conic. Thus $i_{+}(R)\subset C(\cV)$. It follows that there exist an isomorphism $V=\Sym^2 L$ where $\dim L=3$ and a conic  $C'\subset\PP(L^{\vee})$ such that $\Theta=h(C')$. The analogous result with $h$ replaced by $k$ follows by duality - see~(\ref{dualeapiu})-(\ref{dualeakappa}).
\end{proof}
\begin{lmm}\label{lmm:quattrosez}
Let $C$ be an irreducible projective curve with $\omega_C\cong \cO_C$ (an elliptic curve, possibly singular). Let $\cF$ be  a  rank-$2$ vector-bundle  on $C$ such that
\begin{itemize}
\item[(a)]
$\deg\cF=4$,
\item[(b)]
$\cF$ is globally generated,
\item[(c)]
 there is no splitting $\cF\cong\cO_C\oplus\cL$.
\end{itemize}
Then $h^0(\cF)=4$.
\end{lmm}
\begin{proof}
By Riemann-Roch we have $\chi(\cF)=4$ hence it suffices to prove that $h^1(\cF)=0$. Suppose that $h^1(\cF)>0$; we will reach a contradiction. By Serre duality we get that $h^0(\cF^{\vee})>0$ and hence there exists a non-zero $\phi\colon\cF\to\cO_C$. Since $\cF$ is globally generated $\im(\phi)$  is globally generated; it follows that $\im(\phi)=\cO_C$. Let $\cK:=\ker(\phi)$; thus we have an  exact sequence
\begin{equation}\label{gattosacco}
0\lra\cK\lra\cF\lra\cO_C\lra 0.
\end{equation}
 By Serre duality $h^1(\cK)=h^0(\cK^{\vee})$  and furthermore $h^0(\cK^{\vee})=0$ because  $\cK$ is an invertible sheaf of degree $4$; thus $h^1(\cK)=0$ and hence~(\ref{gattosacco}) splits. That contradicts Item~(c).
\end{proof}
\begin{prp}\label{prp:noquartica}
Let  $\Theta\subset\Gr(3,V)$ be an isolated isotropic irreducible curve.  Then Items~($\alpha$),($\beta$) below cannot both hold:
\begin{itemize}
\item[($\alpha$)]
$\Theta\subset J_{v_0}$ for some $[v_0]\in\PP(V)$,
\item[($\beta$)]
$\dim\la\Theta\ra=3$ and $\Theta$ is the intersection of two quadric surfaces in $\la\Theta\ra$. 
\end{itemize}
\end{prp}
\begin{proof}
Let $V_0\in\Gr(5,V)$ be transversal to $[v_0]$. Let $C:= \ov{\rho}_{v_0}(\Theta)$; then $\ov{\rho}_{v_0}$ gives an isomorphism
$g\colon \Theta  \overset{\sim}{\lra} C$.
Let $\cF^{\vee}$ be the restriction to $C$ of the tautological rank-$2$ vector-bundle on $\Gr(2,V_0)$. 
We  have an isomorphism $\cE_{\Theta}\cong\cO_{\Theta}\oplus g^{*}\cF^{\vee}$. It follows from~\Ref{lmm}{zerozerono} that there is no splitting $\cF\cong\cO_C\oplus\cL$. Furthermore  $\deg\cF=4$  by Item~($\beta$) and of course $\cF$ is globally generated. By Item~($\beta$) we have  $\omega_C\cong\cO_C$. Thus~\Ref{lmm}{quattrosez} gives that $h^0(\cF)=4$. It follows that there exists $U\in\Gr(4,V_0)$ such that $C\subset\Gr(2,U)$.
Since $\Gr(2,U)$ is a a smooth quadric in $\PP(\bigwedge^2 U)$ we get that $\la C\ra\cap \Gr(2,U)$ has pure dimension $2$. It follows that $\la\Theta\ra\cap \Gr(3,V)$ has pure dimension $2$ contradicting~\Ref{dfn}{isoiso}.  
\end{proof}
\n
{\it Proof of~\Ref{thm}{tetadimuno}.\/}
Suppose that $\dim\la\Theta\ra=1$; then $\Theta$ is of Type~$
\cF_1$. Thus we may assume that $2\le\dim\la\Theta\ra$. 
By definition $\Theta$ is an irreducible component of $\la\Theta\ra\cap\Gr(3,V)$. Since $\Gr(3,V)$ is cut out by Pl\"ucker quadrics (in $\PP(\bigwedge^3 V)$) it follows that:
\begin{itemize}
\item[(i)]
If $\dim\la\Theta\ra=2$ then $\Theta$ is a smooth conic.
\item[(ii)]
If $\dim\la\Theta\ra=3$ then $\Theta$ is either a cubic rational normal curve or the complete intersection of $2$ quadrics. 
\end{itemize}
If (i) holds then  $\Theta$ is of Type~$\cD$  by~\Ref{prp}{quandoconica}. Suppose that (ii) holds: if $\Theta$ is  a cubic rational normal curve then it is of Type~$\cE_2$ or of Type~$\cE_2^{\vee}$  by~\Ref{prp}{cubgobba}. Thus from now on we may assume that 
\begin{itemize}
\item[(I)]
$\dim\la\Theta\ra=3$ and $\Theta$ is the complete intersection of two quadrics, or
\item[(II)]
$4\le\dim\la\Theta\ra$. 
\end{itemize}
By Morin one of~(a) - (e) of~\Ref{thm}{teomorin} holds. We will perform a case-by-case analysis.
\vskip 2mm
\n
(a): $\Theta\subset F_{\pm}(\cQ)$.
  Let $U$ be a $4$-dimensional complex vector-space and identify $V$ with $\bigwedge^2 U$ so that $\cQ$ gets identified with $\Gr(2,U)$.  We may assume that  $\Theta\subset F_{+}(\cQ)$; thus   $\Theta:=i_{+}(C)$ for an irreducible curve $C\subset\PP(U)$. By definition $\Theta$ is an irreducible component   of $\la\Theta\ra\cap F_{+}(\cQ)$. 
Since $i_{+}$ is given by the complete linear system of quadrics in $\PP(U)$ we get that $C$ is a component of a complete intersection of quadrics. Thus $C$ is a rational curve of degree at most $3$ or the complete intersection of two quadrics. By (I)-(II) above we have  $3\le \dim\la i_{+}(C)\ra$ and hence $C$ is not a line.  Suppose that $C$ is a conic; as is easily verified $R_{\Theta}$ is not a cone, and by duality we get that it is non-degenerate as well. Since $\Theta$ is a degree-$4$ rational normal curve we get that it is of Type~$\bf Q$ by~\Ref{prp}{quartnorm}. If $C$ is a cubic rational normal curve then $\Theta$ is of Type~$\bf R$. Lastly if $C$ is the complete intersection of two quadrics then $\Theta$ is of Type~$\bf S$. 
\vskip 2mm
\n
(b): $\Theta\subset C(\cV)$ or $\Theta\subset T(\cV)$. 
 Let $L$ be a $3$-dimensional complex vector-space and identify $V$ with $\Sym^2 L$ so that $\cV$ gets identified with $\PP((\Sym^2 L)_1)$. We discuss the case  $\Theta\subset C(\cV)$, the other case will follow by duality. There exists a curve $C\subset \PP(L^{\vee})$ such that $\Theta=h(C)$. We recall that $h$ is identified (up to projectivities) with $|\cO_{\PP(L^{\vee})}(3)|$. Arguing as in Case~(a) we get that $\deg C\le 3$. If $C$ is a line then $\Theta$ is a cubic rational normal curve  and hence it is of Type~$\cE_2$ or of Type~$\cE_2^{\vee}$ by~\Ref{prp}{cubgobba} - in fact  of Type~$\cE_2$. If $C$ is a smooth conic then $\Theta$ is of Type~$\bf R$  by~\Ref{clm}{bellacoinc}. If $C$ is a cubic then $\Theta$ if of Type~${\bf T}^{\vee}$.
\vskip 2mm
\n
(c): $\Theta\subset J_{v_0}$.   
By assumption one of  Items~(I), (II) above holds. In fact Item~(I) cannot hold by~\Ref{prp}{noquartica}. 
Hence Item~(II) holds. We claim that  
\begin{equation}\label{eluana}
\dim\la\Theta \ra=4. 
\end{equation}
In fact if $4<\dim\la\Theta \ra$ then every irreducible component of $\la\Theta\ra\cap \Gr(3,V)$ has dimension at least $2$, that contradicts~\Ref{dfn}{isoiso}. By~(\ref{eluana}) we get that
$\Theta$ is of Type~$\cA$. 
\vskip 2mm
\n
(d): $\Theta\subset \Gr(3,E)$ where $E\in\Gr(5,V)$.    
Then $\wt{\delta}_V(\Theta)\subset J_{\phi}$ where $\la\phi\ra=\Ann  (E)$; by the previous case we get that $\Theta$ is of Type~$\cA^{\vee}$. 
\vskip 2mm
\n
(e): $\Theta\subset I_U$ where $U\in\Gr(3,V)$.
Suppose first that $U\in\Theta$. Then $\Theta$ is an irreducible component  of  $\la\Theta\ra\cap I_U$, and since the latter  is a cone with vertex $U$ it follows that $\Theta$ is a cone with vertex $U$.  Thus $\Theta$ is a line and hence it is of Type $\cF_1$. From now on we may assume that $U\notin\Theta$; since $\Theta$ is an irreducible component of $\la\Theta\ra\cap\Gr(3,V)$ it follows that
\begin{equation}\label{novertice}
U\notin\la\Theta\ra.
\end{equation}
Let $\ov{\rho}_U$ be the (rational) map of~\eqref{robaruma}; by~\eqref{novertice}  the restriction of $\ov{\rho}_U$ to $\Theta$ is a regular isomorphism onto  
\begin{equation}\label{nelsegre}
C:=\ov{\rho}_U(\Theta)\subset\PP(U)\times\PP(V/U)\subset
\PP((U)\otimes\PP(V/U)). 
\end{equation}
By assumption one of~(I), (II) above holds. We claim that~(I) cannot hold. In  fact let $f\colon C\to\PP(U)$ and $g\colon C\to\PP(V/U)$ be the two projections. One easily checks that neither $f$ nor $g$ is constant. We have 
\begin{equation}
\deg f^{*}\cO_{\PP(U)}(1)+\deg f^{*}\cO_{\PP(V/U)}(1)=4.
\end{equation}
Since $C$ has arithmetic genus $1$ we get that  
\begin{equation}
2=\deg f^{*}\cO_{\PP(U)}(1)=\deg f^{*}\cO_{\PP(V/U)}(1)
\end{equation}
and moreover $\im(f),\im(g)$ are lines, say $\im(f)=\PP(U_2)$ and $\im(g)=\PP(W_2)$. It follows that 
\begin{equation}
\la C\ra\supset(\PP(U_2)\times \PP(W_2))
\end{equation}
where $\la C\ra$ is the span of $C$ in $\PP((U)\otimes\PP(V/U))$. Thus $C$ is not an irreducible component of $\la C\ra\cap(\PP(U)\times\PP(V/U))$ and hence $\Theta$ is not an irreducible component of $\la\Theta\ra\cap\Gr(3,V)$; that  contradicts~\Ref{dfn}{isoiso}.   Hence Item~(II) above holds i.e.
\begin{equation}\label{sopraquattro}
4\le\dim\la\Theta\ra. 
\end{equation}
On the other hand 
\begin{equation}\label{sottocinque}
\dim\la\Theta\ra\le 5.
\end{equation}
 In fact suppose that $5<\dim\la\Theta\ra$. By~(\ref{novertice}) we get that $5<\dim\la C\ra$ and hence  every irreducible component of $\la C\ra\cap(\PP(U)\times\PP(V/U)$ has dimension at least $2$. It follows that  $\Theta$ is not an irreducible component of $\la\Theta\ra\cap\Gr(3,V)$,  contradiction. 
This proves~(\ref{sottocinque}).
Assume that $\dim\la\Theta\ra=4$. Then 
\begin{equation}\label{kalamazoo}
\dim\la C\ra=4
\end{equation}
by~(\ref{novertice}). Since $\deg(\PP(U)\times\PP(V/U))=6$ it follows that $4\le\deg C\le 5$. If $\deg C=4$ then $C$ is a quartic rational normal curve; as is easily verified $R_{\Theta}$ is not a cone and is non-degenerate thus $\Theta$ is of Type~$\bf Q$ by~\Ref{prp}{quartnorm}. If $\deg C=5$ we will reach a contradiction. First let's prove that $C$ is of arithmetic genus $1$. In fact the intersection of $\PP(U)\times\PP(V/U)$ with a generic $5$-dimensional projective space containing $C$ is a curve of degree $6$ and arithmetic genus $1$ and the component different from $C$ is a line meeting $C$ in a single point (and not tangent to $C$); it follows that $p_a(C)=1$. Arguing as for $C$ satisfying Item~(I) we get that 
\begin{equation}\label{citiciti}
C\subset(\PP(U_2)\times \PP(V/U)),\qquad
U_2\in\Gr(2,U)
\end{equation}
 or 
\begin{equation}\label{bangbang}
C\subset(\PP(U)\times \PP(W_2)),\qquad
W_2\in\Gr(2,V/U).
\end{equation}
 Suppose that~(\ref{citiciti}) holds; since $\dim\la(\PP(U_2)\times \PP(V/U))\ra=5$ we get that 
\begin{equation}
2\le \dim\la C\ra\cap(\PP(U_2)\times\PP(V/U))
\end{equation}
by~(\ref{kalamazoo}). It follows that  $\Theta$ is not an irreducible component of $\la\Theta\ra\cap\Gr(3,V)$,  contradiction. 
If ~(\ref{bangbang}) holds we argue similarly and again we get a contradiction. Thus we are left with the case $\dim\la\Theta\ra=5$; then $\Theta$ is of Type~$\cC_2$. 
\qed
\begin{dfn}
Let $X$ be one of the types listed in Table~\eqref{catundim}: we let $\BB_X\subset\lagr$ be the closure of the set of $A$ such that $\Theta_A$ contains an irreducible component of Type $X$.
\end{dfn}
\begin{prp}\label{atetauno}
Let $A\in\lagr$ and suppose that there exists a $1$-dimensional irreducible component of $\Theta_A$. There exists a Type $X$ in Table~\eqref{catundim} such that $A\in\BB_X$.
\end{prp}
\begin{proof}
Let  $\Theta\subset\Theta_A$ be a $1$-dimensional irreducible component. By~\Ref{rmk}{compisola} we know that $\Theta\subset\Gr(3,V)$ is an isolated isotropic irreducible curve and hence the proposition follows from~\Ref{thm}{tetadimuno}. 
\end{proof}
\begin{prp}\label{prp:genlag}
 Let  $X$ be one of the  Types appearing in Table~\eqref{catundim} and $\Theta$ be generic of Type $X$ (this makes sense because the relevant parameter spaces are irreducible). There exists $A\in\lagr$ such   that $\Theta_A=\Theta$ and moreover $\Theta_A$ is generically reduced. In particular $\BB_X\not=\es$.
\end{prp}
\begin{proof}
By~\Ref{clm}{ducalberti} we may assume that 
\begin{equation}\label{nullaltro}
\la\Theta\ra\cap\Gr(3,V)=\Theta 
\end{equation}
and the scheme-theoretic intersection is reduced. Let $L:=\la\la\Theta\ra\ra$ and $\ell:=\dim L$. We have a bijection
\begin{equation}\label{parigi}
\begin{matrix}
\{A\in\lagr\mid A\supset\la\la\Theta\ra\ra\} & \overset{\sim}{\lra} & \lag(L^{\bot}/L) \\
A & \mapsto & A/L
\end{matrix}
\end{equation}
We will show that there exists $B\in\lag(L^{\bot}/L)$ such that the corresponding $A\in\lagr$ has the stated properties. Let
\begin{equation*}
Z:=\{W_0\in\Gr(3,V)\mid W_0\cap W\not=\{0\}\}.
\end{equation*}
If $A\in\lagr$ contains $\la\la\Theta\ra\ra$ and $W_0\in\Theta_A$ then $W_0\in Z$. On the other hand let $W_0\in)Z\setminus\Theta)$. By~\eqref{nullaltro} we have $\bigwedge^3 W_0\notin L$ and hence its class is non-zero in $\PP(L^{\bot}/L)$; it follows that
\begin{equation}\label{robot}
\cod(\{B\in\lag(L^{\bot}/L)\mid B\supset \bigwedge^3 W_0\},\lag(L^{\bot}/L))=10-\ell.
\end{equation}
(We abuse notation: $\bigwedge^3 W_0$ is actually the image of  $\bigwedge^3 W_0$ in $L^{\bot}/L$.)   Let $\varphi\colon (Z\setminus\Theta)\to \PP(L^{\bot}/L)$ which maps $W_0$ to the image of $\bigwedge^3 W_0$  in $L^{\bot}/L$.  By~\eqref{robot} it suffices to prove that
\begin{equation}\label{megatron}
\dim\varphi (Z\setminus\Theta)<(10-\ell).
\end{equation}
Let $Z=Z_1\cup\ldots\cup Z_r$ be the decomposition into irreducible components. We must prove that
\begin{equation}\label{primeoptimus}
\dim\varphi (Z_i\setminus\Theta)<(10-\ell)
\end{equation}
for all $i$. For each Type $X$ and for $\Theta$ generic of that Type  we will describe the irreducible components $Z_i$ and we will check that~\eqref{primeoptimus} holds. 
\vskip 2mm
\n
$\boxed{\text{\it $\Theta$ of Type $\cF_1$}}$ Notice that $R_{\Theta}$ is a $3$-dimensional linear space. Let $M\subset\PP(V)$ be the intersection of all $\PP(W)$ for $W\in\Theta$, thus $M$ is a line. The decomposition into irreducibles of $Z$ is the following:
\begin{equation*}
Z=\{W_0\mid \dim\PP(W_0)\cap R_{\Theta}\ge 1\}\cup \{W_0\mid \PP(W_0)\cap M\not= \es\}.
\end{equation*}
Let  $Z_i$ be an irreducible component; then  $\dim Z_i=7$ and hence $\dim\varphi(Z_i)\le 7$. Since $\ell=2$ we get that~\eqref{primeoptimus}  holds.
\vskip 2mm
\n
$\boxed{\text{\it $\Theta$ of Type $\cD$}}$ In this case $R_{\Theta}$ is a $3$-dimensional quadric with one singular point $[v_0]$. The variety $F_1(R_{\Theta})$ parametrizing lines on  $R_{\Theta}$ has two connected components (see~\Ref{clm}{fanqusing} for a detailed description),call them $F_1(R_{\Theta})^{\pm}$. The lines contained in the planes parametrized by $\Theta$  belong  to one of the two components, say $F_1(R_{\Theta})^{+}$. The decomposition into irreducibles of $Z$ is the following:
\begin{equation*}
Z=J_{v_0}\cup
\Gr(2, R_{\Theta})
\cup\{W_0\mid \text{$\PP(W_0)\cap R_{\Theta}$ contains a line $J\in F_1(R_{\Theta})^{-}$}\}.
\end{equation*}
Let  $Z_i$ be an irreducible component; then  $\dim Z_i=6$ and hence $\dim\varphi(Z_i)\le 6$. Since $\ell=3$ we get that~\eqref{primeoptimus}  holds. 
\vskip 2mm
\n
$\boxed{\text{\it $\Theta$ of Type $\cE_2$ or $\cE^{\vee}_2$}}$ Suppose first that $\Theta$ is of Type $\cE_2$. Notice that $R_{\Theta}$ is a cone over a smooth normal rational cubic scroll in a $4$-dimensional linear space. Let $[v_0]\in R_{\Theta}$ be the vertex and $\PP(U)\subset R_{\Theta}$ be the plane joining $[v_0]$ to the $(-1)$-line of the cubic scroll. The decomposition into irreducibles of $Z$ is $Z=Z_1\cup I_U\cup J_{v_0}$  where the generic plane in $Z_1$ intersects $R_{\Theta}$ in a smooth conic. The generic plane in $Z_1$ corresponds to an injection $\cO_{\Theta}(-2)\hra\cE_{\Theta}$ and hence $\dim Z_1=\dim\PP(H^0(\cE_{\Theta}(2)))=5$. We also have $\dim J_{v_0}=5$; since $\ell=4$ we get that~\eqref{primeoptimus} holds for $Z_1$ and for $I_U$. Lastly consider  $J_{v_0}$ which has dimension $6$. We notice that $F_{v_0}\supset L$ and that  $\varphi(J_{v_0}\setminus\Theta)\subset\PP(F_{v_0}/L)$; since  $\dim\PP(F_{v_0}/L)=5$ we see that~\eqref{primeoptimus} holds in this case as well. If $\Theta$ is of Type $\cE^{\vee}_2$  the result follows by duality from the case when $\Theta$ is of Type $\cE_2$. 
\vskip 2mm
\n
$\boxed{\text{\it $\Theta$ of Type $\bf Q$}}$ We may choose an isomorphism $V\cong\bigwedge^2 U$ where $\dim U=4$ and a conic $C\subset\PP(U)$ such that $\Theta=i_{+}(C)$. Recall that we have an immersion $i_{-}\colon \PP(U^{\vee})\hra \Gr(3,V)$. Every plane parametrized by $\Theta$ intersects the plane $i_{-}(\la C\ra)$. Let $i_{-}(\la C\ra)=\PP(H)$. 
The decomposition into irreducibles of $Z$ is $Z=Z_1\cup I_H\cup \Theta_{A_{+}(U)}$
  where the generic plane  in $Z_1$ is spanned by the images via $i_{+}$ of the lines in one of the two rulings of a smooth quadric $Q\in|\cI_C(2)|$.  We have $\dim Z_1=4$; since $\ell=5$ we get that~\eqref{primeoptimus} holds for $Z_1$. On the other hand $S_{H}\supset L$, $A_{+}(U)\supset L$ and we have that  $\varphi(I_{H}\setminus\Theta)\subset\PP(S_H/L)$, $\varphi(\Theta_{A_{+}(U)}\setminus\Theta)\subset\PP(A_{+}(U)/L)$; since $4=\dim\PP(S_{H}/L)=\dim\PP(A_{+}(U)/L)$ we see that~\eqref{primeoptimus} holds in this case as well.
\vskip 2mm
\n
$\boxed{\text{\it $\Theta$ of Type $\cA$ or $\cA^{\vee}$}}$ By duality it suffices to consider $\Theta$ of Type $\cA$. Then $R_{\Theta}$ is a cone with vertex $[v_0]$ over a degree-$5$ surface $\ov{R}_{\Theta}$ ruled over an elliptic curve and spanning a $4$-dimensional linear subspace. Clearly $J_{v_0}\subset Z$, we analyze $\Lambda\in(Z\setminus J_{v_0})$. Let $\cO_{\Theta}\subset\cE_{\Theta}$ be the sub line-bundle corresponding to the vertex $[v_0]$ and $\ov{\cE}_{\Theta}:=\cE_{\Theta}/\cO_{\Theta}$; by genericity of $\Theta$ we may assume that $\ov{\cE}_{\Theta}$ is a stable (rank-$2$) vector-bundle. The projection of $\Lambda$ from $[v_0]$ is a plane $\ov{\Lambda}\subset\la \ov{R}_{\Theta}\ra$ intersecting each line of the ruling of $\ov{R}_{\Theta}$. Thus $\ov{\Lambda}$ defines a section of $\PP(\ov{\cE}_{\Theta})$ i.e.~a sub line-bundle 
\begin{equation}\label{puntura}
\cL\hra\ov{\cE}_{\Theta}. 
\end{equation}
By stability of  $\ov{\cE}_{\Theta}$ we have $\deg\cL\le -3$. On the other hand $-4\le\deg\cL$ because $\deg\ov{R}_{\Theta}=5$. We claim that we cannot have $\deg\cL=-4$. In fact  the  transpose of~\eqref{puntura} is a surjection $\ov{\cE}^{\vee}_{\Theta}\to \deg\cL^{\vee}$; since the map on global sections is surjective it follows that the section corresponding to~\eqref{puntura} is a degree-$4$ elliptic curve spanning a $3$-dimensional space, {\bf not} a plane. Thus $\deg\cL=-3$; conversely to each~\eqref{puntura} with $\deg\cL=-3$ there corresponds a plane  $\ov{\Lambda}\subset\la \ov{R}_{\Theta}\ra$ intersecting each line of the ruling of $\ov{R}_{\Theta}$. Given such a plane   all the planes in $\la[v_0]\cup\ov{\Lambda}\ra$ belong to $Z$. Let $Z_1\subset\Gr(3,V)$ be the closure of the locus of such planes; we have proved that the irreducible decomposition of $Z$ is $Z=Z_1\cup J_{v_0}$.
By stability of $\ov{\cE}_{\Theta}$ we have $\dim \Hom(\cL,\ov{\cE}_{\Theta})=1$ for every line-bundle $\cL$ of degree $-3$; it follows that $\dim Z_1=4$. 
  Since $\ell=5$ we get that~\eqref{primeoptimus} holds for $(Z_1\setminus\Theta)$. The argument for the component $J_{v_0}$ is the same as that given for $\Theta$ of Type $\cE_2$.  
\vskip 2mm
\n
$\boxed{\text{\it $\Theta$ of Type $\cC_2$}}$ 
Choose an isomorphism $V\cong\bigwedge^2 \CC^4$. Let $H\subset\CC^4$ be a subspace of codimension $1$ and $C\subset\PP(H)$ a cubic curve. Let $\Theta:=i_{+}(C)$. Let $U\in\Gr(3,V)$ be the subspace such that $i_{-}(H)=\PP(U)$.  One checks easily that $\Theta\in (I_U\setminus \{U\})$ and that $\dim\la\Theta\ra=5$. The proof of~\Ref{thm}{tetadimuno} gives that $\Theta$ is of Type $\cC_2$ - see Case~(e). We will prove that there exists $A\in\lagr$ such   that $\Theta_A=\Theta$; it will follow that the same is true for a  generic $\Theta$ of Type $\cC_2$ (actually the generic $\Theta$ of Type $\cC_2$  is equal to $i_{+}(C)$ as above).  The decomposition into irreducibles of $Z$ is $Z=I_U\cup A_{+}(\CC^4)$; it follows that~\eqref{megatron} holds. 
\vskip 2mm
\n
$\boxed{\text{\it $\Theta$ of Type $\bf R$ or $\bf S$}}$ The decomposition into irreducibles of $Z$ is $Z=Z_1\cup \Theta_{A_{+}(U)}$
  where the generic plane  in $Z_1$ is spanned by the images via $i_{+}$ of the lines in one of the two rulings of a smooth quadric $Q\in|\cI_C(2)|$.  Suppose that $\Theta$ is of Type $\bf R$. Then $\dim Z_1=2$; since $\ell=5$ we get that~\eqref{primeoptimus} holds for $Z_1$. One deals with the component $\Theta_{A_{+}(U)}$ as usual; that proves that~\eqref{megatron} holds for $\Theta$  of Type $\bf R$. Suppose now that $\Theta$ is of Type $\bf S$. Then $\dim Z_1=1$; since $\ell=8$ we get that~\eqref{primeoptimus} holds for $Z_1$ and~\eqref{megatron} follows.
\vskip 2mm
\n
$\boxed{\text{\it $\Theta$ of Type $\bf T$ or ${\bf T}^{\vee}$}}$ By duality it suffices to consider $\Theta$ of Type ${\bf T}^{\vee}$. We have a rational map $\varphi\colon\Gr(3,V)\dashrightarrow \PP(Sym^3 L^{\vee})$ which assigns to a $3$-dimensional subspace $\Lambda\subset \Sym^2 L$ the set of singular points of singular non-zero quadrics $q\in\Lambda$. Since $\dim \Gr(3,V)=\dim \PP(Sym^3 L^{\vee})$ either $\varphi$ is not dominant or  it is dominant with finite generic fiber (in fact it is dominant but we do not need this). Let $\Theta=h(C)$ where $C\subset\PP(L^{\vee})$ is a generic cubic. Then $Z=Z_1\cup \Theta_{A_h(L)}$ where $Z_1=\varphi^{-1}(C)$. Since $Z_1$ is finite and $\ell=9$ we get that~\eqref{primeoptimus} holds for $Z_1$. One deals with the component $\Theta_{A_h(L)}$ as usual.  
\vskip 2mm
\n
It remains to prove that $\Theta_A=\Theta$ is generically reduced.  Let $L=\la\la\Theta\ra\ra$. Given $[W]\in\Theta$ we know that $\dim(L\cap S_W)=2$ because by~\Ref{clm}{ducalberti} the intersection $\PP(L)\cap\Gr(3,V)$ is smooth $1$-dimensional. We must  prove that if $A\in\lagr$ is generic  in the left-hand side of~\eqref{parigi} then 
\begin{equation}\label{valbene}
A\cap S_W=L\cap S_W.
\end{equation}
Since $S_W$ is lagrangian for $(,)_V$ the symplectic form defines an isomorphism $\bigwedge^3 V/S_W\overset{\sim}{\lra} S_W^{\vee}$. Thus $(,)_V$  gives an injection $L/(L\cap S_W)\hra S_W^{\vee}$. It follows that $L^{\bot}\cap S_W/(L\cap S_W)$ is a lagrangian subspace of $L^{\bot}/L$ and hence the generic $B\in\lag(L^{\bot}/L)$ intersects trivially $L^{\bot}\cap S_W/(L\cap S_W)$; the corresponding $A\in\lagr$   in the left-hand side of~\eqref{parigi}  satisfies~\eqref{valbene}.
\end{proof}
A straightforward parameter count gives the dimensions of the $\BB_X$'s;   we listed their codimensions in Table~\eqref{dimbix} (since $\delta_V$ preserves dimension we omitted writing out the codimension of $\BB_{\cE_2^{\vee}},\BB_{\cA^{\vee}},\BB_{{\bf T}^{\vee}}$). 
\begin{crl}\label{crl:sonocomp}
Let $X$ be one of the types listed in Table~\eqref{catundim}. Then  $\BB_X$ is an irreducible component of $\Sigma_{\infty}$. Moreover the $\BB_X$'s are pairwise distinct.
\end{crl}
\begin{proof}
Irreducibility of $\BB_X$ follows from irreducibility of the parameter space for curves of Type $X$. 
By~\Ref{prp}{genlag} we get that $\BB_X$ is an irreducible component of $\Sigma_{\infty}$. It remains to prove that if $X_1\not=X_2$ then $\BB_{X_1}\not\subset\BB_{X_2}$. Suppose that $\BB_{X_1}\subset\BB_{X_2}$. Since for $A$  generic in $\BB_{X_i}$ the scheme $\Theta_A$ is a generically reduced curve we get that $\deg \Theta_1=\Theta_2$ for $\Theta_i$ generic of Type  $X_i$. Looking at the degrees of $\Theta$ appearing in Table~\eqref{catundim} and the codimensions of $\BB_X$'s in Table~\eqref{dimbix} we conclude that the inclusion in question is one of the followings: $\BB_{\cE_2}=\BB_{\cE^{\vee}_2}$, $\BB_{\cA}=\BB_{\cA^{\vee}}$, $\BB_{\bf R}\subset\BB_{\cC_2}$ or $\BB_{\bf T}=\BB_{{\bf T}^{\vee}}$. One proves quickly that the  listed inclusions do not hold except possibly the last one. Suppose that  $\BB_{\bf T}=\BB_{{\bf T}^{\vee}}$. Then the following holds: for $\Theta$ generic of Type $\bf T$ i.e.~such that  $\Theta=k(C)$ for an isomorphism $V\cong \Sym^2 L$ and a cubic $C\subset\PP(L)$ there exist an isomorphism $V\cong \Sym^2 L^{\vee}$ and a cubic $C'\subset\PP(L^{\vee})$ such that $\Theta =h(C')$. Thus $R_{k(C)}=R_{h(C')}$; this is absurd because the closure of the set of multibranch points of $R_{k(C)}$ is isomorphic to  $C^{(2)}$ (points $l_0\cdot l_1$ where $l_0,l_1\in C$) while the closure of the set of multibranch points of $R_{h(C)}$ is isomorphic to  $\PP^{2}$ (points $l^2$). 
\end{proof}
Now assume that $X$ is of calligraphic Type, denote it by $\cX$. We will show that one may characterize the generic point of  $\BB_{\cX}$  by a certain flag condition that one encounters when studying GIT-stability.  Let 
 \begin{equation}\label{basedivu}
 \sF:=\{v_0,\ldots,v_5\} 
\end{equation}
 be a basis of $V$. For each calligraphic $\cX$ appearing in Table~\eqref{catundim} we define $\BB^{\sF}_{\cX}$ to be the set of $A\in\lagr$ satisfying the condition appearing on the second  column of the corresponding row of  Table~\eqref{comericci}; we adopt the notation
\begin{equation}\label{vuconij}
V_{ij}:=\la v_i,v_{i+1},\ldots,v_j\ra,\qquad 0\le i<j\le 5.
\end{equation}
\begin{table}[tbp]
\caption{Flag conditions, I}\label{comericci}
\vskip 1mm
\centering
\renewcommand{\arraystretch}{1.60}
\begin{tabular}{ll}
\toprule
   name  & flag condition    \\
\midrule
 $\BB^{\sF}_{\cA}$ &
$\dim A\cap ([v_0]\wedge\bigwedge  ^2 V_{15})\ge 5$       \\
\midrule
 $\BB^{\sF}_{\cA^{\vee}}$  & $\dim A\cap (\bigwedge  ^3 V_{04})\ge 5$     \\
\midrule
   $\BB^{\sF}_{\cC_2}$ &  $\dim A \cap (\bigwedge  ^3 V_{02}\oplus
(\bigwedge  ^2 V_{02}\wedge V_{35}))\ge 6$   \\
\midrule
 $\BB^{\sF}_{\cD}$   & $\dim A\cap([v_0]\wedge\bigwedge  ^2 V_{14})\ge 3$  \\
\midrule
  $\BB^{\sF}_{\cE_2}$ &  $\dim A\cap ([v_0]\wedge(\bigwedge  ^2 V_{12})\oplus ([v_0]\wedge V_{12}\wedge  V_{35}))\ge 4$  \\ 
\midrule
  $\BB^{\sF}_{\cE^{\vee}_2}$ &     $\dim A\cap(\bigwedge  ^3 V_{02}\oplus (\bigwedge  ^2 V_{02}\wedge V_{34}) )\ge 4$     \\   
\midrule
   $\BB^{\sF}_{\cF_1}$  & $A\supset(\bigwedge  ^2 V_{01}\wedge  V_{23})$  \\ 
\bottomrule 
\end{tabular}
\end{table} 
Let
\begin{equation*}
\BB^{*}_{\cX}:=\bigcup_{\sF}\BB^{\sF}_{\cX}
\end{equation*}
where $\sF$ runs through the set of bases of $V$.
\begin{clm}\label{clm:stelladensa}
Let $\cX$ be one  of the calligraphic Types in Table~\eqref{catundim}. Then $\BB^{*}_{\cX}$ is a constructible dense subset of $\BB_{\cX}$.
\end{clm}
\begin{proof}
It suffices to prove the following two results:
\begin{enumerate}
\item[(I)]
Let $A\in\lagr$ and suppose that $\Theta_A$ contains an irreducible component $\Theta$ of Type $\cX$ according to Table~\eqref{catundim}. There exists a basis $\sF$ of $V$ such that $A\in\BB_{\cX}^{\sF}$.
\item[(II)]
Let $\sF$ be a basis of $V$. If $A\in\BB_{\cX}^{\sF}$ is generic then  $\Theta_A$ contains an irreducible component $\Theta$ of Type $\cX$.
\end{enumerate}
Items~(I), (II) are obvious except possibly for $\cX=\cE_2$ or $\cX=\cE^{\vee}_2$. Suppose that  $\cX=\cE_2$. Let's prove (I). We have $\cE_{\Theta}\cong\cO_{\PP^1}\oplus\cO_{\PP^1}(-1)\oplus 
\cO_{\PP^1}(-2)$ and hence the Harder-Narasimhan filtration of $\cE_{\Theta}$ gives rise to a filtration 
$U_1\subset U_3\subset U_6=V$ where $\dim U_i=i$.  Let $\sF$ be a basis of $V$ such that $[v_0]=U_1$. $[v_0]\oplus V_{12}=U_3$; then $A\in \BB_{\cE_2}^{\sF}$. Let's prove (II). Given $A\in\BB_{\cE_2}^{\sF}$ we let 
\begin{equation*}
F_A:=A\cap([v_0]\wedge(\bigwedge  ^2 V_{12})\oplus ([v_0]\wedge V_{12}\wedge  V_{35})).
\end{equation*}
Suppose that  $A$ is generic. Then $\dim F_A=4$, moreover
$\Theta_A=\PP(F_A)\cap \Gr(3,V)$ and the latter is a curve of Type $\cE_2$. Next suppose that  $\cX=\cE^{\vee}_2$. Let's prove (I). By Table~\eqref{catundim} we may describe $R_{\Theta}$ as follows. Let 
\begin{equation*}
\sigma\colon \PP^1\times\PP^2\hra |\cO_{\PP^1}(1)\boxtimes \cO_{\PP^2}(1)|^{\vee}\cong\PP(V)
\end{equation*}
be Segre's embedding followed by a suitable isomorphism $\CC^2\otimes\CC^3\cong V$. Let
$\cS\subset\PP(V)$ be the image of $\sigma$. Then $R_{\Theta}$ is the projection of $\cS$ from a point $p\notin\cS$; of course the planes that sweep out $R_{\Theta}$ are the projections of the planes $\sigma(\{x\}\times\PP^2)$ for $x\in\PP^1$. Let $H\subset\PP(V)$ be the hyperplane containing $R_{\Theta}$ i.e.~the hyperplane to which we project from $p$. 
One checks easily that there exists a (unique) line $L\subset\PP^2$ such that the span $M:=\la\sigma(\PP^1\times L)\ra$ contains $p$; notice that  $\dim M=3$ because $\sigma(\PP^1\times L)$ is a smooth quadric surface.
Let $P\subset H$ be the projection of  $(M\setminus\{p\})$ from $p$; thus $P$ is a plane. It follows from the definitions that each plane in $\Theta$ intersects $P$ in  a line. Let $\sF$ be a basis of $V$ such that $\PP(V_{02})=P$ and $\PP(V_{04})=H$; then $A\in \BB_{\cE^{\vee}_2}^{\sF}$. The proof of (II) is analogous to the proof of~(II) for $\cX=\cE_2$, we omit details. 
\end{proof}
\begin{rmk}\label{rmk:eccoperche}
Let   $X$ be one of the Types appearing in Table~\eqref{catundim}. In a forthcoming paper we will show that if  $X$ is calligraphic  and $A\in\BB_X$ then $A$ is not GIT-stable (in general it is properly semistable) - calligraphic Types have been  ordered according to the complexity of the destabilizing $1$-PS for generic lagrangians of that Type.  On the other hand we will show that if  $X$ is boldface and $A\in\BB_X$ is generic  then it is GIT-stable.
\end{rmk}
\begin{table}[tbp]
\caption{Codimension of the $\BB_X$'s}\label{dimbix}
\vskip 1mm
\centering
\renewcommand{\arraystretch}{1.60}
\begin{tabular}{l l l l l l l l l l}
\toprule
 $\BB_{\star}$ & $\BB_{\cF_1}$ & $\BB_{\cD}$  & $\BB_{\cE_2}$  & $\BB_{\bf Q}$  & 
 $\BB_{\cA}$  & $\BB_{\cC_2}$  & $\BB_{\bf R}$  & $\BB_{\bf S}$  & $\BB_{\bf T}$ \\
\midrule
 $\cod(\BB_{\star},\lagr)$  & $7$ & 9   & 11  & 9  & 10  & 12  & 17  & 16  & 18  \\
\bottomrule 
\end{tabular}
\end{table} 
\subsection{Two-dimensional components of $\Theta_A$}
\setcounter{equation}{0}
We will analyze those  $A\in\lagr$ such that $\Theta_A$ contains a $2$-dimensional irreducible component.  In order to state our result we will introduce certain closed subsets of $\lagr$.
First we associate to each of a collection of  calligraphic Types $\cX$ (containing all those appearing in Table~\eqref{comericci}) a constructible subset  $\XX_{\cX,+}\subset \lagr$.    Let $\sF$  be a basis of $V$ as in~\eqref{basedivu}. We let $\XX_{\cX,+}^{\sF}$  be the set of $A\in\lagr$ satisfying the condition appearing on the second  column of the corresponding row of  Table~\eqref{nonlascio} (Notation~\eqref{vuconij} is in force).
\begin{table}[tbp]
\caption{Flag conditions, II}\label{nonlascio}
\vskip 1mm
\centering
\renewcommand{\arraystretch}{1.60}
\begin{tabular}{ll}
\toprule
   name  & flag condition    \\
\midrule
 $\XX^{\sF}_{\cA_{+}}$ &
$\dim A\cap ([v_0]\wedge\bigwedge  ^2 V_{15})\ge 6$       \\
\midrule
 $\XX^{\sF}_{\cA^{\vee}_{+}}$  & $\dim A\cap (\bigwedge  ^3 V_{04})\ge 6$     \\
\midrule
   $\XX^{\sF}_{\cC_{1,+}}$ &  $A\supset\bigwedge  ^3 V_{02}$ and  $\dim A\cap (\bigwedge  ^2 V_{02}\wedge V_{35})\ge 4$   \\
\midrule
   $\XX^{\sF}_{\cC_{2,+}}$ &  $\dim A \cap (\bigwedge  ^3 V_{02}\oplus
\bigwedge  ^2 V_{02}\wedge V_{35})\ge 7$   \\
\midrule
 $\XX^{\sF}_{\cD_{+}}$   & $\dim A\cap([v_0]\wedge\bigwedge  ^2 V_{14})\ge 4$  \\
\midrule
  $\XX^{\sF}_{\cE_{2,+}}$ &    $\dim A\cap ([v_0]\wedge\bigwedge  ^2 V_{12}\oplus [v_0]\wedge V_{12}\wedge  V_{35})\ge 5$  \\ 
\midrule
  $\XX^{\sF}_{\cE^{\vee}_{2,+}}$ &  $\dim A\cap(\bigwedge  ^3 V_{02}\oplus \bigwedge  ^2 V_{02}\wedge V_{34})\ge 5$  \\ 
\midrule
   $\XX^{\sF}_{\cF_{1,+}}$  & $A\supset(\bigwedge  ^2 V_{01}\wedge  V_{23})$ and 
   $\dim A\cap( \bigwedge  ^2 V_{01}\wedge  V_{45}\oplus V_{01}\wedge  
   \bigwedge  ^2 V_{23})\ge 1$ \\ 
\bottomrule 
\end{tabular}
\end{table} 
Let
\begin{equation*}
\XX^{*}_{\cX,+}:=\bigcup_{\sF}\BB^{\sF}_{\cX,+},\qquad
 \XX_{\cX,+}:=\ov{\XX}^{*}_{\cX,+}.
\end{equation*}
($\sF$ runs through the set of bases of $V$.)
\begin{dfn}
Let $U$ be a complex vector-space of dimension $4$ and choose an isomorphism $V\cong \bigwedge^2 U$.  Let $i_{+}\colon\PP(U)\hra \Gr(3,V)$ be the map given by~\eqref{piumenomap}. 
\begin{enumerate}
\item[(a)]
Let $\XX_{\cW}\subset\lagr$ be the closure of the set of $PGL(V)$-translates of those $A$ such that $\PP(A)$ contains  $i_{+}(Z)$ where $Z$ is a smooth quadric.
\item[(b)]
Let $\XX_{\cY}\subset\lagr$ be the set of $PGL(V)$-translates of those $A$ such that $\PP(A)$ contains  $i_{+}(Z)$ where $Z$ is either a quadric cone or a plane.
\end{enumerate}
\end{dfn}
Lastly let $L$ be a complex vector-space of dimension $3$. Choose an isomorphism $V\cong \Sym^2 L$ and let $A_k(L),A_h(L)\in\lagr$ be the lagrangians defined by~\eqref{aconkconh}. We let
\begin{equation*}
\XX_{k}:=\ov{PGL(V)A_k(L)},\qquad \XX_{h}:=\ov{PGL(V)A_h(L)}.
\end{equation*}
\begin{rmk}\label{rmk:maggiorato}
If $\cX$ appears in Table~\eqref{comericci} then $\XX^{*}_{\cX,+}\subset\BB^{*}_{\cX}$. Moreover  $\XX_{\cY},\XX_{\cW}\subset(\BB_{\bf R}\cap\BB_{\bf S})$, $\XX_{k}\subset \BB_{\bf T}$ and $\XX_{h}\subset \BB_{{\bf T}^{\vee}}$. 
\end{rmk}
\begin{thm}\label{thm:tetadimdue}
Suppose that $A\in\lagr$ and  that $\Theta_A$ contains a $2$-dimensional irreducible component. 
Then 
\begin{equation}\label{maggiolino}
A\in (\XX_{\cA+}\cup \XX_{\cA^{\vee}_{+}}\cup \XX_{\cC_{1,+}}\cup \XX_{\cC_{2,+}}\cup \XX_{\cD_{+}}\cup  \XX_{\cE_{2,+}}\cup \XX_{\cE^{\vee}_{2,+}}\cup \XX_{\cF_{1,+}}\cup \XX_{\cY}\cup \XX_{\cW}\cup \XX_{h}\cup \XX_{k}).
\end{equation}
\end{thm}
In a forthcoming paper we will prove that each of the $\XX_{\cX_{+}}$ appearing in~\eqref{maggiolino}  is contained in the GIT-unstable locus of $\lagr$ and hence is irrelevant when considering moduli of (double)EPW-sextics. That is the reason why the statement of~\Ref{thm}{tetadimdue} is not as detailed as that of~\Ref{thm}{tetadimuno}.  In fact the set in~\eqref{maggiolino} contains strictly the locus of $A\in\lagr$ such that $\Theta_A$ contains an irreducible component of dimension $2$:
 if  $A$ is generic in  $ \XX_{\cC_{1,+}}$  then $\dim\Theta_A=0$. The proof of~\Ref{thm}{tetadimdue} will be given at the end of the present subsection:   we will first prove a series of auxiliary results. 
For the rest of this subsection  $V_0\subset V$ will be a $5$-dimensional subspace.  
\begin{prp}\label{prp:cocaina}
Keeping notation as above let $L\subset\PP(\bigwedge^2 V_0)$ be a linear subspace not contained in $\Gr(2,V_0)$. Suppose that $L\cap \Gr(2,V_0)$ contains a  hypersurface which is not  a hyperplane. Then there exists $U\in\Gr(4,V_0)$ such that  
$L\subset \PP(\bigwedge^2 U)$.
\end{prp}
\begin{proof}
Let $Z\subset (L\cap \Gr(2,V_0))$ be the hypersurface which exists by hypothesis. Then $L\cap\Gr(2,V_0)$ is an intersection of quadrics (in $L$) because 
$\Gr(2,V_0)$ is the intersection of Pl\"ucker quadrics in $\PP(\bigwedge^2 V_0)$. Since $L\cap\Gr(2,V_0)$ contains the non-degenerate hypersurface $Z$  we get that $Z$ is a quadric and it equals the scheme-theoretic intersection  $L\cap\Gr(2,V_0)$.
Let us consider  the rational map
\begin{equation}
f\colon\PP(\bigwedge^2 V_0)\dashrightarrow |I_{\Gr(2,V_0)}(2)|^{\vee}=\PP(V_0^{\vee}).
\end{equation}
The restriction of $f$ to  $L$ is regular and constant  because $Z$ is a quadric hypersurface in $L$. On the other hand $f$ has the following geometric interpretation by~\Ref{lmm}{geomappa}: if $[\alpha]\in(\PP(\bigwedge^2 V_0)\setminus \Gr(2,V_0))$ then $f([\alpha])$ is canonically equal to the span of $\alpha$. Since $f$ is constant on $L$  the proposition follows.
\end{proof}
Given a quadric $Q_0\subset\PP(V_0)$ let
\begin{equation}
F_1(Q_0):=\{\ell\subset Q_0\mid \text{$\ell$ a line}\}
\end{equation}
be the variety parametrizing lines contained in $Q_0$.
We will need an explicit description of $F_1(Q_0)$ for $Q_0$ of  corank at most $1$. We start by recalling how one describes $F_1(Q)$ for $Q\subset\PP(V)$ a smooth quadric. Let $U$ be a $4$-dimensional complex vector-space; choose an isomorphism $\PP(V)\cong \PP(\bigwedge^2 U)$ taking  $Q$ to $\Gr(2,U)$ (embedded in $\PP(\bigwedge^2 U)$ by Pl\"ucker). Let $Z\subset\PP(U)\times\PP(U^{\vee})$ be the incidence subvariety of couples $([u],[f])$ such that $f(u)=0$. We have an isomorphism
\begin{equation}\label{riorio}
\begin{matrix}
Z & \overset{\mu}{\lra} & F_1(Q) \\
([u],[f]) & \mapsto & \{K\in\Gr(2,U)\mid u\in K\subset \ker(f)\}
\end{matrix}
\end{equation}
Furthermore
\begin{equation}
\mu^{*}\cO_{F_1(Q)}(1)\cong\cO_{\PP(U)}(1)\boxtimes
\cO_{\PP(U^{\vee})}(1)
\end{equation}
Let 
\begin{equation}
\PP(U)\overset{\pi_1}{\longleftarrow} 
Z\overset{\pi_2}{\lra}\PP(U^{\vee}) 
\end{equation}
be the projections - via Isomorphism~(\ref{riorio}) we will also view $\pi_1,\pi_2$ as maps with domain $F_1(Q)$.  
  Now we are ready to describe $F_1(Q_0)$. We have 
\begin{equation}\label{winnypooh}
Q_0=Q\cap V(\sigma),\qquad \sigma\in\bigwedge^2 V^{\vee}
\end{equation}
where
\begin{itemize}
\item[(a)]
 $\sigma$ is non-degenerate if $Q_0$ is smooth,
\item[(b)]
 $\dim\ker(\sigma)=2$ if $\cork Q_0=1$.
\end{itemize}
 The proof of the following result is an easy exercise that we leave to the reader.
\begin{clm}\label{clm:fanquliscia}
Keep notation as above and suppose that $Q_0$ is smooth. The restrictions of $\pi_1$ and $\pi_2$ to $F_1(Q_0)$ are isomorphisms onto $\PP(U)$ and $\PP(U^{\vee})$ respectively; let $\nu_1\colon\PP(U)\overset{\sim}{\lra} F_1(Q_0)$ and $\nu_2\colon\PP(U^{\vee})\overset{\sim}{\lra} F_1(Q_0)$ be the inverses. Then
\begin{equation}
\nu_1^{*}\cO_{F_1(Q_0)}(1)\cong\cO_{\PP(U)}(2),\qquad
\nu_2^{*}\cO_{F_1(Q_0)}(1)\cong\cO_{\PP(U^{\vee})}(2).
\end{equation}
\end{clm}
Now suppose that $\cork Q_0=1$ i.e.~Item~(b) holds. Let
\begin{equation}
L_1:=\PP(\ker\sigma)\subset\PP(U),\qquad
L_2:=\PP(\Ann  (\ker\sigma))\subset\PP(U^{\vee}),
\end{equation}
where $\sigma$ is as in~(\ref{winnypooh}). Thus  $L_1,L_2$ are lines.  Let $T_1,T_2\subset Z$ be the closed subsets defined by
\begin{equation}
T_i  =\pi_i^{-1}L_i.
\end{equation}
We leave the easy proof of the following result to the reader.
\begin{clm}\label{clm:fanqusing}
Keep notation as above and suppose that $\cork Q_0=1$. Then:
\begin{itemize}
\item[(1)]
The irreducible components of $F_1(Q_0)$ are $\mu(T_1)$ and $\mu(T_2)$.
\item[(2)]
 Let $\rho_i:=\pi_{3-i}|_{T_i}$. The map $\rho_i$ is the blow-up of $L_{3-i}$.
\item[(3)]
Let $E_i\subset T_i$ be the exceptional divisor of $\rho_i$; then
\begin{eqnarray}
\cO_{T_1}(1) \cong & \rho_1^{*}\cO_{\PP(U^{\vee})}(2)(-E_1), \label{corto}\\
\cO_{T_2}(1)  \cong & \rho_2^{*}\cO_{\PP(U)}(2)(-E_2)
\label{maltese}
\end{eqnarray}
(We view $T_i$ as a subset of $F_1(Q_0)$ via Item~(1) and we let $\cO_{T_i}(1)$ be the restriction to $T_i$ of the Pl\"ucker line-bundle on $F_1(Q_0)$.) 
\item[(4)]
The maps between spaces of global sections induced by~(\ref{corto})-(\ref{maltese}) are surjective.
\end{itemize}
\end{clm}
In the following lemma we think of $\Gr(2,V_0)$ as embedded in $\PP(\bigwedge^2 V_0)$ by Pl\"ucker: 
given $W\subset\Gr(2,V_0)$ we denote by  $\la W\ra$  the span of $W$ in $\PP(\bigwedge^2 V_0)$.
\begin{lmm}\label{lmm:mercury}
Let $Q_0\subset\PP(V_0)$ be an irreducible quadric. Suppose that $S\subset F_1(Q_0)$ is a projective surface and that 
\begin{itemize}
\item[(a)]
$S$ is a two-dimensional cubic rational normal scroll (possibly singular), or
\item[(b)]
$\dim\la S\ra=4$ and $S$ is the intersection of two quadric surfaces in $\la S\ra$.
\end{itemize}
Then there exists a line $\ell_0\subset\PP(V_0)$ which intersects all the lines  parametrized by $S$.
\end{lmm}
\begin{proof}
Since $Q_0$ is irreducible the corank of $Q_0$ is $0$, $1$ or $2$. We claim that  $Q_0$ can not be smooth. In fact suppose that $Q_0$ is smooth. By~\Ref{clm}{fanquliscia} we get that $\cO_S(1)$ is divisible by $2$; this is absurd because  a surface $S$ satisfying Item~(a) or Item~(b) contains lines. Next let us suppose that $\cork Q_0=1$. We adopt the notation of~\Ref{clm}{fanqusing}. Since $S$ is irreducible we have $S\subset \mu(T_i)$ for $i=1$ or $i=2$. By simmetry we may suppose that $S\subset\mu(T_2)$.  Since $\dim\la S\ra=4$ it follows easily from Item~(4) of~\Ref{clm}{fanqusing} that $\rho_2(S)$ is a plane not containing $L_1$ - in particular $S$ is a smooth normal cubic scroll. Let $q:=L_1\cap \rho_2(S)$; then $(q,\rho_2(S))\in T_1$ and hence  $\mu(q,\rho_2(S))$ is a line $\ell_0\subset Q_0$. By construction every line in $S$ intersects $\ell_0$. 
Finally suppose that $\cork Q_0=2$. Then $\sing Q_0$ is a line and every line contained in $Q_0$ intersetcts $\sing Q_0$; thus the lemma holds in this case as well.
\end{proof}
\begin{prp}\label{prp:supgrass}
Keep notation as above. Suppose that $\Lambda\subset\PP(\bigwedge^2 V_0)$ is a linear subspace and that
there exists  a $2$-dimensional irreducible component $S$ of $\Lambda\cap \Gr(2,V_0)$. 
 Then one of the following holds:
\begin{itemize}
\item[(1)]
$S$ is a plane.
\item[(2)]
$\dim\la S\ra=3$ and there exists $W\in\Gr(4,V_0)$  such that 
$S\subset \Gr(2,W)$.
\item[(3)]
$\dim\la S\ra=4$ and there exists a line $\ell_0\subset\PP(V_0)$ which intersects all the lines  parametrized by $S$.
\item[(4)]
$5\le\dim\la S\ra$.
\end{itemize}
\end{prp}
\begin{proof}
 Suppose that $\dim\la S\ra=2$. Then $S$ is a plane and hence
 Item~(1) holds. On the other hand if $5\le\dim\la S\ra$ holds then Item~(4) holds.
 Thus from now  we may assume that 
 \begin{equation}
 3\le \dim\la S\ra\le 4. 
\end{equation}
Suppose that $\deg S=2$; then  Item~(2) holds by~\Ref{prp}{cocaina}. Thus we may suppose that $\deg S\ge 3$.
The intersection  $\Lambda\cap\Gr(2,V_0)$ is cut out by quadrics  and $S$ is one of its irreducible components:  it follows that  
\begin{equation}\label{maxquattro}
\dim\la S\ra=4,\qquad \deg S\le 4.
\end{equation}
 Notice also that if $\deg S=4$ then necessarily $S$ is an intersection of two quadrics. Thus $S$ is one of the following:
  \begin{itemize}
\item[($\alpha$)]
A normal cubic scroll (possibly singular).
\item[($\beta$)]
An intersection  of two quadrics (in $\Lambda$) which is not a cone.
\item[($\gamma$)]
An intersection  of two quadrics (in $\Lambda$) which is a cone over a degree-$4$ curve of arithmetic genus $1$. 
\end{itemize}
If~($\gamma$) holds then Item~(3) holds with $\ell_0$ the line corresponding to the vertex of the cone. Thus we may suppose that either~($\alpha$) or~($\beta$) holds. 
Let  $T\to S$ be a desingularization of $S$; since either~($\alpha$) or~($\beta$) holds $T$ is rational. In particular we have
 \begin{equation}\label{chiuno}
\chi(\cO_{T})=1. 
\end{equation}
Composing the desingularization map $T\to S$ with the inclusion $S\hra\Gr(2,V_0)$ we get a map $g\colon T\to\Gr(2,V_0)$. Let $\cE_T$ be the pull-back to $T$
of the tautological rank-$2$ vector-bundle on $\Gr(2,V_0)$, $f_T\colon\PP(\cE_T)\to\PP(V_0)$ be the tautological map and $R_T=\im(f_T)$.
 Let $\xi:=\cO_{\PP(\cE_T)}(1)$. Let $\theta\colon\PP(\cE_T)\to T$ be the bundle map. We let $F:=\cE_T^{\vee}$; thus $F$ is globally generated.
The relation  
\begin{equation}
c_1(\xi)^2-(\theta^{*} c_1(F)) c_1(\xi)+\theta^{*} c_2(F)=0
\end{equation}
gives the equation
\begin{equation}
\int_{T} (c_1(F)^2-c_2(F))=\int_{\PP(\cE_T)}c_1(\xi)^3.
\end{equation}
Thus we get
\begin{equation}\label{tciupapa}
\int_{T} (c_1(F)^2-c_2(F))=
\begin{cases}
\deg f_T\cdot \deg R_T & \text{if $\dim R_T=3$,}\\
0 & \text{if $\dim R_T<3$,}
\end{cases}
\end{equation}
and in particular 
\begin{equation}\label{tempomele}
0\le c_2(F)\le c_1(F)^2.
\end{equation}
(The first inequality holds because $F$ is globally generated.) 
We notice that $\dim R_T=3$ unless $R_T$ is a plane.
Since~($\alpha$) or~($\beta$) holds $R_T$ is not a plane and hence $\dim R_T=3$.  By~(\ref{tciupapa})-(\ref{tempomele}) we have $\deg R_T\le 4$. We claim that $2\le\deg R_T$. In fact suppose that $\deg R_T=1$. Then $R_T=\PP(W)$ where $W\in\Gr(4,V)$ and hence $S\subset\Gr(2,W)$. Since  $\Gr(2,W)$ is a smooth quadric hypersurface in $\PP(\bigwedge^2 W)$ and $\dim\la S\ra=4$ the intersection $\la S\ra\cap\Gr(2,W)$ is a $3$-dimensional irreducible quadric containing $S$; this contradicts the hypothesis that $S$ is an irreducible component of $\Lambda\cap\Gr(2,V_0)$.  
If $\deg R_T=2$ then  $R_T$ is an irreducible $3$-dimensional quadric and hence Item~(3) holds by~\Ref{lmm}{mercury}. It remains to prove the proposition under the assumption that $\dim R_T=3$ and $3\le \deg R_T\le 4$.  Let's prove that
 \begin{equation}\label{moltend}
h^0(End F)\ge 3.
\end{equation}
 By~(\ref{tciupapa}) it follows that if $S$ is a cubic normal scroll then $c_1(F)^2=3$, $c_2(F)=0$ and if $S$ is the intersection of two quadrics then $c_1(F)^2=4$,
 $0\le c_2(F)\le 1$. 
 In both cases we have
 \begin{equation}\label{piccolocidue}
0\le c_2(F)\le 1.
\end{equation}

 Thus Hirzebruch-Riemann-Roch gives that
 \begin{equation}
\chi(End F)=4\chi(\cO_{T})-(4c_2(F)-c_1(F)^2)\ge 4\chi(\cO_{T})=4. 
\end{equation}
(The last equality holds by~(\ref{chiuno}).) By Serre duality we get 
that 
\begin{equation}\label{endemol}
h^0(End F)+h^0(End F(K_{T}))\ge 4.
\end{equation}
One easily checks that there exists a non-zero $\sigma\in H^0(-K_{T})>0$. Thus $\sigma$ defines a map  $F(K_{T})\to F$ which is an isomorphism away from the (non-empty) zero-set of $\sigma$; it follows that we have an inclusion $H^0(End F(K_{T}))\subset H^0(End F)$ which does not contain the subsapce of homotheties.
Hence $h^0(End F)\ge 1+h^0(End F(K_{T}))$ and thus~\eqref{moltend} follows from~\eqref{endemol}. 
\begin{clm}\label{clm:teletubbies}
There exist divisors $C,D$ on $T$ and a zero-dimensional subscheme $Z\subset T$ such that $F$ fits into an exact sequence
\begin{equation}
0\lra \cO_{T}(C)\lra F\lra I_Z(D)\lra 0.
\end{equation}
and the following hold:
\begin{itemize}
\item[(I)]
$\cO_{T}(C)$ and $I_Z(D)$ are globally generated, in particular
\begin{equation}\label{positivo}
C\cdot C\ge 0,\qquad C\cdot D\ge 0,\qquad D\cdot D\ge 0.
\end{equation}
\item[(II)]
$h^0(\cO_{T}(C-D))>0$ but $\cO_{T}(C-D)\not\cong\cO_{T}$.
\end{itemize}
\end{clm}
\n
{\it Proof of the claim.}
By~(\ref{moltend}) there exists $\psi\in H^0(End F)$ which is not a scalar. Let $\lambda$ be an eigen-value of $\psi$ (notice that the characteristic polynomial of $\psi$ has constant coefficients) and $\phi:=(\psi-\lambda \Id _F)$. Let $\cK:=\ker\phi$; then $\cK$ is a rank-$1$ subsheaf of $F$ and hence $\cK\cong\cO_{T}(C')$ for a certain divisor (class) $C'$. The  quotient $F/\cK$ is locally-free away from a finite set and hence it is isomorphic to $I_Z(D')$ for a certain divisor (class) $D'$. Thus we have an exact sequence
\begin{equation}\label{befana}
0\lra \cO_{T}(C')\lra F\overset{\pi}{\lra} I_Z(D')\lra 0
\end{equation}
and an inclusion of vector-bundles $\iota\colon\cO_{T}(D')\hra F$ (injective on fibers) such that $\phi=\iota\circ\pi$. Suppose that $\pi\circ\iota\circ\pi\not=0$; then~(\ref{befana}) splits and hence $F\cong  \cO_{T}(C')\oplus \cO_{T}(D')$. By~(\ref{moltend}) we get that $h^0(\cO_{T}(C'-D')>0$ or $h^0(\cO_{T}(D'-C')>0$; if the former holds we let $C:=C'$ and $D:=D'$, if the latter holds we let
 $C:=D'$ and $D:=C'$. With these choices the claim holds except possibly for the assertion that $\cO_{T}(C-D)\not\cong\cO_{T}$. We have
 \begin{equation}\label{aussie}
\cO_{T}(C+D)\cong g^{*}\cO_S(1).
\end{equation}
Suppose  that $\cO_{T}(C-D)\cong\cO_{T}$; then $g^{*}\cO_S(1)$ is divisible by two, this is absurd because $S$ contains lines which are not contained in $\sing S$. This proves the claim under the assumption that $\pi\circ\iota\circ\pi\not=0$. Now suppose that $\pi\circ\iota\circ\pi=0$. Then $\phi$ defines a non-zero map $I_Z(D')\to\cO_{T}(C')$ and hence we get that $h^0(\cO_{T}(C'-D')>0$. Let
 $C:=C'$ and $D:=D'$. 
 The claim holds with these choices.
\qed
\vskip 2mm
We resume the proof of~\Ref{prp}{supgrass}. Let $C,D$ and $Z$ be as in~\Ref{clm}{teletubbies}; we have~(\ref{aussie}) and by Whitney's formula we get that
\begin{eqnarray}
C\cdot D+\ell(Z)  = & c_2(F), \label{ciddue}\\
C\cdot C+2 C\cdot D +D\cdot D  = & \deg S. \label{gradisca}
\end{eqnarray}
By~(\ref{piccolocidue}), (\ref{positivo}) and~(\ref{ciddue}) we have $0\le C\cdot D\le 1$. 
Suppose that $C\cdot D=0$. By~\Ref{clm}{teletubbies} and Hodge Index  we get that $D=0$ and $F\cong\cO_{T}(C)\oplus\cO_{T}$. Thus $T$ is a cone and hence Item~(3) holds with $\ell_0$ the line corresponding to the vertex of $T$.  Now suppose that $C\cdot D=1$. By~(\ref{piccolocidue}) and~(\ref{ciddue}) we have $Z=\emptyset$. By~(\ref{aussie}),
Item~(II) of~\Ref{clm}{teletubbies} and Hodge Index we have 
\begin{equation}
0<(C-D)\cdot (C+D)=C\cdot C-D\cdot D.
\end{equation}
By~(\ref{positivo})  and~(\ref{gradisca}) we get that one of the following holds:
\begin{itemize}
\item[(i)]
$C\cdot C=1$ and $D\cdot D=0$ ($\deg S=3$).
\item[(ii)]
$C\cdot C=2$ and $D\cdot D=0$ ($\deg S=4$).
\end{itemize}
Since $\cO_{T}(D)$ is globally generated and $C\cdot D=1$ it follows that $h^0(\cO_{T}(D))=2$. Furthermore we get  surjectivity of  the map $\epsilon\colon V_0^{\vee}\to  H^0(\cO_{T}(D))$ given by the composition $V_0^{\vee}\to H^0(F)\to H^0(\cO_{T}(D))$. Thus $\cod(\ker\epsilon,V_0^{\vee})=2$. Let $\ell_0:=\PP(\Ann  (\ker\epsilon))$; then Item~(3) holds with this choice of $\ell_0$.     
\end{proof}
Below is the  analogue of~\Ref{prp}{supgrass} obtained upon replacing $\Gr(2,V_0)$  by $\PP^2\times\PP^2$. Let $W_1,W_2$ be $3$-dimensional complex vector spaces; then $\PP(W_1)\times\PP(W_2)\subset\PP(W_1\otimes W_2)$ via  Segre's embedding. 
\begin{prp}\label{prp:supsegre}
Keep notation as above. Suppose that $\Lambda\subset\PP(W_1\otimes W_2)$ is a linear subspace such that there exists a $2$-dimensional irreducible component $S$ of
$\Lambda\cap (\PP(W_1)\times\PP(W_2))$.
 Then one of the following holds possibly after exchanging $W_1$ with $W_2$:
\begin{itemize}
\item[(1)]
$S=\{[v_0]\}\times\PP(W_2)$ for some $[v_0]\in\PP(W_1)$.
\item[(2)]
$S=\PP(U_1)\times\PP(U_2)$ where $U_i\in\Gr(2,W_i)$ i.e.~$S$ is a smooth quadric surface.
\item[(3)]
$S$ is a smooth hyperplane section of  $\PP(U_1)\times \PP(W_2)$ where $U_1\in\Gr(2,W_1)$, i.e.~$S$ is a 
smooth normal cubic scroll.
\item[(4)]
$S$ is the graph of an isomorphism $\PP(W_1)\overset{\sim}{\lra}\PP(W_2)$ (and hence is a Veronese surface).
\item[(5)]
$6\le\dim\la S\ra$.
\end{itemize}
\end{prp}
\begin{proof}
Let $f_i\colon S\to\PP(W_i)$  be the restriction of projection for $i=1,2$. For $i=1,2$ let $C_i$ be a divisor  on $S$ such that $\cO_S(C_i)\cong f_i^{*}\cO_{\PP(W_i)}(1)$. Then $\cO_S(C_1+C_2)\cong\cO_S(1)$ and hence we have
\begin{equation}\label{robzampa}
C_1\cdot C_1+2 C_1\cdot C_2+ C_2\cdot C_2=\deg S.
\end{equation}
Since $\cO_S(C_i)$ is globally generated we have
\begin{equation}\label{maggiore}
C_i\cdot C_j\ge 0,\qquad 1\le i,j\le 2.
\end{equation}
Suppose that $C_1\cdot C_2=0$; applying Hodge index to a desingularization of $S$ we get that one of $C_1,C_2$ is linearly equivalent to $0$ and hence Item~(1) holds (possibly after exchanging $W_1$ with $W_2$). Thus we may assume that
\begin{equation}
C_1\cdot C_2>0,\qquad \deg S\ge 2.
\end{equation}
 Suppose that $\deg S=2$. Then $C_1\cdot C_2=1$ and $C_i\cdot C_i=0$ for $i=1,2$; it follows easily that Item~(2) holds. Thus we are left with the case $\deg S\ge 3$. Since  $\Lambda\cap(\PP(W_1)\times\PP(W_2))$ is cut out by quadrics it follows that $4\le \dim\la S\ra $. On the other hand if $6\le \dim\la S\ra$ then Item~(5) holds and hence from now on we may assume that
 \begin{equation}\label{abbasta}
4\le\dim\la S\ra \le 5
\end{equation}
 We claim that
 \begin{equation}\label{maxcinque}
\deg S\le 5.
\end{equation}
In fact by~(\ref{abbasta}) there exists a $6$-dimensional linear space $\wt{\Lambda}\subset\PP(W_1\otimes W_2)$ containing $\Lambda$. If $\wt{\Lambda}$ is a generic such linear space then  $\wt{\Lambda}\cap(\PP(W_1)\times \PP(W_2))$ is of pure dimension $2$ and it contains an irreducible component other than $S$; thus 
\begin{equation}
6=\deg\PP(W_1)\times \PP(W_2)=
\deg (\wt{\Lambda}\cap(\PP(W_1)\times \PP(W_2)))>\deg S.
\end{equation}
This proves~(\ref{maxcinque}). By~(\ref{robzampa}) and~(\ref{maggiore}) we get that  one of the following holds possibly after exchanging $W_1$ with $W_2$:
\begin{itemize}
\item[($\alpha_k$)]
$C_1\cdot C_1=0$, $C_1\cdot C_2=1$  and $C_2\cdot C_2=k$ where $1\le k\le 3$.
\item[($\beta_m$)]
$C_1\cdot C_1=1$, $C_1\cdot C_2=1$ and $C_2\cdot C_2=m$ where $1\le m\le 2$.
\item[($\gamma_n$)]
$C_1\cdot C_1=0$, $C_1\cdot C_2=2$ and $C_2\cdot C_2=n$  where $0\le n\le 1$.
\end{itemize}
Suppose that~($\alpha_k$) holds. Then $f_1(S)$ is a curve in $\PP(W_1)$, in fact it is a line $L$ because $C_1\cdot C_2=1$. Thus $S\subset L\times\PP(W_2)$; since $C_1\cdot C_2=1$ and $C_2\cdot C_2=k$ we have $S\in|\cO_L(k)\boxtimes\cO_{\PP(W_2)}(1)|$. We claim that $\dim\la S\ra=4$. In fact suppose the contrary; then $\dim\la S\ra=5$ by~(\ref{abbasta}) and hence $\la S\ra= L\times\PP(W_2)$, contradicting the hypothesis that $S$ is an irreducible component of $\Lambda\cap\PP(W_1)\times\PP(W_2)$. Since $\dim\la S\ra=4$ we must have $k=1$ and hence Item~(2) holds.  
Next suppose that~($\beta_m$) holds: it follows  that $m=1$ and that Item~(4) holds. Lastly suppose that~($\gamma_n$) holds; we will reach a contradiction. If $n=0$ then $f_i(S)$ is a curve for $i=1,2$; since $C_1\cdot C_2=2$ it follows that $S=L\times Z$ where $L\subset\PP(W_1)$ is a line and $Z\subset\PP(W_2)$ is a smooth conic, possibly after 
 exchanging $W_1$ with $W_2$. Then $\la S\ra= L\times\PP(W_2)$, this contradicts the hypothesis that $S$ is an irreducible component of $\Lambda\cap\PP(W_1)\times\PP(W_2)$. If $n=1$ then $f_1(S)$ is a curve and since $C_1\cdot C_2=2$ we get that $f_1(S)$ is either  a line or a smooth conic. Suppose that $f_1(S)$ is a line $L$: arguing as in Case~($\alpha_k$) for $k=2,3$ one gets that $\la S\ra\supset L\times\PP(W_2)$, this contradicts the hypothesis that $S$ is an irreducible component of $\Lambda\cap\PP(W_1)\times\PP(W_2)$.
Lastly suppose that $f_1(S)$ is a (smooth) conic. Then one gets that  
$f_2\colon S\to\PP(W_2)$ is the blow-up of a point and that the linear system cut out on $S$ by $|\cO_{\PP(W_1\otimes W_2)}(1)|$ is equal to $|f_2^{*}\cO_{\PP(W_2)}(3)(-2E)|$ where $E\subset S$ is the exceptional divisor of $f_2$. It follows that $\dim\la S\ra=6$, that contradicts~(\ref{abbasta}).  
\end{proof}
The following is our last preliminary result. 
\begin{clm}\label{clm:veronsegre}
Let $U\in\Gr(3,V)$. Suppose that  $\Theta\subset I_U$ is a projective surface such that
\begin{itemize}
\item[(a)]
$U\notin \Theta$,
\item[(b)]
$\ov{\rho}_U(\Theta)$ is the graph of an isomorphism $\PP(\bigwedge^2 U)\overset{\sim}{\lra}\PP(V/U)$. (See~\eqref{robaruma} for the definition of $\ov{\rho}_U$.)
\end{itemize}
Then there exist an identification $V=\bigwedge^2 \CC^4$ and   a plane $Z\subset\PP(\CC^4)$  such that 
$\Theta=  i_{+}(Z)$ where $i_{+}$ is given by~(\ref{piumenomap}).
\end{clm}
\begin{proof}
Let $Z\subset\PP^3$ be  a plane and  $\Theta=  i_{+}(Z)$.   Let $U\subset V$ be such that $\PP(U)=i_{-}(Z)$ (we recall that $i_{-}\colon (\PP^3)^{\vee}\hra \Gr(3,V)$); then Items~(a) and~(b) hold. The result follows because
 $SL(V)$ acts transitively on the family of couples $(U,\Theta)$ such that Items~(a) and~(b) hold. 
\end{proof}

\noindent
{\it Proof of~\Ref{thm}{tetadimdue}.\/} 
 By Morin one of~(a) - (e) of~\Ref{thm}{teomorin} holds. We will perform a case-by-case analysis. 
\vskip 2mm
\n
 (a): $\Theta\subset F_{\pm}(\cQ)$.  
  Let $U$ be a $4$-dimensional complex vector-space and identify $V$ with $\bigwedge^2 U$ so that $\cQ$ gets identified with $\Gr(2,U)$.  We may assume that  $\Theta\subset F_{+}(\cQ)$; thus   $\Theta:=i_{+}(Z)$ for an irreducible surface $Z\subset\PP(U)$.  By~(\ref{paccheri}) the map $i_{+}$ is defined by the linear system of quadrics in $\PP(U)$. It follows that $Z$ is contained in a quadric of  $\PP(U)$ and hence is a plane or a quadric, thus $A\in(\XX_{\cY}\cup \XX_{\cW})$. 
\vskip 2mm
\n
(b): $\Theta\subset C(\cV)$ or $\Theta\subset T(\cV)$.  
Then $A\in(\XX_k\cup\XX_h)$. 
\vskip 2mm
\n
 (c): $\Theta\subset J_{v_0}$. 
 Let $V_0\in\Gr(5,V)$ 
be transversal to $[v_0]$ 
and  $\ov{\rho}_{v_0}$ be given by~(\ref{robarvu}). Let $S:=\ov{\rho}_{v_0}(\Theta)\subset\Gr(2,V_0)$ and $\Lambda:=\la S\ra$. Then the hypotheses of~\Ref{prp}{supgrass} are satisfied and hence one of Items~(1)-(4) of that proposition holds. If Item~(1)  holds then $\Theta$ is plane: as is easily checked
\begin{equation}\label{sepiano}
\text{if $\Theta_A$ contains a plane then $A\in(\XX_{\cF_{1,+}}\cup\XX_{\cF_{2,+}})$.} 
\end{equation}
If Item~(2) of~\Ref{prp}{supgrass} holds then  $A\in\XX_{\cD_{+}}$. 
Suppose that Item~(3) of~\Ref{prp}{supgrass}  holds. Let $V_{12}\subset V_0$ be the  subspace of dimension $2$ such that $\ell_0=\PP(V_{12})$  and $V_{35}\subset V_0$ be a subspace complementary to $V_{12}$. Then 
$S\subset\PP(\bigwedge^2 V_{12}\oplus V_{12}\wedge V_{35})$
and hence
\begin{equation*}
\Theta\subset\PP([v_0]\wedge\bigwedge^2 V_{12}\oplus 
[v_0]\wedge V_{12}\wedge V_{35}).
\end{equation*}
  Since $\dim\la \Theta\ra=\dim\la S\ra=4$ we get that $A\in\XX_{\cE_{2,+}}$.
If Item~(4) of~\Ref{prp}{supgrass} holds then $A\in\XX_{\cA_{+}}$. 
\vskip 2mm
\n
(d): $\Theta\subset \Gr (3,E)$ where $E\in\Gr(5,V)$. 
By Case~(c) and duality we get that 
\begin{equation*}
A\in \XX_{\cF_{1,+}}\cup\XX_{\cF_{2,+}}\cup \XX_{\cD_{+}}\cup 
\XX_{\cE^{\vee}_{2,+}}\cup\XX_{\cA^{\vee}_{+}}.
\end{equation*}
\vskip 2mm
\n
(e): $\Theta\subset I_U$ where $U\in\Gr(3,V)$. 
 Since $\dim\Theta=2$ we have $2\le\dim\la\Theta\ra$. If   $\dim\la\Theta\ra=2$ then $\Theta$ is  plane: by~\eqref{sepiano} we have
$A\in(\XX_{\cF_{1,+}}\cup\XX_{\cF_{2,+}})$.  Thus from now on we may suppose that
\begin{equation}\label{difretta}
3\le\dim\la\Theta\ra.
\end{equation}
We distinguish between the two cases:
\begin{itemize}
\item[(e1)]
$U\in\la\Theta\ra$,
\item[(e2)]
$U\notin\la\Theta\ra$.
\end{itemize}
Suppose that~(e1) holds; we will prove that
\begin{equation}\label{ahiservaitalia}
A\in(\XX_{\cC_{1,+}}\cup \XX_{\cD_{+}}).
\end{equation}
Since $U\in\la\Theta\ra$ and $\Theta$ is an irreducible component  of  $\la\Theta\ra\cap I_U$ we get that
 $\Theta$ is a cone with vertex $U$. 
If $4\le \dim\la\Theta\ra$ then $A\in\XX_{\cC_{1,+}}$. Thus we may assume that
\begin{equation}\label{cenetta}
 \dim\la\Theta\ra= 3.
\end{equation}
Let $\ov{\rho}_U$ be the map of~\eqref{robaruma}. Then $C:=\ov{\rho}_U(\Theta)$ is a $1$-dimensional irreducible component of $\la C\ra\cap(\PP(\bigwedge^2 U)\times\PP(V/U))$.
By~(\ref{cenetta}) we have $\dim\la C \ra= 2$;
thus $C$ is a smooth conic 
because $\PP(\bigwedge^2 U)\times\PP(V/U)$ is cut out by quadrics (in $\PP(\bigwedge^2 U\otimes(V/U))$). Let $f\colon C\to \PP(\bigwedge^2 U)$ and $g\colon C\to \PP(V/U)$ be the projection maps. Neither $f$ nor $g$ is constant because $C$ is a $1$-dimensional irreducible component of $\la C\ra\cap(\PP(\bigwedge^2 U)\times\PP(V/U))$; thus 
\begin{equation}
f^{*}\cO_{\PP(\bigwedge^2 U)}(1)\cong
g^{*}\cO_{\PP(V/U)}(1)\cong\cO_C(1).
\end{equation}
It follows that $\im(f)$ and $\im(g)$ are lines and hence there exist $[v_0]\in\PP(U)$ and $W_2\in\Gr(2,V/U)$ such that 
\begin{equation}
\im(f)=\PP(\{v_0\wedge u\mid u\in U\}),\qquad
\im(g)=\PP(W_2).
\end{equation}
Let $U_0\subset U$ be complementary to $[v_0]$. Let 
$V_{14}\subset V$ be the $4$-dimensional subspace containing $U_0$ and
projecting to $W_2$ under the quotient map $V\to V/U$.
Let $\sF:=\{v_0,v_1,\ldots,v_5\}$ be a basis of $V$ adapted to $V_{14}$ i.e.~$V_{14}=\la v_1,\ldots,v_4\ra$:  then  $A\in\XX^{\sF}_{\cD_{+}}$. This finishes the proof that if Item~(e1) above holds then~\eqref{ahiservaitalia} holds. Next suppose that Item~(e2) holds.
 Let $W\subset V$ be a  subspace complementary to $U$; thus 
 \begin{equation}
\la\Theta\ra\subset \PP(\bigwedge^3 U\oplus \bigwedge^2 U\wedge W).
\end{equation}
Let $W_1:=\bigwedge^2 U$, $W_2:=W$, $S:=\ov{\rho}_U(\Theta)\subset \PP(W_1)\times\PP(W_2)$ and $\Lambda:=\la S\ra$; 
then the hypotheses of~\Ref{prp}{supsegre} are satisfied. Moreover since Item~(e2) holds we have $\dim\Lambda=\dim\la\Theta\ra$. By~\eqref{difretta} one of Items~(2)-(5) of~\Ref{prp}{supsegre} holds.   If Item~(2)  holds then $A\in\XX_{\cD_{+}}$.  If Item~(3) of~\Ref{prp}{supsegre} holds then $A\in\XX_{\cE_{2,+}}$. 
 If Item~(4) of~\Ref{prp}{supsegre} holds then  $A\in\XX_{\cY}$ by~\Ref{clm}{veronsegre}. If Item~(5) of~\Ref{prp}{supsegre} holds then $A\in_{\cC_{2,+}}$.
\qed
\begin{prp}\label{prp:anchedue}
Let $A\in\lagr$ and suppose that $\Theta_A$ contains an irreducible component of dimension $2$. Then there exists a Type $X$ appearing in Table~\eqref{catundim} such that $A\in\BB_X$.  
\end{prp}
\begin{proof}
By~\Ref{thm}{tetadimdue} we know that $A$ belongs to the right-hand side of~\eqref{maggiolino}. By~\Ref{rmk}{maggiorato} there exists a Type $X$ in Table~\eqref{catundim} such that $A\in\BB_X$ except possibly if $A\in\BB_{\cC_{1,+}}$. It is not true that $\BB_{\cC_{1,+}}$ is contained in on of the $\BB_X$ - in fact $\Theta_A$ is a singleton for generic $A\in  \BB_{\cC_{1,+}}$. However going through the proof of~\Ref{thm}{tetadimdue}   we see that the only instance in which $A\in  \BB_{\cC_{1,+}}$ corresponds to Item~(e1), see~\eqref{ahiservaitalia}. Let $A_0$ be as in Item~(e1) with $A_0\in\BB_{\cC_{1,+}}$; thus $\Theta_A$ contains a $2$-dimensional irreducible component $\Theta\subset I_U$, where $U\in\Gr(3,V)$, $\Theta$ is a cone with vertex $U$ and
$4\le \dim\la\Theta\ra$.  Let $H\subset\la\Theta\ra$ be a generic codimension-$1$ linear subspace, in particular it does not contain $U$ and $\dim(H\cap\Theta)=1$. Then $A_0$ is in  the family $\cF$ of  $A\in\lagr$ containing $\la\la H\ra\ra$; if $A\in\cF$ is generic then $H\cap \Theta$ is a $1$-dimensional irreducible component of  $\Theta_A$ and hence $A$ belongs to right-hand side of~\eqref{maggiolino}. Since the $\BB_X$ are closed we get that $A_0$ is contained in one of the $\BB_X$. Actually going through the proof of\Ref{thm}{tetadimuno} we get that
$A_0\in \BB_{\cE_{2}}\cup \BB_{\cE^{\vee}_{2}}\cup\BB_{\bf Q}\cup\BB_{\cC_2}$. 
\end{proof}
\subsection{Higher-dimensional components of $\Theta_A$}
\setcounter{equation}{0}
We let
\begin{equation*}
\XX_{+}:=\ov{PGL(V)A_{+}(U)}
\end{equation*}
 where $A_{+}(U)$ is given by~\eqref{aconiconj}.   
\begin{thm}\label{thm:tetadimtre}
Suppose that $A\in\lagr$ and  that $\dim\Theta_A\ge 3$. then
\begin{equation}\label{emozione}
A\in(\XX_{\cA_{+}}\cup\XX_{\cA^{\vee}_{+}}\cup \XX_{\cC_{1,+}}\cup \XX_{\cC_{2,+}}\cup \XX_{\cD}\cup  \XX_{\cF_{1,+}}\cup \XX_{+}).
\end{equation}
\end{thm}
\begin{proof}
Since $3\le\dim\Theta$ we have $3\le\dim\la\Theta\ra$. If  the latter is an equality then $\Theta$ is a $3$-dimensional projective  space  and thus $A\in\XX_{\cF_{1,+}}$. Thus from now on we may suppose that 
\begin{equation}\label{almamater}
4\le\dim\la\Theta\ra.
\end{equation}
By Morin one of~(a), (c), (d), (e) of~\Ref{thm}{teomorin} holds (notice that $\dim C(\cV)=\dim T(\cV)=2$). We will perform a case-by-case analysis.
\vskip 2mm
\n
(a): $\Theta\subset F_{\pm}(\cQ)$. 
In this case $A\in\XX_{+}$.
\vskip 2mm
\n
(c): $\Theta\subset J_{v_0}$. 
  If $5\le\dim\la\Theta\ra$ then $A\in\XX_{\cA_{+}}$.  Thus we may assume that $\dim\la\Theta\ra=4$.
Let $V_0\in\Gr(5,V)$ be transversal to $[v_0]$ and $\ov{\rho}_{v_0}$ be as in~(\ref{robarvu}). By~(\ref{almamater}) we have $\dim\la\ov{\rho}_{v_0}(\Theta)\ra= 4$. 
 Since $\Gr(2,V_0)$ contains no projective space of dimension $4$ we get that $\ov{\rho}_{v_0}(\Theta)$ is $3$-dimensional.
By~\Ref{prp}{cocaina}  there exists $U\in\Gr(4,V_0)$ such that $\ov{\rho}_{v_0}(\Theta)\subset \Gr(2,U)$; it follows that $A\in\XX_{\cD_{+}}$. 
\vskip 2mm
\n
(d): $\Theta\subset \Gr(3,E)$ where $E\in\Gr(5,V)$.
By duality and Case (d) we get that $\in(\XX^{\vee}_{\cA_{+}}\cup \XX_{\cD_{+}})$.
\vskip 2mm
\n
 (e): $\Theta\subset I_U$ where $U\in\Gr(3,V)$.
 We distinguish  the two cases:
 \begin{itemize}
\item[(i)]
$U\in\la\Theta\ra$,
\item[(ii)]
$U\notin\la\Theta\ra$.
\end{itemize}
 Suppose that (i) holds. Since $\Theta$ is an irreducible component of $\la\Theta\ra\cap \Gr(3,V)$ we get that  $U\in\Theta$ and moreover $\Theta$ is a cone with vertex $U$. Since we are under the assumption that~\eqref{almamater} holds it follows that  $A\in\XX_{\cC_{1,+}}$.  Next suppose that~(ii) holds. 
There exists a subspace $W\subset V$ complementary to $U$ and such that 
 \begin{equation}\label{arsenal}
\la\Theta\ra\subset \PP(\bigwedge^3 U\oplus \bigwedge^2 U\wedge W). 
\end{equation}
Moreover since $U\notin\la\Theta\ra$ we have $\dim\la\ov{\rho}_U(\Theta)\ra=\dim\la\Theta\ra$ where  $\ov{\rho}_U$ be as in~\eqref{robaruma}.
If $6\le\dim\la\Theta\ra$ we get that $A\in\XX_{\cC_{2,+}}$ by~\eqref{arsenal}.
Suppose that  $\dim\la\Theta\ra\le 5$.
Then $\ov{\rho}_U(\Theta)$  is an effective Cartier divisor on $PP(\bigwedge^2 U)\times\PP(W)$ contained in  3 linearly independent   divisors of the linear system $ |\cO_{\PP(\bigwedge^2 U)}(1)\boxtimes\cO_{\PP(W)}(1) |$.
 It follows that
 \begin{itemize}
\item[($\alpha$)]
 $\ov{\rho}_U(\Theta)\in |\cO_{\PP(\bigwedge^2 U)}(1)\boxtimes\cO_{\PP(W)} |$, or
\item[($\beta$)]
$\ov{\rho}_U(\Theta)\in |\cO_{\PP(\bigwedge^2 U)}\boxtimes\cO_{\PP(W)}(1) |$.
\end{itemize}
Notice that in both cases $\Theta$ is isomorphic to $\PP^1\times\PP^2$ embedded by Segre and hence  
\begin{equation}\label{carogna}
5=\dim\la\ov{\rho}_U(\Theta)\ra=\dim\la\Theta\ra. 
\end{equation}
Suppose that~($\alpha$) holds. We have a canonical identification $\PP(U)=|\cO_{\PP(\bigwedge^2 U)}(1)|$; 
let $[v_0]\in \PP(U)$ be the point giving the divisor in $|\cO_{\PP(\bigwedge^2 U)}(1)|$ corresponding to $\Theta$. Then $\Theta\subset J_{v_0}$; since~(\ref{carogna}) holds we get that $A\in\XX_{\cA_{+}}$. Similarly one shows that if~($\beta$) holds then $A\in\XX_{\cA^{\vee}_{+}}$.   
\end{proof}
\begin{thm}\label{thm:sigminf}
The irredundant irreducible decomposition of $\Sigma_{\infty}$  is given by
\begin{equation}\label{decosigma}
\Sigma_{\infty}=\BB_{\cA}\cup\BB_{\cA^{\vee}}\cup\BB_{\cC_2}\cup
\BB_{\cD}\cup\BB_{\cE_2}\cup\BB_{\cE_2^{\vee}}\cup
\BB_{\cF_1}\cup\BB_{\bf Q}\cup\BB_{\bf R}\cup\BB_{\bf S}\cup\BB_{\bf T}\cup\BB_{{\bf T}^{\vee}}
\end{equation}
\end{thm}
\begin{proof}
By~\Ref{crl}{sonocomp} all we have to prove is that $\Sigma_{\infty}$ is contained in the right-hand side of~\eqref{decosigma}. Let $A\in\Sigma_{\infty}$; thus $\dim\Theta_{A}\ge 1$. If $\dim\Theta_{A}= 1$ then $A$ belongs to the right-hand side of~\eqref{decosigma} by~\Ref{thm}{tetadimuno}.  If $\dim\Theta_{A}= 2$ then $A$ belongs to the right-hand side of~\eqref{decosigma} by~\Ref{thm}{tetadimdue} and~\Ref{prp}{anchedue}. Lastly suppose that $\dim\Theta_{A}\ge 3$; then $A$ belongs to the right-hand side of~\eqref{emozione} by~\Ref{thm}{tetadimtre}. By~\Ref{rmk}{maggiorato} we get that $A$ belongs to the right-hand side of~\eqref{decosigma} except possibly if $A\in\BB_{\cC_{1,+}}$; arguing as in the proof of~\Ref{prp}{anchedue} we get that $A$ belongs to the right-hand side of~\eqref{decosigma} in that case as well.
\end{proof}
\section{Appendix: Quadratic forms}
 \setcounter{equation}{0}
Let $U$ be a complex vector-space of finite dimension  $d$. We view $\Sym^2 U^{\vee}$ as the vector-space of quadratic forms on $U$. 
 We will recall a few standard results regarding the loci
\begin{equation}
(\Sym^2 U^{\vee})_r:=\{q\in \Sym^2 U^{\vee}\mid\rk q\le r\},
\qquad 0\le r\le d.
\end{equation}
The following is   well-known.
\begin{prp}\label{prp:tanquad}
Keep notation as above. Then $(\Sym^2 U^{\vee})_r$ is a closed irreducible set of codimension ${d-r+1}\choose{2}$ in $\Sym^2 U^{\vee}$, smooth away from  
$(\Sym^2 U^{\vee})_{r-1}$. Let $q_{*} \in((\Sym^2 U^{\vee})_r\setminus (\Sym^2 U^{\vee})_{r-1})$ and $K:=\ker q_{*}$; the (embedded) tangent space to $(\Sym^2 U^{\vee})_r$ at $q_{*}$ is
\begin{equation}
T_{q_{*}}(\Sym^2 U^{\vee})_r=\{q\in \Sym^2 U^{\vee}\mid q|_K=0\}\,.
\end{equation}
\end{prp} 
By~\Ref{prp}{tanquad}  the normal cone of 
$(\Sym^2 U^{\vee})_r$ in $\Sym^2 U^{\vee}$ is a vector-bundle away from $(\Sym^2 U^{\vee})_{r-1}$ and if $q_{*}\in((\Sym^2 U^{\vee})_r\setminus (\Sym^2 U^{\vee})_{r-1})$ 
  then we have a canonical isomorphism
\begin{equation}\label{norquad}
\begin{matrix}
(C_{(\Sym^2 U^{\vee})_r}\Sym^2 U^{\vee})_{q_{*}} & \overset{\sim}{\lra} &
\Sym^2 K^{\vee} \\
[q] & \mapsto & q|_K.
\end{matrix} 
\end{equation}
Next we consider $(\Sym^2 U^{\vee})_{d-1}$ i.e.~the locus defined by vanishing of the determinant of $q$. 
Given $q_{*}\in \Sym^2 U^{\vee}$
we let $\Phi$ be the polynomial on the vector-space $\Sym^2 U^{\vee}$ defined by $\Phi(q):=\det(q_{*}+q)$. Of course $\Phi$  is  defined up to  multiplication  by a non-zero scalar, moreover it depends on $q_{*}$ although that does not show up in the notation. Let 
\begin{equation}\label{granfi}
\Phi=\Phi_0+\Phi_1+\ldots +\Phi_d,
\qquad \Phi_i\in \Sym^i (\Sym^2 U)
\end{equation}
be the decomposition into homogeneous components. 
The result below follows from a straightforward  computation.
\begin{prp}\label{prp:conodegenere}
Let $q_{*}\in\Sym^2 U^{\vee}$
 and 
\begin{equation}
K:=\ker(q_{*}), \qquad k:=\dim K. 
\end{equation}
  Let $\Phi_i$ be the polynomials appearing in~\eqref{granfi}. Then
\begin{itemize}
\item[(1)]
$\Phi_i=0$ for $i<k$, and
\item[(2)]
there exists $c\not=0$ such that $\Phi_k(q)=c\det(q|_K)$. 
\end{itemize}
\end{prp}
\end{document}